\documentclass[11pt]{amsart} 
\usepackage[mathscr]{eucal}
\usepackage{amsmath,amsfonts}

\parskip=\smallskipamount
\newtheorem{theorem}{Theorem}[section]

\newtheorem{proposition}[theorem]{Proposition}

\newtheorem{corollary}[theorem]{Corollary}
\newtheorem{lemma}[theorem]{Lemma}

\newtheorem{remark}[theorem]{Remark}

\makeatletter
\@addtoreset{equation}{section}
\makeatother

\newcommand{\BB}{{\mathbb B}}
\newcommand{\CC}{{\mathbb C}}
\newcommand{\NN}{{\mathbb N}}

\newcommand{\ZZ}{{\mathbb Z}}
\newcommand{\DD}{{\mathbb D}}

\newcommand{\FF}{{\mathbb F}}

\newcommand{\PP}{{\mathbb P}}

\newcommand{\cA}{{\mathcal A}}

\newcommand{\cC}{{\mathcal C}}
\newcommand{\cD}{{\mathcal D}}
\newcommand{\cE}{{\mathcal E}}
\newcommand{\cF}{{\mathcal F}}
\newcommand{\cG}{{\mathcal G}}
\newcommand{\cH}{{\mathcal H}}
\newcommand{\cK}{{\mathcal K}}
\newcommand{\cL}{{\mathcal L}}
\newcommand{\cM}{{\mathcal M}}
\newcommand{\cN}{{\mathcal N}}
\newcommand{\cP}{{\mathcal P}}

\newcommand{\cT}{{\mathcal T}}

\newcommand{\cV}{{\mathcal V}}
\newcommand{\cY}{{\mathcal Y}}
\newcommand{\cX}{{\mathcal X}}
\newcommand{\cW}{{\mathcal W}}
\newcommand{\cZ}{{\mathcal Z}}

\newdimen\expt
\def\boxit#1{\setbox0\hbox{$\displaystyle{#1}$}
      \hbox{\lower.4\expt
 \hbox{\lower3\expt\hbox{\lower\dp0
      \hbox{\vbox{\hrule height.4\expt
 \hbox{\vrule width.4\expt\hskip3\expt
      \vbox{\vskip3\expt\box0\vskip2\expt}%
 \hskip3\expt\vrule width.4\expt}\hrule height.4\expt}}}}}}
 
\begin{document}

 \pagestyle{myheadings}
\markboth{ Gelu Popescu}{ Entropy   and multivariable interpolation }


\title [  Entropy   and multivariable interpolation  ] 
{ Entropy    and multivariable interpolation} 
  \author{Gelu Popescu}
\date{May 20, 2003}
\thanks{Research supported in part by an NSF grant}

\address{Department of Mathematics, The University of Texas 
at San Antonio \\ San Antonio, TX 78249, USA}
\email{\tt gpopescu@math.utsa.edu}

\begin{abstract}
We define a new notion of  entropy for operators
 on Fock spaces and positive definite multi-Toeplitz kernels on free semigroups. 
 This is studied in connection with
 factorization theorems for (multi-Toeplitz, multi-analytic, etc.) operators on Fock
 spaces.
 These results lead to entropy inequalities  and entropy formulas for positive definite
 multi-Toeplitz kernels on free semigroups (resp.~multi-Toeplitz operators) and 
 consequences concerning the extreme points of the unit ball of the noncommutative
 analytic Toeplitz algebra $F_n^\infty$.
 
 We obtain several geometric characterizations of the multivariable central 
 intertwining lifting, a maximum principle, and a permanence principle for the 
 noncommutative commutant lifting theorem.
 Under certain natural conditions, we find explicit forms for the maximal
  entropy
 solution (and its entropy) for this multivariable commutant lifting theorem.
 
 All these results are used to solve maximal entropy interpolation problems in several
 variables. We obtain explicit forms for the maximal entropy
 solution (as well as its entropy) of the Sarason, Carath\' eodory-Schur, 
 and Nevanlinna-Pick type interpolation problems for the 
 noncommutative (resp.~commutative) analytic Toeplitz algebra $ F_n^\infty$
 (resp.~$W_n^\infty$) and their tensor products  with $B(\cH, \cK)$.
 In particular, we provide explicit forms for the maximal entropy
 solutions of  several    interpolation (resp.~optimization)  problems  
     on the unit  ball of $\CC^n$.

\end{abstract}

\maketitle

\clearpage

\tableofcontents

 \clearpage
\section*{Introduction}

Let $H_n$ be an $n$-dimensional complex  Hilbert space with orthonormal basis
$e_1$, $e_2$, $\dots,e_n$, where $n\in \{1,2,\dots\}$ or $n=\infty$.
  We consider the full Fock space  of $H_n$ defined by
$$F^2(H_n):=\bigoplus_{k\geq 0} H_n^{\otimes k},$$ 
where $H_n^{\otimes 0}:=\CC 1$ and $H_n^{\otimes k}$ is the (Hilbert)
tensor product of $k$ copies of $H_n$.
Define the left creation 
operators $S_i:F^2(H_n)\to F^2(H_n), \  i=1,\dots, n$,  by
$$
 S_i\psi:=e_i\otimes\psi, \  \psi\in F^2(H_n).
$$
 The  noncommutative analytic Toeplitz algebra   $F_n^\infty$ 
  and  its norm closed version
  (the noncommutative disc
 algebra  $\cA_n$)  were introduced by the author   \cite{Po-von} in connection
   with a multivariable noncommutative von Neumann inequality.
$F_n^\infty$  is the algebra of left multipliers of the full Fock space
$F^2(H_n)$  and  can be identified with 
 the
  weakly closed  (or $w^*$-closed) algebra generated by the left creation operators
   $S_1,\dots, S_n$ on the full Fock space $F^2(H_n)$,
    and the identity.
     The noncommutative disc algebra $\cA_n$ is 
    the  norm closed algebra generated by the left creation operators
   $S_1,\dots, S_n$  
    and the identity.
     When $n=1$, $F_1^\infty$ 
   can be identified
   with $H^\infty(\DD)$, the algebra of bounded analytic functions
    on the open unit disc. The algebra $F_n^\infty$ can be viewed as a
     multivariable noncommutative 
    analogue of $H^\infty(\DD)$.
We should add that the algebra $F_n^\infty$ shares many
 properties with $H^\infty(\DD)$. There are many analogies with the invariant
  subspaces of the unilateral 
 shift on $H^2(\DD)$, inner-outer factorizations,
  analytic operators, Toeplitz operators, $H^\infty(\DD)$--functional
   calculus, bounded (resp.~spectral) interpolation, etc.
The  noncommutative analytic Toeplitz algebra   $F_n^\infty$ 
  has  been studied
    in several papers
\cite{Po-charact},  \cite{Po-multi},  \cite{Po-funct}, \cite{Po-analytic},
\cite{Po-disc}, \cite{Po-poisson}, 
 \cite{ArPo}, and recently in
  \cite{DP1}, \cite{DP2},   \cite{DP},
  \cite{ArPo2}, \cite{Po-curvature},  \cite{DKP},  \cite{PPoS}, and \cite{Po-similarity}.

We established a strong connection between the algebra $F_n^\infty$
  and the function theory on the open unit ball  
  $$\BB_n:=\{(\lambda_1,\ldots, \lambda_n)\in \CC^n:\  
  |\lambda_1|^2+\cdots +|\lambda_n|^2<1\},
  $$
   through the noncommutative von Neumann inequality \cite{Po-von} 
   (see also  \cite{Po-funct}, \cite{Po-disc}, \cite{Po-poisson}, and \cite{Po-tensor}).
    In particular, we proved that there is a completely contractive
     homomorphism $\Phi:F_n^\infty\to H^\infty(\BB_n)$
     defined by
 $$
       [\Phi(f(S_1,\dots, S_n))](\lambda_1,\dots, \lambda_n)=f(\lambda_1,\dots,
       \lambda_n)$$ for any
 $ f(S_1,\dots, S_n)\in F_n^\infty$ and
  $ (\lambda_1,\dots, \lambda_n)\in \BB_n.
$
 A characterization of the analytic functions in the range of the map
   $\Phi$ was obtained  in \cite{ArPo2}, and independently in \cite{DP}.   
   Moreover, it was proved that the quotient $F_n^\infty/{\ker\Phi}$ is an
  operator algebra which can be identified with
   $W_n^\infty:=P_{F_s^2(H_n)} F_n^\infty|_{F_s^2(H_n)}$,
  the compression  of $ F_n^\infty$ to the symmetric Fock space
   ${F_s^2(H_n)} \subset F^2(H_n)$. 
  In  \cite{Po-poisson}, \cite{Arv1}, \cite{Arv2},  \cite{ArPo2}, \cite{ArPo1}, \cite{DP},  
  \cite{Po-tensor}, \cite{Po-curvature}
  and  \cite{PPoS}, a good case is made 
  that the appropriate
        multivariable commutative analogue
   of $H^\infty(\DD)$  is the algebra $W_n^\infty$, which 
    was   also proved to be  the $w^*$-closed algebra generated by  the operators
   $B_i:=P_{F_s^2(H_n)} S_i|_{F_s^2(H_n)}, \ i=1,\dots, n$, and the identity. 
    
  Moreover, Arveson showed in  \cite{Arv1}
   that $W_n^\infty$ can be seen as the algebra of all analytic multipliers
    of $F_s^2(H_n)$, when  $F_s^2(H_n)$ is identified with a 
    class of analytic functions in $\BB_n$.
    More precisely, the range of the homomorphism  $\Phi$ is the multiplier algebra of the 
reproducing kernel Hilbert space 
with reproducing kernel $K_n: \BB_n\times \BB_n\to \CC$ defined by
 $$
 K_n(z,w):= {\frac {1}
{1-\langle z, w\rangle_{\CC^n}}}, \qquad z,w\in \BB_n.
$$ 

Interpolation problems for the noncommutative analytic Toeplitz algebra
 $F_n^\infty$ were first considered in \cite{Po-analytic}, where   we  
 obtained a
Carath\'eodory  \cite{Ca} type interpolation theorem in this setting. In  1997,
 Arias and
 the author \cite{ArPo2} (see also \cite{ArPo1} and \cite{Po-interpo}) extended Sarason's result
  \cite{S} and obtained a distance formula to 
an arbitrary WOT-closed ideal in $F_n^\infty$ as well as    
a Nevanlinna-Pick type interpolation
theorem (see \cite{N}) for the  noncommutative analytic Toeplitz algebra
 $F_n^\infty$. 
Using different methods, Davidson and Pitts  proved   these
  results  independently  in  \cite{DP}.
Let us mention that, recently, interpolation problems for 
 $F_n^\infty$ (resp.~$W_n^\infty$)  and interpolation problems on 
  the unit  ball 
$\BB_n$  were  also considered in  
\cite{AMc2}, \cite{Po-spectral},  \cite{Po-tensor}, \cite{Po-lifting},
  \cite{BTV}, \cite{Po-central},  \cite{Po-meromorphic}, \cite{Po-nehari}, 
   \cite{BB},  and \cite{EP} .

Due to the connection between  the algebras $F_n^\infty$, $W_n^\infty$, 
and $H^\infty(\BB_n)$, 
the interpolation problems for the noncommutative analytic Toeplitz algebra
 $F_n^\infty$   have appropriate versions
 for the commutative  Toeplitz algebra $W_n^\infty$ and  the Hardy space 
 $H^\infty(\BB_n)$ (resp. some classes of bounded
 analytic functions in $\BB_n$). 
 This  claim is supported by  the following  papers: \cite{ArPo1}, \cite{ArPo2}, 
 \cite{Po-interpo},
   \cite{Po-spectral}, \cite{Po-tensor}, \cite{Po-lifting}, 
 \cite{Po-central},  \cite{Po-meromorphic},  and \cite{Po-nehari}.

\smallskip
   
Let
 $\Theta: F^2(H_n)\otimes \cH\to F^2(H_n)\otimes \cK $ be  a  multi-analytic
  operator,
 i.e., 
   $$
   \Theta(S_i\otimes I_\cH)= (S_i\otimes I_\cK) \Theta,\qquad i=1,\ldots, n.
   $$
   We recall that (see \cite{Po-analytic}, \cite{Po-central}) that 
   $\Theta\in R_n^\infty\bar\otimes B(\cH, \cK)$, where $R_n^\infty$ is the weakly
   closed algebra generated by the right creation operators 
   on the full Fock space, and the 
   identity. Moreover, we have $R_n^\infty=U^* F_n^\infty U$, where $U$ is a unitary operator
   (see Section \ref{entro}).
  If $\|\Theta\|\leq 1$ and   $\dim \cH<\infty$,  then  we define  the
   prediction entropy of
   the multi-analytic operator $\Theta$   by setting
   $$
   E(\Theta):= \ln \det \Delta(\Theta),
   $$
   where 
 \begin{equation*} 
 \left< \Delta(\Theta) x,x\right>:=\inf \{ \left< (I-\Theta^*\Theta) (x-p), x-p\right>:
 \ p\in F^2(H_n)\otimes \cH,\  
 p(0)=0\}
 \end{equation*}
 for any $x\in \cH$.
 It turns out (using Szeg\" o's theorem)  that in the particular case
  when $n=1$ and  $\cH= \cK=\CC$ we have 
 $$
 E(f)= \frac{1} {2\pi} \int_{-\pi}^{\pi} \ln (1-|f(e^{it})|^2)~dt,
 $$
  which is the classical 
 definition for the entropy
 of  $f\in H^\infty(\DD)$ with $\|f\|\leq 1$.
 
 There is an extensive literature 
 (see   \cite{Bu}, \cite{Aro}, \cite{DG1}, \cite{DG2},
   \cite{FF-book}, \cite{EGL}, \cite{FFG},  \cite{FFGK-book}, \cite{RSW}, etc.)
 concerning the maximal 
 entropy solutions for  the classical interpolation
   problems of Carath\' eodory-Schur (\cite{Ca}, 
 \cite{Sc}),
  Nevanlinna-Pick (\cite{N}), Nehari (\cite{Ne}), Sarason (\cite{S}),  and the  more general commutant lifting theorem of 
  Sz.-Nagy and Foia\c s  (\cite{SzF}). 
 Maximal entropy solutions have also played an important role in control theory 
(see \cite{GM1}, \cite{GM2}, \cite{Fra}, etc.). 

 The main goal of this paper is to find  the maximal entropy solution 
 in the noncommutative commutant lifting theorem (see \cite{Po-isometric}, \cite{Po-intert})
 and to use it in order to get  the maximal entropy solutions  for the Sarason,
  Carath\' eodory-Schur, and
 Nevanlinna-Pick type interpolation problems for the noncommutative (resp.~commutative)
 analytic Toeplitz algebra $F_n^\infty$ (resp.~$W_n^\infty$). Moreover, under certain 
 natural conditions, we obtain explicit forms for 
  the maximal entropy solutions  and their entropy. In particular,  we solve maximal entropy 
  interpolation problems on the unit ball of $\CC^n$. 
  
  Our investigation is based on multivariable noncommutative dilation theory
 (see  \cite{B}, \cite{F},  \cite{Po-models},
  \cite{Po-isometric}, \cite{Po-charact},  \cite{Po-intert},  \cite{Po-structure}),
  harmonic analysis on Fock spaces (see
  \cite{Po-charact},  \cite{Po-multi},  \cite{Po-intert},  \cite{Po-analytic},
  \cite{Po-curvature},
 \cite{ArPo}, 
  \cite{DP1},  and \cite{DP2}), 
  the results 
   of Foia\c s, Frazho, and Gohberg  \cite{FFG}
 (see also  \cite{FF-book}, \cite{FFGK-book}, and \cite{SzF-book}) 
 on maximal entropy interpolation (case $n=1$), 
  and the 
  classical results on interpolation and commutant lifting theorem
  (see \cite{Ca},  \cite{Sc}, \cite{N}, \cite{Ne}, \cite{S}, 
  \cite{SzF}, \cite{SzF-book}, \cite{FF-book}, \cite{FFGK-book}, etc.)

  The paper is organized  in three chapters.
In the first chapter, 
we define a new notion of  entropy for operators
 on Fock spaces and positive definite  multi-Toeplitz kernels on free semigroups. 
 This is studied in connection with
 factorization theorems for (multi-Toeplitz, multi-analytic, etc.) operators on Fock
 spaces.
 These results lead to entropy inequalities  and entropy formulas for positive definite
 multi-Toeplitz kernels on free semigroups (resp.~multi-Toeplitz operators) and 
 consequences concerning the extreme points of the unit ball of the noncommutative
 analytic Toeplitz algebra $F_n^\infty$.

More precisely, 
in Section \ref{entro},
 we prove the existence of maximal outer factors 
for arbitrary positive multi-Toeplitz operators on Fock spaces. A connection between the 
Szeg\" o infimum  and maximal outer factors of  multi-Toeplitz operators
  via the prediction-error 
operators is considered. 
 These results extend some classical 
 results (see \cite{Sz1}, \cite{Sz2}, \cite{Sz3}, \cite{GSz}, \cite{HL1}, \cite{HL2}, 
 \cite{K1}, \cite{K2}, \cite{SzF-book}, \cite{W}, \cite{WM1}, and  \cite{WM2}) as well as 
 some extensions to Fock spaces from (\cite{Po-analytic} and \cite{Po-structure}). 
  We introduce the notion of prediction entropy for 
 positive multi-Toeplitz operators $T\in B(F^2(H_n)\otimes \cE)$, where $\cE$ is a Hilbert
 space  with  $\dim \cE<\infty$.
 In particular,  we prove that
  the entropy of $T$ satisfies the equation 
 $$
 e(T)= \ln \det [\varphi(0)^* \varphi(0)],
 $$
 where $M_\varphi$ is the maximal outer factor of $T$.
  Finally, we provide   an explicit form  for   the  square 
 outer spectral factor $M_\varphi$  corresponding to a   strictly positive
  multi-Toeplitz operator $T$, i.e., $T=M_\varphi^* M_\varphi$.

In the next section,
%
%
 we provide factorization results for  bounded linear operators in 
$B(\cH_1, F^2(H_n)\otimes \cH_2)$, where $\cH_1$   is 
a finite dimensional Hilbert space, 
and for multi-analytic operators, generalizing 
classical results from
 \cite{HL1}, \cite{HL2}, \cite{RR}, and  \cite{W},  as well as some extensions 
 to Fock spaces 
 (see \cite{Po-multi}, \cite{Po-analytic},  and \cite{Po-structure}). 
A noncommutative multivariable analogue of Robinson's minimum energy
 delay principle (see  \cite{R} and \cite{FFG})
for outer operators on Fock spaces is obtained.

  These results are used to
  %
   prove entropy inequalities for positive definite
 multi-Toeplitz kernels
on the free semigroup with $n$ generators $\FF_n^+$ and  multi-Toeplitz operators on Fock spaces.
 We extend the classical result (see \cite{H}) which stated for
  the Hardy space  $H^2(\DD)$ says that 
 if $f\in H^2(\DD)$, then $\ln|f(e^{it})|$ is integrable and 
 \begin{equation}\label{eni}
 \frac {1}{2\pi} \int_{-\pi}^ \pi \ln|f(e^{it})| ~dt \geq \ln |f(0)|.
 \end{equation}
 Next, we give a characterization for the outer operators in 
 $B(\cE, F^2(H_n)\otimes \cE)$ if $\dim \cE<\infty$.
In particular, we find a  
noncommutative  multivariable analogue of the   classical result saying that
a function $f\in H^2(\DD)$ is outer if and only if $f(0)\neq 0$   and 
the equality
 holds in \eqref{eni}.

 It is well-known  \cite{H} that a function $f\in H^\infty(\DD)$ is an extreme point of 
 the unit ball of $H^\infty(\DD)$ if and only if
 $$
 \frac{1}{2\pi} \int_{-\pi}^\pi \ln(1-|f(e^{it})|^2)~dt=-\infty.
 $$
In    Section  
\ref{extreme points},
 we prove some results concerning the extreme points
 of the unit ball of $F_n^\infty$. In particular, we show that if $\varphi\in F_n^\infty$,
 $\|\varphi\|\leq 1$, and the entropy $E(\varphi)=-\infty$, then $\varphi$ is 
 an extreme point of the unit ball of $F_n^\infty$.

 In the second chapter of this paper, we obtain several geometric 
 characterizations for the  multivariable central 
 intertwining lifting, a maximum principle, and a permanence principle for the 
 noncommutative commutant lifting theorem.
 Under certain natural conditions, we find explicit forms for the maximal entropy
 solution of this multivariable commutant lifting theorem, and concrete formulas 
 for its
 entropy.
 
 Let us recall the noncommutative commutant lifting theorem for row contractions
  \cite{Po-isometric}, \cite{Po-intert} 
  (see \cite{SzF} for the classical case $n=1$).
  Let $\cT:=[T_1~\cdots~T_n]$, \ $T_i\in B(\cH)$, and 
  $\cT':=[T_1'~\cdots~T_n']$, \ $T_i'\in B(\cH')$, be row contractions
   and let 
  $\cV:=[V_1~\cdots~V_n]$, \ $V_i\in B(\cK)$, and 
  $\cV':=[V_1'~\cdots~V_n']$, \ $V_i'\in B(\cK')$, be the corresponding
   minimal isometric 
  dilations.
 The noncommutative commutant lifting theorem states that if
  $A\in B(\cH, \cH')$
 is an operator satisfying
 $$
 AT_i=T_i' A,\qquad i=1,\ldots, n,
 $$
 then there exists an operator $B\in B(\cK,\cK')$ with the following properties:
 \begin{enumerate}
 \item[(i)] $BV_i=V_i'B$ ~for any $i=1,\ldots, n$;
 \item[(ii)] $B^*|\cH'= A^*$;
 \item[(iii)] $\|B\|= \|A\|$.
 \end{enumerate}
 As in the classical case, the general setting of the 
 noncommutative commutant lifting theorem  can be reduced to the case when 
 $\cT:=[T_1~\cdots~T_n]$ is a row isometry (see \cite{Po-nehari}).

In Section \ref{geom}, we present some results concerning the
 geometric structure of the
 intertwining liftings in the
noncommutative commutant lifting theorem.
It is shown that there is a one-to-one correspondence between the set of all 
multivariable
intertwining liftings $B$  with tolerance $t>0$ (i.e., $\|B\|\leq t$) and
 certain families of  contractions
$\{C_k\}_{k=1}^\infty$ and $\{\Lambda_\alpha\}_{\alpha\in \FF_n^+}$. 
 Moreover, we prove  that the intertwining lifting corresponding to
   the parameters $C_k=0$, $k=1,2,\ldots$, 
   (resp.~$\Lambda_\alpha=0$, $\alpha\in \FF_n^+$)  coincides with the central 
   intertwining lifting  with tolerance $t$,
   for which we have an explicit form.
   The geometric structure of the  central 
   intertwining lifting will play a very important role in our investigation.
  We prove a maximum principle  for
 the noncommutative commutant lifting theorem, which also provides
  a new characterization for the central intertwining lifting. 
  This   result is  used   to prove 
   a  permanence principle
 for the central intertwining lifting, which generalizes the permanence principle
  for the Carath\' eodory
 interpolation problem  (see \cite{EGL}  and \cite{FFG}) (case $n=1$).
 Applications of this principle  will be considered in  the  last chapter.

  The next step is to   obtain  explicit formulas for the  quasi outer  spectral
   factor  for the defect operator  $t^2I- B_c^* B_c$ of the central 
   intertwining lifting $B_c$ with tolerance $t>0$.
    This leads, 
   in the next section,  to concrete formulas for 
  the entropy of $B_c$   as well as to
   a maximum principle for the noncommutative commutant
    lifting
  theorem  with respect to non-minimal isometric liftings.

  In Section \ref{max-entro-sol}, we present  one of the main results of this paper. 
 Using the maximum principle,  we prove that the 
 central intertwining lifting  $B_c$ is the maximal entropy solution for 
 the noncommutative commutant lifting theorem, when 
 $\cT:=[S_1\otimes I_\cE~\cdots ~ S_n\otimes I_\cE]$ with $\dim \cE<\infty$.
 Based on several results of this paper, we are  led to concrete formulas for 
  the entropy of $B_c$   and, under a certain condition of stability,  to
   a maximum principle
   and  a characterization (in terms of  entropy) of the   central intertwining lifting
  $\tilde{B_c}$ with respect to non-minimal isometric liftings.

 The results of the first two chapters  are used 
   to solve maximal entropy interpolation problems in several
 variables, in the last chapter  of this paper.
  We obtain explicit forms for the maximal entropy
 solution (as well as its entropy) of the Sarason \cite{S}, 
 Carath\' eodory-Schur \cite{Ca},
 \cite{Sc}, 
 and Nevanlinna-Pick  \cite{N} type interpolation problems for the 
 noncommutative (resp.~commutative) analytic Toeplitz algebra $ F_n^\infty$
 (resp.~$W_n^\infty$) and their tensor products  with $B(\cH, \cK)$, 
 the set of all bounded linear operators 
 acting on Hilbert spaces.
 Moreover, under certain conditions, we also find explicit forms for
  the corresponding classical optimization problems, 
  in our multivariable noncommutative (resp.~commutative)
  setting.
  
  In particular, we provide explicit forms for the maximal entropy
 solutions of  several    interpolation problems  
     on the unit  ball of $\CC^n$. Finnaly, we apply our permanence principle
     to the Nevanlinna-Pick interpolation problem on the unit ball.

\clearpage
 \section{Operators on Fock spaces and their entropy}

 A new notion of  entropy for operators
 on Fock spaces and positive multi-Toeplitz kernels on free semigroups 
 is defined 
  and  studied in connection with
 factorization theorems for (multi-Toeplitz, multi-analytic, etc.) 
 operators on Fock
 spaces.
 These results lead to entropy inequalities  and entropy formulas for positive 
 multi-Toeplitz kernels on free semigroups (resp.~multi-Toeplitz operators) and 
 consequences concerning the extreme points of the unit ball of 
 the noncommutative
 analytic Toeplitz algebra $F_n^\infty$.

\subsection{Entropy  and spectral factorization for multi-Toeplitz operators  }
\label{entro}
 
 In this section, we define  the  notion of  prediction entropy 
 for positive multi-Toeplitz 
 operators $T\in B(F^2(H_n)\otimes \cE)$ with $\dim \cE<\infty$.
 We prove that   there is a multi-analytic
 operator $M_\varphi\in  R_n^\infty\bar \otimes B(\cE)$  such that
 $M_\varphi ^* M_\varphi \leq T$ if and only if the entropy of $T$ satisfies 
 $e(T)>-\infty$. This is based on Theorem \ref{sub-fac} which proves the existence of maximal
 outer factors for arbitrary positive  multi-Toeplitz operators on Fock spaces.
 Moreover, we prove a Szeg\" o type infimum theorem for 
  arbitrary positive  multi-Toeplitz operators (see Theorem \ref{Sze}), and 
  provide an explicit form for the 
 square outer  spectral factor  $M_\varphi$  corresponding to  a strictly positive
  multi-Toeplitz operator $T$, i.e., $M_\varphi^* M_\varphi= T$ 
  (see Theorem \ref{outer}).  All these results are needed 
  in the next sections.

Let $\FF_n^+$ be the unital free semigroup on $n$ generators 
$g_1,\dots,g_n$, and the identity $g_0$.
The length of $\alpha\in\FF_n^+$ is defined by
$|\alpha|:=k$, if $\alpha=g_{i_1}g_{i_2}\cdots g_{i_k}$, and
$|\alpha|:=0$, if $\alpha=g_0$.
We also define
$e_\alpha :=  e_{i_1}\otimes e_{i_2}\otimes \cdots \otimes e_{i_k}$  
 and $e_{g_0}= 1$.
It is  clear that 
$\{e_\alpha:\alpha\in\FF_n^+\}$ is an orthonormal basis of the full Fock space
 $F^2(H_n)$.

 Let $\cE$ be a Hilbert space and let $T$ be a positive multi-Toeplitz operator
  on $F^2(H_n)\otimes \cE$, i.e.,
 $$
 (S_i\otimes I_\cE)^* T (S_j\otimes I_\cE)=\delta_{ij} I
 $$ 
 for any $i,j=1,\ldots, n$.
 Define the positive operator $\Delta_T: \cE\to \cE$ by setting
 \begin{equation}\label{DT}
 \left< \Delta_T x,x\right>:=\inf \{ \left< T(x-p), x-p\right>:\ p\in \cP(\cE), 
 \ p(0)=0\}
 \end{equation}
 for any $x\in \cE$, where $\cP(\cE)$ denotes the set of  polynomials
 ~$p=\sum\limits_{|\alpha|\leq m} e_\alpha \otimes h_\alpha$, \ $m=0,1,\ldots$, 
   in  $F^2(H_n)\otimes \cE$ and 
 $p(0):=h_{g_0}$.
 We remark that if $T$ is a normalized positive multi-Toeplitz operator, i.e., 
 $P_\cE T|\cE= I_\cE$, then $\Delta_T$ coincides with the prediction-error
  operator of the stationary process determined by $T$ (see \cite{Po-structure}).

 When $\cE$ is finite dimensional, we define the prediction entropy of
 the positive multi-Toeplitz operator $T$ by
 \begin{equation}\label{ET}
 e(T):= \ln \det\Delta_T.
 \end{equation} 
 Note that the prediction entropy is different from the   entropy
 defined in \cite{Po-structure}.

We need to recall from  \cite{Po-charact},
\cite{Po-multi}, \cite{Po-von},  \cite{Po-funct}, and  \cite{Po-analytic}  a few facts
 concerning multi-analytic operators on Fock spaces.
   We say that 
 a bounded linear
  operator 
$M\in B(F^2(H_n)\otimes \cK, F^2(H_n)\otimes \cK')$ is  multi-analytic
if 
\begin{equation}
M(S_i\otimes I_\cK)= (S_i\otimes I_{\cK'}) M\quad 
\text{\rm for any }\ i=1,\dots, n.
\end{equation}
Note that $M$ is uniquely determined by the operator
$\theta:\cK\to F^2(H_n)\otimes \cK'$, which is   defined by ~$\theta k:=M(1\otimes k)$, \ $k\in \cK$, 
 and
is called the  symbol  of  $M$. We denote $M=M_\theta$. Moreover, 
$M_\theta$ is 
 uniquely determined by the ``coefficients'' 
  $\theta_{\alpha}\in B(\cK, \cK')$, which are given by
 $$
\left< \theta_{\tilde\alpha}k,k'\right>:= \left< \theta k, e_\alpha 
\otimes k'\right>=\left< M_\theta(1\otimes k), e_\alpha 
\otimes k'\right>,\quad 
k\in \cK,\ k'\in \cK',\ \alpha\in \FF_n^+,
$$
where $\tilde\alpha$ is the reverse of $\alpha$, i.e., $\tilde\alpha= g_{i_k}\cdots g_{i_1}$ if
$\alpha= g_{i_1}\cdots g_{i_k}$.
Note that 
$$
\sum\limits_{\alpha \in \FF_n^+} \theta_{\alpha}^*   \theta_{\alpha}\leq 
\|M_\theta\|^2 I_\cK.
$$
 If $T_1,\dots,T_n\in B(\cH)$  (the algebra of all bounded  
  linear operators on the Hilbert space $ \cH$), define 
$T_\alpha :=  T_{i_1}T_{i_2}\cdots T_{i_k}$,
if $\alpha=g_{i_1}g_{i_2}\cdots g_{i_k}$ and 
$T_{g_0}:=I_\cH$.  
 We can associate with $M_\theta$ a unique formal Fourier expansion 
\begin{equation}\label{four}
M_\theta\sim \sum_{\alpha \in \FF_n^+} R_\alpha \otimes \theta_{\alpha},
\end{equation}
where $R_i:=U^* S_i U$, \ $i=1,\ldots, n$, are the right creation operators
on $F^2(H_n)$ and 
$U$ is the (flipping) unitary operator on $F^2(H_n)$ mapping 
  $e_{i_1}\otimes e_{i_2}\otimes\cdots\otimes e_{i_k}$ into
 $e_{i_k}\otimes\cdots\otimes e_{i_2}\otimes e_{i_1} $.
 Since $M_\theta$ acts like its Fourier representation on ``polynomials'', 
  we will identify them for simplicity.
The set of multi-analytic operators in 
$B(F^2(H_n)\otimes \cK,
F^2(H_n)\otimes \cK')$  coincides  with   
$R_n^\infty\bar \otimes B(\cK,\cK')$, where $R_n^\infty=U^* F_n^\infty U$
(see \cite{Po-analytic}  and \cite{Po-central}).
 A multi-analytic  operator $M_\theta$ (resp. its symbol $\theta$) is called 
 inner if
 $M_\theta$ is an isometry.
 We call $M_\theta$ (resp. its symbol $\theta$ )   outer if 
$$
\bigvee\{(S_\alpha\otimes I_{\cK'})\theta k:\ k\in \cK, \alpha\in \FF_n^+\}
=F^2(H_n)\otimes \cK'.
$$

 We say that a positive multi-Toeplitz operator $T\in B(F^2(H_n)\otimes \cE)$
 has a maximal outer factor if there is  a Hilbert space $\cG$ and an outer 
 multi-analytic operator
 $M_\varphi\in 
 R_n^\infty\bar \otimes B(\cE, \cG)$  with the properties:
 \begin{enumerate}
 \item[(i)]
 $M_\varphi^* M_\varphi \leq T$;
 \item[(ii)]
 If $\cN$ is a Hilbert space and $M_\theta \in R_n^\infty\bar \otimes B(\cE, \cN)$
 satisfies 
 $$
 M_\theta^* M_\theta \leq T,
 $$
 then 
 $M_\varphi^* M_\varphi\geq M_\theta^* M_\theta$.
 \end{enumerate}

 In what follows we  prove the existence of maximal outer factors for arbitrary
  positive multi-Toeplitz operators on Fock spaces.

 \begin{theorem}
 \label{sub-fac}
 If $\cE$ is an arbitrary Hilbert space and  $T\in B(F^2(H_n)\otimes \cE)$ is
  a positive multi-Toeplitz operator, then  $T$  has  a  maximal outer factor
    $M_\varphi\in 
 R_n^\infty\bar \otimes B(\cE, \cG)$ such that
 \begin{equation}\label{M*MT}
 M_\varphi^* M_\varphi \leq T.
 \end{equation}
 Moreover, the  maximal outer factor is  
   uniquely determined up to a unitary diagonal 
  multi-analytic operator.
 \end{theorem}

 \begin{proof}
 Since $T$ is a positive multi-Toeplitz operator, we have
 \begin{equation*}\begin{split}
 \left\|T^{1/2}\left(\sum_{i=1}^n (S_i\otimes I_\cE)h_i\right)\right\|^2&=
  \sum_{i,j=1}^n \left< (S_j^*\otimes I_\cE )T (S_i\otimes I_\cE)h_i, h_j\right>\\
  &= \sum_{i,j=1}^n\left< \delta_{ij} T h_i, h_j\right>= \sum_{i=1}^n \|T^{1/2} h_i\|^2
 \end{split}
 \end{equation*}
 for any $h_i\in F^2(H_n)\otimes \cE)$, \ $i=1,\ldots, n$.
Let $\cY:= \overline{ T^{1/2}(F^2(H_n)\otimes \cE)}$ and note that there are unique isometries
$V_i\in B(\cY)$, \ $i=1,\ldots, n$, such that 
\begin{equation}\label{ViT}
V_iT^{1/2} =T^{1/2} (S_i\otimes I_\cE),\quad i=1,\ldots, n.
\end{equation} 
The above calculations show that
$$
\left\|\sum_{i=1}^n V_i T^{1/2} h_i\right\|^2= \sum_{i=1}^n \|T^{1/2}h_i\|^2
$$
for any $h_i\in \cY$. Hence, $V_1,\ldots, V_n$ are isometries with orthogonal ranges.
 According to the Wold type decomposition
 for isometries with orthogonal ranges (see \cite{Po-isometric}), we have
 $V_i=U_i\oplus W_i$, \ $i=1,\ldots, n$,
 with respect to the orthogonal decomposition
 $\cY= \cY_0\oplus \cY_1$, where 
 $\cY_0:= \oplus_{\alpha\in \FF_n^+} V_\alpha \cM_T$ and 
 \begin{equation}\label{M_T}
 \cM_T=\cY\ominus \bigvee_{i=1}^n T^{1/2} (S_i\otimes I_\cE)(F^2(H_n)\otimes \cE).
 \end{equation}
 The subspaces $\cY_0$ and $\cY_1$ are reducing for each 
  operator $V_i$, \ $i=1,\ldots, n$. Moreover,  
  $$
  \sum_{i=1}^n W_iW_i^*=I_{\cY_1}
  $$ 
  and $\{U_i\}_{i=1}^n$ is unitarily equivalent 
 to the orthogonal shift $\{S_i\otimes I_{\cM_T}\}_{i=1}^n$.
 Hence, we have
 $$
 \Phi U_i= (S_i\otimes I_{\cM_T}) \Phi, \quad i=1,\ldots, n,
 $$
 where $\Phi: \cY_0\to F^2(H_n)\otimes \cM_T$ is the Fourier transform
  defined by 
  \begin{equation}\label{fouri}
  \Phi(V_\alpha \ell_\alpha)
 =e_\alpha\otimes \ell_\alpha,\quad
  \ell_\alpha\in \cM_T, \ \alpha\in \FF_n^+.
  \end{equation}
Since $\cY_0$ reduces each operator $V_i$, $i=1,\ldots, n$, it follows that
$$
U_iP_{\cY_0}=V_iP_{\cY_0}=P_{\cY_0}V_i, \quad i=1,\ldots, n.
$$
Hence, we infer that
\begin{equation*}\begin{split}
\Phi P_{\cY_0}T^{1/2} (S_i\otimes I_\cE)&= \Phi P_{\cY_0}V_i T^{1/2}= \Phi U_i P_{\cY_0}T^{1/2}\\
&= (S_i\otimes I_{\cM_T})\Phi P_{\cY_0}T^{1/2}
\end{split}
\end{equation*}
 for any $i=1,\ldots, n$.
 This shows that  the operator $\Phi P_{\cY_0}T^{1/2}$ is  multi-analytic   and, 
 according to \cite{Po-analytic},
 there exists $M_\varphi\in R_n^\infty\bar \otimes B(\cE, \cM_T)$ such that 
 \begin{equation}
 \label{PHIP}
 \Phi P_{\cY_0}T^{1/2} =M_\varphi.
 \end{equation}
 Since  the operator $P_{\cY_0}T^{1/2} $ has dense range in $\cY_0$,
  it follows that $M_\varphi$ 
 is outer.
   Finally, relation \eqref{PHIP} impies $M_\varphi^* M_\varphi\leq T$.
    
  To prove the maximality property,  let $\cN$ be a Hilbert space and let 
   $M_\theta \in R_n^\infty\bar \otimes B(\cE, \cN)$
 be such that 
 $
 M_\theta^* M_\theta \leq T.
 $  
   Define $X:\cY\to F^2(H_n)\otimes \cE$ by setting
   $$
   X\left(\sum_{|\sigma| \leq m} V_\sigma T^{1/2} y_\sigma\right):= 
   \sum_{|\sigma| \leq m} (S_\sigma\otimes I_\cN)\theta y_\sigma, \quad y_\sigma\in \cE.
   $$
   Note that, for any  $y_\sigma\in \cE$, we have 
   \begin{equation*}
   \begin{split}
   \left\|  X\left(\sum_{|\sigma| \leq m} V_\sigma T^{1/2} y_\sigma\right)\right\|^2
   &=
   \left< M_\theta^* M_\theta \left(\sum_{|\sigma|\leq m} e_\sigma\otimes y_\sigma\right),
   \sum_{|\sigma|\leq m} e_\sigma\otimes y_\sigma\right>\\
   &\leq 
   \left< T\left(\sum_{|\sigma|\leq m} e_\sigma\otimes y_\sigma\right),
   \sum_{|\sigma|\leq m} e_\sigma\otimes y_\sigma\right>\\
   &=
   \left\|\sum_{|\sigma| \leq m} V_\sigma T^{1/2} y_\sigma\right\|^2.
   \end{split}
   \end{equation*}
   This shows that $X$ extends to a contraction from $\cY$ to  $F^2(H_n)\otimes \cN$.
   One can easily check that 
   $
   XV_i= (S_i\otimes I_\cN) X$, \ $ i=1,\ldots, n$.
   Hence, and using  the Wold decomposition  from  \cite{Po-isometric}, we deduce
   \begin{equation*}
   \begin{split}
   X\cY_1\subseteq \bigcap_{k=0}^\infty X\left(
   \operatornamewithlimits\bigoplus\limits_{|\alpha|=k} V_\alpha \cY\right)&= 
   \bigcap_{k=0}^\infty \operatornamewithlimits\bigoplus_{|\alpha|=k} (S_\alpha\otimes I_\cE) X\cY\\
   &\subseteq \bigcap_{k=0}^\infty \operatornamewithlimits\bigoplus_{|\alpha|=k} (S_\alpha\otimes I_\cE)
   (F^2(H_n)\otimes \cE)= \{0\}.
   \end{split}
   \end{equation*}
   Therefore, $X|\cY_1=0$.
    Now, taking into account the above considerations, we have
   \begin{equation*}
   \begin{split}
   \left< M_\theta^* M_\theta \left(\sum_{|\sigma|\leq m} 
   e_\sigma\otimes y_\sigma\right)\right.,
   &\left.\sum_{|\sigma|\leq m} e_\sigma\otimes y_\sigma\right>
   =
   \left\|X\left(\sum_{|\sigma|\leq m} P_{\cY_0}V_\sigma
    T^{1/2}y_\sigma \right)\right\|^2\\
    &\leq
    \left\|\sum_{|\sigma|\leq m}\Phi P_{\cY_0}V_\sigma
    T^{1/2}y_\sigma \right\|^2\\
    &=
    \left\|\sum_{|\sigma|\leq m}\Phi P_{\cY_0}
    T^{1/2}(S_\sigma \otimes I_\cE)(1\otimes y_\sigma) \right\|^2\\
    &=
    \left\|\sum_{|\sigma|\leq m} M_\varphi
    (S_\sigma \otimes I_\cE)(1\otimes y_\sigma) \right\|^2\\
    &=
    \left< M_\varphi^* M_\varphi \left(\sum_{|\sigma|\leq m} e_\sigma\otimes y_\sigma\right),
   \sum_{|\sigma|\leq m} e_\sigma\otimes y_\sigma\right>.
   \end{split}
   \end{equation*}
   Therefore, $M_\theta^* M_\theta \leq M_\varphi^* M_\varphi$. The proof is complete. 
 \end{proof}

   Let us remark that  the equality holds in \eqref{M*MT}
if and only if $\cY_0=\cY$.  Another characterization of this fact can be found in 
\cite{Po-multi}.
On the other hand,
 we have a concrete form for the maximal outer factor of $T$, that is,
 $M_\varphi= \Phi P_{\cY_0} T^{1/2}$, where the Fourier transform 
 $\Phi$   and the subspace $\cY_0$ are defined  in  the proof of Theorem \ref{sub-fac}.

 \begin{corollary}\label{entro-fact}
 Let $T\in B(F^2(H_n)\otimes \cE)$
  be a positive multi-Toeplitz operator and assume that $\dim \cE<\infty$.
   Then the entropy $e(T)>-\infty$ if and only if 
  there exists an outer multi-analytic operator $M_\varphi\in 
 R_n^\infty\bar \otimes B(\cE)$ such that
 \begin{equation}\label{M*M}
 M_\varphi^* M_\varphi \leq T.
 \end{equation}
 \end{corollary}
 
 \begin{proof}
 
Assume now that $T=M_\varphi^* M_\varphi$ with
 $M_\varphi\in R_n^\infty\otimes B(\cE)$ an outer operator.
 Then there exists a unitary operator $Z:\cY\to F^2(H_n)\otimes \cE$ satisfying
 $$
 ZT^{1/2}=M_\varphi.
 $$
 Since $M_\varphi$ commutes with each operator  $S_i\otimes I_\cE$, $i=1,\ldots, n$, relation 
 \eqref{ViT} implies
 \begin{equation}\begin{split}
 ZV_iT^{1/2}&= ZT^{1/2} (S_i\otimes I_\cE)= M_\varphi (S_i\otimes I_\cE)\\
 &= (S_i\otimes I_\cE) M_\varphi = (S_i\otimes I_\cE)ZT^{1/2}
 \end{split}
 \end{equation}
 for any $i=1,\ldots, n$.
Hence,  $Z V_i=(S_i\otimes I_\cE) Z$ for any $i=1,\ldots, n$, and  consequently 
$\{V_i\}_{i=1}^n$
 is an orthogonal shift having the same multiplicity as 
$\{S_i\otimes I_\cE\}_{i=1}^n$. Therefore, we have  $\dim \cE= \dim \cM_T$, where 
\begin{equation}
\label{MTY}
\cM_T=\cY\ominus \bigvee_{i=1}^n T^{1/2} (S_i\otimes I_\cE)(F^2(H_n)\otimes \cE)
\end{equation}
is the wandering subspace of $\{V_i\}_{i=1}^n$ (see \cite{Po-isometric}).
 Note that, if  $\dim \cE<\infty$, then we have
\begin{equation*}
\begin{split}
\cM_T&=
\overline { P_{\cM_T}T^{1/2} (F^2(H_n)\otimes \cE)}\\
&=(P_{\cM_T}T^{1/2}\cE) \bigvee (\bigvee_{i=1}^n P_{\cM_T}T^{1/2}
(S_i\otimes \cE )(F^2(H_n)\otimes \cE)\\
&= P_{\cM_T}T^{1/2}\cE.
\end{split}
\end{equation*}
Hence, it is clear that
   the following statements are equivalent:
 \begin{enumerate}
 \item[(i)]  $e(T)>-\infty$;
 \item[(ii)] $\det \Delta_T\neq 0$;
 \item[(iii)] $\dim \cE =\dim \cM_T$.
 \end{enumerate}
Since $e(T)= \ln\det\Delta_T$,  we deduce $e(T)> -\infty$.
Now, if $T\geq M_\varphi^* M_\varphi$, then  we have
$$e(T)\geq e(M_\varphi^* M_\varphi)>-\infty.
$$

Conversely, assume that the entropy $e(T)> -\infty$. Then $\dim \cE= \dim \cM_T$ and 
the result follows from Theorem \ref{sub-fac}.
 \end{proof}

It is well-known that in the classical case ($n=1$)  (see \cite{FFGK-book}),
 we have equality in \eqref{M*M}.
Whether or not  this is  true if $n\geq 2$   remains an open problem.

 In what follows, we prove a   Szeg\" o  type infimum theorem for arbitrary
 positive multi-Toeplitz operators. If $M_\theta$ is a multi-analytic operator with Fourier
 expansion \eqref{four}, then we denote $\theta(0):=\theta_{g_0}$.
 
 \begin{theorem}\label{Sze}
 Let $\cE$ be  a  Hilbert space and  let  $T\in B(F^2(H_n)\otimes \cE)$ be 
  a   positive  multi-Toeplitz operator. 
 Then, for any $h\in \cE$,
 $$
 \inf_{p\in \cP(\cE), ~p(0)=0} \left< T(h-p), h-p\right> =\left<\varphi(0)^* 
 \varphi(0) h,h\right>,
 $$
 where $M_\varphi$ is the maximal outer factor of $T$.
 In particular, if $\dim \cE<\infty $, then the entropy of $T$ satisfies the equality
 $$
 e(T)=\ln \det [\varphi(0)^* \varphi(0)].
 $$
 \end{theorem}

 \begin{proof}
   First we  prove  that 
  \begin{equation}\label{DTo}
  \Delta_T= P_\cE T^{1/2} P_{\cM_T} T^{1/2}|\cE,
  \end{equation}
  where  the operator $\Delta_T$ is defined by relation \eqref{DT}, and
  $$
  \cM_T:= \overline{T^{1/2}(F^2(H_n)\otimes \cE)}\ominus 
  \bigvee_{i=1}^n T^{1/2}(S_i\otimes I_\cE)(F^2(H_n)\otimes \cE).
  $$
  Indeed, we have 
 \begin{equation*} \begin{split}
 \left< \Delta_Th, h\right>&= 
\inf\{\|T^{1/2}( h-p)\|^2: p\in \cP(\cE),\  p(0)=0\}\\
&=\inf\left\{\|T^{1/2}h-y\|^2:\ y\in \bigvee_{i=1}^n T^{1/2} (S_i\otimes I_\cE)
 (F^2(H_n)\otimes \cE)\right\}\\
&= \|P_{\cM_T}T^{1/2} h\|^2=\left< (P_\cE T^{1/2} P_{\cM_T} T^{1/2}|\cE)h,
 h\right>
 \end{split}
 \end{equation*}
 for any $h\in \cE$.
 Hence, and using relation \eqref{PHIP}, we  obtain 
 \begin{equation*}\begin{split}
  \left< \Delta_Th, h\right>&= \|P_{\cM_T}T^{1/2} h\|^2\\
  &=\|P_{1\otimes\cM_T}\Phi P_{\cY_0}T^{1/2} h\|^2= \|P_{1\otimes\cM_T} M_\varphi h\| \\
  &= \left<\varphi(0)^* 
 \varphi(0) h,h\right>,
 \end{split}
 \end{equation*}
  for any $h\in \cE$, where $M_\varphi$ is the maximal outer factor of $T$.
  The last part of the theorem  is now obvious. 
 \end{proof}
   
Given a multi-Toeplitz operator $T\in B(F^2(H_n)\otimes \cH)$,
we say that $F\in R_n^\infty\bar \otimes B(\cH,\cY)$ is a spectral factor of $T$ if 
$T= F^* F$. If, in addition,  $\cH=\cY$ and $F$ is outer, then $F$ is called 
a square
 outer spectral factor of $T$.

In \cite{Po-multi}, we proved that any   strictly positive multi-Toeplitz operator $T$
admits a spectral factor. Based on the proof of Theorem \ref{sub-fac}, we can deduce 
 a stronger result, namely that
$T$ has a  square outer spectral factorization.
 
 \begin{corollary}
 If $\cE$ is a  Hilbert space and   $T\in B(F^2(H_n)\otimes \cE)$ is 
  a strictly positive multi-Toeplitz operator, 
 then 
 there exists  an outer multi-analytic operator $M_\varphi\in 
 R_n^\infty\bar \otimes B(\cE)$ such that
 \begin{equation*}
 M_\varphi^* M_\varphi = T.
  \end{equation*}
 \end{corollary}

 \begin{proof}
 As in the proof of Theorem \ref{sub-fac},  we can define the isometries
  $V_i$, $i=1,\ldots, n$, by
 relation \eqref{ViT}. Since $T^{1/2}$ is invertible, Theorem 2.1 from \cite{Po-funct}
 implies that there exists a unitary operator $U: F^2(H_n)\otimes \cE\to \cY$ such that
 $$
 U^* V_i U= S_i\otimes I_\cE, \quad i=1,\ldots, n.
 $$
 This shows that $\cY= \cY_0$ and relation \eqref{PHIP} becomes
  $\Phi T^{1/2}= M_\varphi$, where $\Phi$ is the unitary operator
   defined be \eqref{fouri}.
 Hence,  $ M_\varphi^* M_\varphi = T$ and the proof is complete.
 \end{proof}

 In the next theorem, we provide   an explicit form  for   the  square 
 outer spectral factor corresponding to a   strictly positive
  multi-Toeplitz operator. 
 
\begin{theorem}\label{outer}
Let    $T\in B(F^2(H_n)\otimes \cE)$ be a strictly positive multi-Toeplitz 
operator.
Then 
$$T=\Theta^* \Theta,
$$
where $\Theta\in R_n^\infty\bar \otimes B(\cE)$ is a square outer spectral 
factor of $T$
given by
$$
\Theta :=(I\otimes N)M_\psi^{-1},
$$
 where $M_\psi$ is an invertible multi-analytic operator with 
  symbol $\psi:\cE\to F^2(H_n)\otimes \cE$  defined by
 $$
 \psi h:= T^{-1}(1\otimes h), \quad h\in \cE,
 $$
 and $N:= (P_\cE T^{-1} |\cE)^{1/2}$.
\end{theorem}

 \begin{proof}
 First we show that  $\psi \cE$ is cyclic for 
 $\{S_i\otimes I_\cE\}_{i=1}^n$ on  $F^2(H_n)\otimes \cE$.
 Let $x\in F^2(H_n)\otimes \cE$  be such  that 
  ~$x\perp (S_\alpha\otimes I_\cE) T^{-1} \cE$~ for any
 $\alpha\in \FF_n^+$.  Since $T$ is invertible, there exists
  $y\in F^2(H_n)\otimes \cE$
  such that $Ty=x$.  Hence,  we have 
 \begin{equation}\label{s*}
 (S_\alpha^*\otimes I_\cE) Ty\perp T^{-1} \cE, \quad \alpha\in \FF_n^+.
 \end{equation}
 In particular, if $\alpha=e$, then  we have $Ty\perp T^{-1} \cE$. Since $T$ is positive, 
 this implies $y\perp \cE$,
 and  therefore 
 \begin{equation}
 \label{y=}
 y=\left(\sum_{i=1}^n S_i S_i^*\otimes I_\cE\right)y.
 \end{equation}
 Since $T$ is  a multi-Toeplitz  operator and using relations \eqref{y=} 
 and \eqref{s*}, 
 we infer that
 $$
   T^{-1}\cE \perp
 (S_j^*\otimes I_\cE)Ty=(S_j^*\otimes I_\cE)
 T\left(\sum_{i=1}^n S_i S_i^*\otimes I_\cE\right)y= T(S_j^*\otimes I_\cE)y.
 $$
 Therefore, $(S_j^*\otimes I_\cE)y\perp \cE$ 
  for any $j=1,\ldots, n$. Hence,
 $$
 (S_j^*\otimes I_\cE)y=\left(\sum_{i=1}^n S_i S_i^*\otimes I_\cE\right)(S_j^*\otimes I_\cE)y,
 $$
 which together with \eqref{y=} imply 
 $$
 y= \left(\sum_{|\alpha|=2} S_\alpha S_\alpha^*\otimes I_\cE\right)y.
 $$
   By induction, we can prove that
  $$
 y= \left(\sum_{|\alpha|=k} S_\alpha S_\alpha^*\otimes I_\cE\right)y
 $$
  for any $k=1,2,\ldots$.
  This shows that
   $y\perp \left(S_\alpha\otimes I_\cE\right) (1\otimes\cE)$ 
  for any $\alpha\in \FF_n^+$.
  Hence, $y=0$ and $x=0$. This proves that $M_\psi$ is an outer operator.
  
  Now, let us show that  $M_\psi$ is a bounded multi-analytic operator.
  Since  $T$ is a strictly  positive multi-Toeplitz operator, we have
  \begin{equation*}\begin{split}
  \left< T(S_\alpha \otimes I_\cE)\psi h, (S_\alpha \otimes I_\cE) \psi h\right>
  &= \left< T \psi h,  \psi h\right>=\left< TT^{-1} h, T^{-1} k\right>\\
  &=\left< T^{-1} h,  k\right>= \left< N^2 h,  k\right>,
  \end{split}
  \end{equation*}
  for any $h,k\in \cE$.
 On the other hand, if $\alpha, \beta\in \FF_n^+$ and   $\beta>\alpha$, i.e.,
 there exists $\tau\in \FF_n^+$ such that $\beta=\alpha \tau$,  then 
 we have 
 \begin{equation*} \begin{split}
 \left< T(S_\alpha \otimes I_\cE)\psi h, (S_\beta \otimes I_\cE)\psi k\right>
 &=
 \left< (S_{\beta\backslash \alpha}^* \otimes I_\cE)T\psi h, \psi k\right>\\
 &=
 \left< (S_{\beta\backslash \alpha}^* \otimes I_\cE)(1\otimes h), \psi k\right>=0,
 \end{split}
 \end{equation*}
 where $\beta\backslash \alpha\in \FF_n^+$ is uniquely determined by the equation 
 $\beta=\alpha (\beta\backslash \alpha)$. 
 If $\beta<\alpha$, then 
 \begin{equation*} \begin{split}
 \left< T(S_\alpha \otimes I_\cE)\psi h, (S_\beta \otimes I_\cE)\psi k\right>
 &=
 \left< T(S_{\alpha\backslash \beta} \otimes I_\cE)\psi h, \psi k\right>\\
 &=
 \left< (S_{\alpha\backslash \beta} \otimes I_\cE) \psi h,  1\otimes k\right>=0.
 \end{split}
 \end{equation*}
 Since $T$ is multi-analytic, 
   if $\alpha$ and $\beta $ are not comparable, then 
 $$
 \left< T(S_\alpha \otimes I_\cE)\psi h, 
 (S_\beta \otimes I_\cE)\psi k\right>=0.
 $$
  Therefore, if $f:=\sum_{|\alpha|\leq m} e_\alpha \otimes h_\alpha$ and 
  $g:=\sum_{|\beta|\leq p} e_\beta \otimes h_\beta$ are in 
  $F^2(H_n)\otimes \cE$, then we have
  \begin{equation*}\begin{split}
  \left<T^{1/2} M_\psi f, T^{1/2} M_\psi g\right>
  &=\sum_{|\alpha|\leq m, |\beta|\leq p} 
  \left< T(S_\alpha \otimes I_\cE)\psi h_\alpha, (S_\beta \otimes I_\cE)\psi k_\beta\right>\\
  &=\sum_{|\alpha|\leq m, |\beta|\leq p} \delta_{\alpha, \beta}
  \left< N^2 h_\alpha, k_\beta\right>\\
  &=\sum_{|\alpha|\leq m, |\beta|\leq p}\left< (I\otimes N^2)(e_\alpha 
  \otimes h_\alpha), 
  e_\beta\otimes k_\beta\right>\\
  &= \left< (I\otimes N^2)f, g\right>.
  \end{split}
  \end{equation*}
  In particular, by choosing $f=g$, we get
  $$
  \|T^{1/2} M_\psi f\|\leq \|N\| \|f\|.
  $$
  Since $T^{1/2}$ is invertible, it follows that $M_\psi$ can be extended
  to a bounded operator on $F^2(H_n)\otimes \cE$. Therefore, $M_\psi$ 
  is a multi-analytic operator.
  The above computations, show that
  \begin{equation}\label{M*TM}
  M_\psi^* T M_\psi= I\otimes N^2.
  \end{equation}
  
  Since $T$ is a strictly positive operator, we
     infer that $P_\cE T^{-1} |\cE$ is  an invertible operator on $\cE$.
    The equation 
   $
     P_\cE T^{-1}|\cE= N^2
   $
    shows that $N$  and   $I\otimes N$  are invertible.
  Hence, and using relation \eqref{M*TM}, we deduce that $M_\psi$ is injective.
  Define $\Lambda_0: 
  \text{\rm range}~ M_\psi\to F^2(H_n)\otimes \cE$ by setting
  $$
  \Lambda_0(M_\psi f):=f, \quad f\in  F^2(H_n)\otimes \cE.
  $$
  Using relation \eqref{M*TM}, we have 
  $$ 
  \|(I\otimes N)f\|^2=\left< T M_\psi f, M_\psi f\right>\leq \|T\|\|M_\psi f\|^2.
  $$
 On the other hand, since  $I\otimes N$  is  invertible, there exists a constant
  $K>0$ such that
 $$
 \|(I\otimes N)f\|^2\geq K\|f\|^2= K\|\Lambda_0(M_\psi f)\|^2.
 $$
 Combining these inequalities, we infer
 $$
 \|\Lambda_0(M_\psi f)\|\leq \sqrt{\frac{\|T\|} {K}} \|M_\psi f\|
 $$
 for any $f\in F^2(H_n)\otimes \cE$.
 Since $M_\psi$ is outer, $\Lambda_0$ can be extended to a bounded
  operator on 
 $F^2(H_n)\otimes \cE$.
 It is clear that
 \begin{equation}\label{LM}
 \Lambda M_\psi=I.
 \end{equation}
 On the other hand, if $g\in F^2(H_n)\otimes \cE$  and
  $\{f_k\}_{k=1}^\infty$ is a sequence of elements 
 in $F^2(H_n)\otimes \cE$ such that $M_\psi f_k\to g$, as $k\to\infty$, then,
 using \eqref{LM}, we obtain
 $$
 M_\psi \Lambda g= \lim_{k\to\infty} M_\psi \Lambda M_\psi f_k= 
 \lim_{k\to\infty} M_\psi f_k=g.
 $$
 Therefore, $M_\psi \Lambda =I$. Now, we can draw the conclusion that $M_\psi$ 
 is invertible and 
  $\Lambda$ is its inverse, which is also an outer multi-analytic operator.
  Moreover, relation
  \eqref{M*TM} shows that
  $$
  T=(M_\psi^{-1})^* (I\otimes N^2) M_\psi^{-1}.
  $$
 Therefore, $(I\otimes N) M_\psi^{-1}$ is an outer spectral factor for $T$. 
 The proof is complete.
 \end{proof}

 \begin{corollary}\label{to-entro}
 If   $\dim \cE<\infty$ and $T\in B(F^2(H_n)\otimes \cE)$ is  a strictly positive
  multi-Toeplitz operator, then its entropy  $e(T)$ satisfies the equality
  $$
  e(T)= -\ln \det [P_\cE T^{-1} |\cE].
  $$
 \end{corollary}
 
 \begin{proof}
 According to  relations \eqref{DTo} and \eqref{ET}, 
 we  have 
 $$
 e(T)=\ln\det \Delta_T=\ln \det [P_\cE T^{1/2} P_{\cM_T} T^{1/2}|\cE].
 $$
 Therefore, 
 it is enough to prove that
 \begin{equation}
 \label{PET}
 P_\cE T^{1/2} P_{\cM_T} T^{1/2}|\cE= [P_\cE T^{-1}|\cE]^{-1}.
 \end{equation}
 According to \eqref{M_T},  $h\in \cM_T$ if and only if 
 $$T^{1/2} h\perp (S_i\otimes I_\cE) (F^2(H_n)\otimes \cE), \quad i=1,\ldots, n.
 $$
 Hence, $h\in \cM_T$ if and only if $T^{1/2} h\in \cE$. Therefore, we have 
 \begin{equation} \label{MT}
 \cM_T= \text{\rm range}~(T^{-1/2}|\cE).
 \end{equation}
 On the other hand, it is well-known that if $X:\cX_1\to \cX_2$
  is an injective operator with closed range, then $X(X^*X)^{-1} X^*$ is equal
   to the orthogonal projection of $\cX_2$ onto the range of $X$. Applying this
    result to  the operator 
    $$
    X:= (T^{-1/2}|\cE):\cE\to F^2(H_n)\otimes \cE,
    $$
    and taking into account relation 
    \eqref{MT}, we obtain
    $$
    P_{\cM_T}= (T^{-1/2}|\cE)[P_\cE T^{-1}|\cE]^{-1}(P_\cE T^{-1/2}).
    $$
    This clearly implies  relation \eqref{PET}. The proof is complete.
 \end{proof}

 \bigskip

\subsection{Operators on Fock spaces and factorizations}
\label{factorizations}

In this section, we provide factorization results for operators in 
$B(\cH_1, F^2(H_n)\otimes \cH_2)$, where $\cH_1$   is a finite dimensional 
Hilbert space, 
and for multi-analytic operators (see Theorem \ref{in-out}), 
generalizing 
classical results from
 \cite{HL1}, \cite{HL2}, \cite{RR},  and \cite{W},  as well as some extensions 
 to Fock spaces 
 (see \cite{Po-multi}, \cite{Po-analytic}, and \cite{Po-structure}). 
A noncommutative multivariable analogue of Robinson's minimum energy
 delay principle (see  \cite{R})
for outer operators on Fock spaces is obtained (see Theorem \ref{outer21}).

A positive definite kernel on   $\FF_n^+$ is a map
$
K:\FF_n^+\times \FF_n^+\to B(\cH)
$
with the property that 
$$
\sum\limits_{i,j=1}^k\langle K(\sigma_i,\sigma_j)h_j,h_i\rangle\ge 0
$$
for any   
 $h_1,\ldots, h_k\in \cH$,  
$\sigma_1,\dots,\sigma_k\in\FF_n^+$, and  $k\in\NN$.
A kernel $K$ on $\FF_n^+$ is called multi-Toeplitz  if 
$$
K(\sigma, \omega):=
\begin{cases} K(\alpha,g_0), & \text{ if  }  
\sigma=\omega \alpha \text { for some }  \alpha\in \FF_n^+;\\ 
K(g_0,\alpha), & \text{ if  }  
\omega=\sigma \alpha  \text { for some }  \alpha\in \FF_n^+;\\ 
 0,  &  \text{  otherwise.  }  
\end{cases}
$$
If $K(g_0, g_0)=I_\cH$, ($g_0$ is the neutral element in $\FF_n^+$),
then the  kernel is called normalized.
  Let $\theta\in B(\cH, F^2(H_n)\otimes \cK)$, i.e.,
   $$
   \theta h= \sum_{\sigma\in \FF_n^+} e_{\tilde\sigma} \otimes \theta_\sigma h\quad  
  \text{  for some } \quad \theta_\sigma\in B(\cH, \cK)
  $$
   with the property that there
    is $c>0$ such that 
$$
\sum_{\sigma\in \FF_n^+}\|\theta_\sigma h\|^2\leq c\|h\|^2\quad  
  \text{  for any } \quad h\in \cH.
  $$
Denote by $\cP(\cH)$ the set of all polynomials in $F^2(H_n)\otimes \cH$.
 Define the linear operator $M_\theta: \cP(\cH)\to F^2(H_n)\otimes \cK$ by
$M_\theta(1\otimes h):=\theta h$ and 
$$
M_\theta(e_\omega\otimes h):= (S_\omega\otimes I_\cK)\theta h \ 
\text{ for any } \ h\in \cH, \omega\in \FF_n^+. 
$$
Since 
$$
M_\theta(S_i\otimes I_\cH)|\cP(\cH)= (S_i\otimes I_\cK)M_\theta |\cP(\cH), 
\quad  i=1,\ldots, n,
$$
 we can view $M_\theta$ as an unbounded generalized multiplier.
In general, $M_\theta$ cannot be extended to a bounded linear  operator from
$F^2(H_n)\otimes \cH$ to $ F^2(H_n)\otimes \cK$.
However,
its matrix representation 
$$
M_\theta:=[M_{\sigma,\omega}], \quad M_{\sigma,\omega}:= P_\cK(S_\sigma^* 
\otimes I_\cK)M_\theta (S_\omega\otimes I_\cH)|_\cH\in B(\cH, \cK)
$$ 
has each column bounded as an operator from $\cH$ to $F^2(H_n)\otimes \cK$.
It makes sense to define the kernel $K_\theta  :\FF_n^+\times \FF_n^+\to B(\cH)$
  by setting 
$$
K_\theta(\sigma, \omega):= \sum_{\alpha\in \FF_n^+}
 M_{\sigma, \alpha}^* M_{\alpha, \omega}, 
$$
where the convergence is in the SO-topology.
 It is easy to see that $K_\theta$ is a
 positive definite multi-Toeplitz kernel which is not normalized in general,
   i.e.,  $K_\theta(g_0,g_0)\neq I_\cH$.
Notice that if $M_\theta$ can be extended to a bounded operator, 
then the operator matrix $  [K_\theta(\sigma, \omega)]_{\sigma, \omega \in 
\FF_n^+}$
represents a multi-Toeplitz operator  on $F^2(H_n)\otimes \cH$
which is equal
to $M_\theta^* M_\theta$.

 In \cite{Po-structure}, we found an inner-outer factorization for
  any bounded  linear operator 
 $\theta\in B(\cH_1 , F^2(H_n)\otimes \cH_2)$ with $K_\theta(g_0,g_0)=I$.
 In what follows, we show that the latter condition can be removed if 
 ~$\dim \cH_1<\infty$~ or, more generally, if the operator
  $K_\theta(g_0,g_0) $ has closed range.
 Moreover, if $\cH_1$ is an arbitrary Hilbert space, we  
  prove the existence of an inner-outer
 factorization for any bounded multi-analytic operator   $M_\theta\in R_n^\infty
 \bar\otimes B(\cH_1, \cH_2)$.    When $\cH_1= \cH_2$, we obtain an explicit form 
   of the inner-outer
 factorization provided in 
 \cite{Po-multi}.   Our proof here is   based on the
  existence of maximal outer factors for arbitrary multi-Toeplitz operators (see Theorem
  \ref{sub-fac}).

 \begin{theorem}\label{in-out}
 Let $\cH_1$ and  $\cH_2$ be Hilbert spaces.
 \begin{enumerate}
 \item[(i)]
  If  $\dim \cH_1<\infty$, then any
  operator 
 $\theta\in B(\cH_1 , F^2(H_n)\otimes \cH_2)$ admits a factorization
 $$
 \theta=M_{\chi} \psi,
 $$
 where
  $\psi\in B(\cH_1, F^2(H_n)\otimes \cH_3)$   is outer  and 
  $ M_{\chi}\in R_n^\infty\bar \otimes B(\cH_3, \cH_2)$ is
    an  inner operator. Moreover, the factorization is uniquely determined up to a
     diagonal  unitary multi-analytic operator.
     
    \item[(ii)] If $\cH_1$ is an arbitrary Hilbert space and 
     $M_\theta\in R_n^\infty\bar{\otimes} B(\cH_1, \cH_2)$ is  a bounded 
 multi-analytic operator, then   there exist  multi-analytic operators
 $M_\psi \in R_n^\infty\bar{\otimes} B(\cH_1, \cH_3)$    and 
 $M_\chi \in R_n^\infty\bar{\otimes} B(\cH_3, \cH_2)$   such that $M_\psi$ is outer, 
 $M_\chi$ is inner, and 
 $$
 M_\theta=M_\chi M_\psi.
 $$
 Moreover, the inner-outer factorization of $M_\theta$ is  
   uniquely determined up to a unitary diagonal 
  multi-analytic operator.
  \end{enumerate}
 \end{theorem}
 
 \begin{proof}
 Consider the representation 
 $$\theta h:= \sum_{\alpha\in \FF_n^+}e_{\tilde\alpha} \otimes  \theta_\alpha h, \quad h\in \cH_1,
 $$
 where $\theta_\alpha\in B(\cH_1, \cH_2)$, and note that 
 $$
 X:= K_\theta (g_0,g_0)= \sum_{\alpha\in \FF_n^+}\theta_\alpha^* 
 \theta_\alpha \in B(\cH_1).
 $$
 If $X$ is an invertible operator, then the multi-Toeplitz kernel
  $K_{\theta X^{-1/2}}$ is normalized, i.e., 
 $K_{\theta X^{-1/2}}(g_0, g_0)=I_{\cH_1}$.
 In this case, we can apply   Theorem 3.3 from  \cite{Po-structure} 
 to the kernel $K_{\theta X^{-1/2}}$
  and get the desired factorization.
 
 Now, assume that $X$ is not invertible. Let $\cN_0:= \ker X$ 
 and $\cH_1= \cN_0\oplus \cN_1$ 
 be the corresponding decomposition.
 For each $\alpha \in \FF_n^+$, denote 
 $$
 \varphi_\alpha:= \theta_\alpha|\cN_1\in B(\cN_1, \cH_2).
 $$
  Note that 
  since $\cH_1$ is finite dimensional, the operator 
 $\sum_{\alpha\in \FF_n^+}\varphi_\alpha^* \varphi_\alpha\in B(\cN_1)$ is invertible.
 Apply now the first part of the proof to the operator 
 $\varphi\in B(\cN_1, F^2(H_n)\otimes \cH_2)$ defined by
 $$
 \varphi h :=\sum_{\alpha\in \FF_n^+}e_{\tilde\alpha} \otimes  \varphi_\alpha h, \qquad h\in \cN_1,
 $$
 and get the factorization
 \begin{equation}\label{psi-prime}
 \varphi=M_{\chi} \varphi',
 \end{equation}
 where $\varphi'\in B(\cN_1, F^2(H_n)\otimes \cH_3)$ is outer and 
 $\chi\in B(\cH_3, F^2(H_n)\otimes \cH_2)$ is inner.
 If $\varphi'$ has the representation 
 $$
 \varphi'h =\sum_{\alpha\in \FF_n^+} e_{\tilde\alpha} \otimes  \varphi_\alpha' h, 
 \quad  h\in \cN_1,
 $$
  with $\varphi_\alpha '\in B(\cN_1, \cH_3)$, define
 the operator
 $\psi\in B(\cH, F^2(H_n)\otimes \cH_3)$ by
 \begin{equation}
 \label{psi0}
 \psi h:=\sum_{\alpha\in \FF_n^+} e_{\tilde\alpha} \otimes \varphi_\alpha'P_{\cN_1} k,
  \quad k\in \cH_1.
 \end{equation}
 Since $\varphi'$ is outer, it follows that $\psi$ is also outer.
 On the other hand,  relations \eqref{psi-prime}, \eqref{psi0},
  and $\theta_\alpha |\cN_0=0$, $\alpha \in \FF_N^+$,
  imply $\theta=M_{\chi} \psi$. For the uniqueness, see the proof of Theorem 3.3 from 
  \cite{Po-structure}.  
  
  Now, let us  prove the second part of the theorem.  
  Let 
  $M_\psi\in R_n^\infty\bar\otimes B(\cH_1, \cH_3)$ be the maximal 
   outer factor of the multi-Toeplitz operator $M_\theta^* M_\theta$. According 
   to Theorem 
   \ref{sub-fac}, we have
   \begin{equation}\label{th-ps}
   M_\theta^* M_\theta= M_\psi^* M_\psi.
   \end{equation}
  Define the operator $Y:F^2(H_n)\otimes \cH_3\to F^2(H_n)\otimes \cH_2$ by
 \begin{equation}\label{bigY}
 Y\left[\sum_{|\sigma|\leq m} (S_\sigma \otimes I_{\cH_3}) \psi h_\sigma\right]:=
 \sum_{|\sigma|\leq m} (S_\sigma \otimes I_{\cH_2}) \theta h_\sigma, \quad h_\sigma\in 
 \cH_1.
 \end{equation}
 Since $M_\theta$ and $M_\psi$ are multi-analytic operators satisfying relation
  \eqref{th-ps}
  and $M_\psi$ is outer, it is clear that $Y$ can be extended to a unique isometry  
  from $F^2(H_n)\otimes \cH_3$ to $F^2(H_n)\otimes \cH_2$. 
  We also have $Y M_\psi= M_\theta$. Due to relation \eqref{bigY}, one can check that
  \begin{equation}\label{ana}
  Y(S_i\otimes I_{\cH_3})= (S_i\otimes I_{\cH_2})Y, \quad i=1,\ldots, n.
  \end{equation}
 Indeed, for any $f\in F^2(H_n)\otimes \cH_1$, we have
 \begin{equation*}\begin{split}
 Y(S_i\otimes I_{\cH_3})M_\psi f &= Y M_\psi (S_i\otimes I_{\cH_3})f= 
 M_\theta (S_i\otimes I_{\cH_1}) f\\
 &= (S_i\otimes I_{\cH_1}) M_\theta f= (S_i\otimes I_{\cH_2})Y M_\psi f.
 \end{split}
 \end{equation*}
 Since $M_\psi$ is outer, relation \eqref{ana} follows.
 According to \cite{Po-analytic}, there exists an inner  multi-analytic operator
 $ M_\chi\in R_n^\infty\bar\otimes B(\cH_3, \cH_2)$ such that $Y=M_\psi$.
 Therefore, we have $M_\theta= M_\chi M_\psi$, which is the desired inner-outer factorization.
 
 To prove the uniqueness, let $M_\theta= M_{\chi '} M_{\psi'}$ be another
  inner-outer factorization with $ M_{\chi'}\in R_n^\infty\bar\otimes B(\cH'_3, \cH_2)$
 inner and $ M_{\psi'} \in R_n^\infty\bar\otimes B(\cH_1, \cH'_3)$ outer.
 Then we have
 $$
 M_\theta^* M_\theta = M_\psi^* M_\psi = M_{\psi'}^* M_{\psi '}.
 $$
 As we did earlier in the proof, we can find an inner multi-analytic operator $Z\in 
 R_n^\infty\bar\otimes B(\cH_3, \cH'_3)$ such that $ZM_\psi = M_{\psi'}$.
 Since $M_{\psi'}$ is outer, we deduce that
 $Z$ is unitary. According to \cite{Po-charact}, $Z=I\otimes  U$ where 
  $U\in B(\cH_3, \cH_3')$ is unitary.
 Now  the equation  $M_{\chi } M_{\psi}=M_{\chi '} M_{\psi'}$  implies $M_{\chi } M_{\psi}=
 M_{\chi '}Z M_{\psi}$.  Since $M_\psi$ is outer, we obtain $M_\chi = M_{\chi'}Z$. 
 This completes the proof.
 \end{proof}

 \begin{remark}\label{closed}
   Theorem {\rm \ref{in-out}} part {\rm (i)} remains true if $\cH_1$ 
 is an arbitrary Hilbert space  provided that the operator 
 \ $ K_\theta(g_0,g_0)$
 has closed  range.
 \end{remark}

Given a Hilbert space $\cE$, let   $\cP_k(\cE):=
\operatornamewithlimits{\oplus}\limits_{|\alpha|\leq k} e_\alpha\otimes \cE$ 
be the set of all polynomials 
 of degree $\leq k$. 
 For any  Hilbert space $\cE$,
denote by $\PP_k$ the orthogonal projection of $F^2(H_n)\otimes \cE$ onto 
$\cP_k(\cE)$.
 We say that $\theta\in B(\cE, F^2(H_n)\otimes \cG)$ is outer 
 up to a constant inner operator on the left if $\theta$ admits
  a factorization $\theta=M_{\chi} \psi$ where
   $\psi\in B(\cE, F^2(H_n)\otimes \cE_1)$   is outer  and 
  $ M_{\chi}\in R_n^\infty\bar \otimes B(\cE_1, \cG)$ is
    a constant inner operator, i.e., $M_\chi=I\otimes V$, where $V\in B(\cE,\cG)$ is an
    isometry.
 Note that if  $M_\theta\in R_n^\infty\bar{\otimes} B(\cE, \cG)$ then 
 $M_\psi\in R_n^\infty\bar{\otimes} B(\cE, \cE_1)$.

 The next theorem is a noncommutative multivariable analogue 
 of Robinson's minimum energy delay principle for outer functions
 (see \cite{R}, \cite{FFG}).

\begin{theorem}\label{outer21} Let $k$ be a fixed nonnegative integer.
\begin{enumerate}
\item[(i)]
 If  $\cE$ be a finite dimensional Hilbert space,
 then a bounded   operator
$\theta\in  B(\cE,F^2(H_n)\otimes \cG)$ is outer up to 
 a  constant inner operator 
 on the left if and only if 
 \begin{equation}\label{MM2}
 \|\PP_kM_\psi p\|\leq  \|\PP_kM_\theta p\|, \quad p\in \cP_k(\cE),
 \end{equation}
 for any   operator 
  $\psi\in  B(\cE, F^2(H_n)\otimes\cY)$ satisfying 
 $K_\psi = K_\theta $.
 \item[(ii)]
 If  $\cE$ is an arbitrary Hilbert space, then 
 a multi-analytic operator
$M_\theta\in R_n^\infty \bar{\otimes} B(\cE, \cG)$ is outer up to 
 a  constant inner operator 
 on the left if and only if  the inequality \eqref{MM2} holds
 for any  multi-analytic operator 
  $M_\psi\in R_n^\infty \bar{\otimes} B(\cE, \cY)$ satisfying 
 $M_\psi^* M_\psi= M_\theta^* M_\theta$.
 \end{enumerate}
\end{theorem}

 \begin{proof}
 Assume that 
 $\theta\in  B(\cE, F^2(H_n) \otimes\cG)$ is outer and 
  $\psi\in  B(\cE, F^2(H_n) \otimes\cY)$
 is  an   operator 
  satisfying 
 $K_\psi= K_\theta$.
 Define $Q:F^2(H_n)\otimes \cG\to F^2(H_n)\otimes \cY$ by
 $$
 Q(\sum_{|\sigma | \leq m}(S_\sigma\otimes I_{\cG})\theta h_\sigma):=
 \sum_{|\sigma|\leq m}(S_\sigma\otimes I_\cY)\psi h_\sigma, \qquad h_\sigma \in \cE.
 $$
 Since $K_\psi= K_\theta$ and $\theta$ is an  outer operator, $Q$
  extends to an isometry  from $F^2(H_n)\otimes \cG$ to 
  $F^2(H_n)\otimes \cY$ such that $QM_\theta|\cP(\cE)=M_\psi|\cP(\cE)$.
   Since 
   $$
   Q(S_i\otimes I_{\cG})=
  (S_i\otimes I_\cY)Q,\qquad i=1,\ldots, n,
  $$
   $Q$ is 
  an inner multi-analytic operator.
  Therefore,
  $M_\psi |\cP(\cE)=M_f M_\theta |\cP(\cE)$, 
 where 
 $M_f\in R_n^\infty \bar{\otimes} B(\cG, \cY)$ is an inner operator.
  Note that
 $$\PP_k M_\theta |\cP(\cE)=\PP_k M_\theta \PP_k |\cP(\cE)\quad \text{ \rm and }\quad 
 \PP_k M_f=\PP_k M_f \PP_k.
 $$
 Consequently, for any $p\in \cP_k( \cE)$, we have
 \begin{equation}\label{PPk}
 \|\PP_k M_\psi p\|=\|\PP_k M_f M_\theta p\|=\|\PP_k M_f\PP_k M_\theta p\|\leq
 \|\PP_k  M_\theta p\|.
 \end{equation}
  Hence, we deduce relation 
 \eqref{MM2}.
 
  Conversely, assume that  
   $\theta\in  B(\cE, F^2(H_n) \otimes\cG)$ is an  
   operator
  such that relation 
 \eqref{MM2} holds.
 According to Theorem \ref{in-out} part (i),   we have
  $\theta=M_\chi \varphi$, where 
  $\varphi\in  B(\cE, F^2(H_n)\otimes \cE_1)$ is outer and 
  $M_\chi\in R_n^\infty \bar{\otimes} B(\cE_1, \cG)$ is inner.
  Since $K_\theta=  
  K_\varphi$, we can apply the first part of the proof to
   the outer operator $\varphi$,
   to infer that
  \begin{equation}\label{PAP}
  \|\PP_kM_\theta p\|\leq  \|\PP_k M_\varphi p\|, \quad p\in \cP_k(\cE).
 \end{equation}
 Combining  relation \eqref{MM2} (when $\psi=\varphi$)  with   relation \eqref{PAP}, 
 we obtain 
 \begin{equation}\label{PMP}
 \|\PP_kM_\theta p\|=  \|\PP_k M_\varphi p\|, \quad p\in \cP_k(\cE).
 \end{equation}
 Since $\PP_k M_\chi=\PP_k M_\chi \PP_k$ and $M_\chi$ is inner, we have
 $$\|\PP_k M_\theta p\|=\|\PP_k M_\chi M_\varphi p\|
 =\|\PP_k M_\chi\PP_k M_\varphi p\|
 \leq \|\PP_k M_\varphi p\|, \quad p\in \cP(\cE).
 $$
 On the other hand, 
 since 
   $\varphi$ is outer, we can use relation \eqref{PMP} to deduce that
    $\PP_k M_\chi\PP_k$ is an isometry
  from $\cP_k(\cE_1)$ to $\cP_k(\cE)$.
  Let 
  $$
  M_\chi=\sum_{\alpha \in \FF_n^+}  R_\alpha\otimes \chi_{(\alpha)},\quad  
  \chi_{(\alpha)}\in B(\cE_1, \cG),
  $$ 
   be
   the Fourier representation of
  the multi-analytic operator $M_\chi$, and let $x\in \cE_1$.
  Note that if $\alpha \in \FF_n^+$, $|\alpha|=k$, then 
  $$
  \PP_k M_\chi(S_\alpha\otimes I_{\cE_1})(1\otimes x)=
  \PP_k (S_\alpha\otimes I_{\cG}) M_\chi(1\otimes x)=
   e_\alpha\otimes \chi_{(0)}x.
  $$
 Since $\PP_k M_\chi\PP_k$ is an isometry  
  from $\cP_k(\cE_1)$ to $\cP_k(\cE)$, we have
 $$
 \|\chi_{(0)}x\|=\|\PP_k M_\chi(S_\alpha\otimes I_{\cE_1})(1\otimes x)\|=
 \|(S_\alpha\otimes I_{\cE_1})(1\otimes x)\|=\|x\|
 $$
 for any $x\in \cE_1$. Thus $\chi_{(0)}$ is an isometry from $\cE_1$ to $\cG$.
 Moreover, since $M_\chi$ is inner, we deduce that
 $$
 \|\chi_{(0)}x\|^2= \|x\|^2=\|M_\chi(1\otimes x)\|^2=
 \sum_{\alpha\in \FF_n^+} \|\chi_{(\alpha)} x\|^2
 $$ 
 for any $x\in \cE_1$.  Hence, $M_\chi=I\otimes \chi_{(0)}$,
  i.e., a constant inner operator. 
  
 To prove part (ii) of the theorem, note that if  
  $M_\theta\in R_n^\infty \bar{\otimes} B(\cE, \cG)$, then 
  $[K_\theta(\sigma, \omega)]_{\sigma, \omega \in \FF_n^+}$  is the matrix representation
  of the multi-Toeplitz operator $M_\theta^* M_\theta$. The proof  of part (ii) is very
   similar to the proof of part (i). The only difference is that, in this case, 
    we have to use  Theorem \ref{in-out} part (ii) for 
    the inner-outer factorization of $M_\theta$.
 \end{proof}
 
 Let us  mention  that,  using Remark \ref{closed}, one can show that  Theorem
  \ref{outer21} part (i) remains true if $\cH_1$ 
 is an arbitrary Hilbert space  and the operator 
 \ $ K_\theta (g_0, g_0)$
 has closed  range.

 \bigskip
 
 \subsection{Prediction entropy for positive  definite  multi-Toeplitz 
 kernels on free semigroups}\label{pred-entro}

  In this section, we 
  define the notion of prediction entropy and
  prove entropy inequalities for positive definite
 multi-Toeplitz kernels
on the free semigroup $\FF_n^+$ and multi-analytic operators
 (see Theorem \ref{integr-gen}).
 We extend the classical result (see \cite{H}) which stated for $H^2(\DD)$ says that 
 if $f\in H^2(\DD)$, then $\ln|f(e^{it})|$ is integrable and 
 $$
 \frac {1}{2\pi} \int_{-\pi}^ \pi \ln|f(e^{it})| ~dt \geq \ln |f(0)|.
 $$
 Next, we give a characterization for the outer operators in 
 $B(\cE, F^2(H_n)\otimes \cE)$ if $\dim \cE<\infty$.
In particular, we find a  
noncommutative  multivariable analogue of the following classical result. 
A function $f\in H^2(\DD)$ is outer if and only if $f(0)\neq 0$ and 
$$
 \frac {1}{2\pi} \int_{-\pi}^ \pi \ln|f(e^{it})| ~dt = \ln |f(0)|.
 $$

 Let  $K:\FF_n^+\times \FF_n^+\to B(\cE)$ be a positive definite multi-Toeplitz kernel (not necessarily
 normalized).  Define the positive operator $\Delta_K\in B(\cE)$ by setting
 \begin{equation}\label{DEK}
 \left< \Delta_Kh, h\right>:= \inf_{h_{g_0}=h,~ h_\sigma\in \cE}
  \sum \left< K(\omega,\sigma)h_\sigma, 
  h_\omega\right>, \quad h\in \cE,
  \end{equation}
  where the sum  is taken over all finitely supported sequences 
  $\{h_\sigma\}_{\sigma\in \FF_n^+}\subset \cE$.
  When  $\dim\cE< \infty$, we define the prediction entropy  of the kernel  $K$ by
  \begin{equation}\label{entK}
  e(K):= \ln \det \Delta_K.
  \end{equation}

Let $\cT:=[T_1~\cdots ~T_n]$, $T_i\in B(\cH)$, and define the kernel
$K_\cT:\FF_n^+\times \FF_n^+\to  B(\cH)$ by setting
$$
K_\cT(\sigma, \omega):=
\begin{cases} T_\alpha, & \text{ if  }  
\omega=\sigma\alpha  \text { for some }  \alpha\in \FF_n^+;\\ 
T_\alpha^*, & \text{ if  }  
\sigma=\omega\alpha  \text { for some }  \alpha\in \FF_n^+;\\ 
 0,  &  \text{  otherwise.  }  
\end{cases}
$$
We proved  in \cite{Po-structure} that the multi-Toeplitz kernel $K_\cT$ is positive
 definite if and only if 
$[T_1~\cdots ~T_n]$ is a row contraction. In this case, if $\dim \cH<\infty$, we can use Theorem
1.3 and Theorem 4.1 from \cite{Po-structure} to deduce that  the entropy of $K_\cT$ satisfies
$$
e(K_\cT)= \ln \det (I-T_1T_n^*-\cdots -T_nT_n^*).
$$

 In what follows,  we provide entropy inequalities for positive  definite
  multi-Toeplitz kernels on free semigroups and
  multi-Toeplitz operators on Fock spaces.

 \begin{theorem}\label{integr-gen} Let $\cE$ be a Hilbert space and 
  $\theta\in  B(\cE, F^2(H_n)\otimes \cE)$
  be a nonzero   operator.
 \begin{enumerate}
 \item[(i)]
 If $\dim\cE<\infty$, 
  then   $\Delta_{K_\theta}\neq 0$ and
 \begin{equation}
 \label{infint-gen}
 \left<\Delta_{K_\theta}h, h\right>
 \geq   \sum_{|\alpha|=k}  \left<\theta_\alpha^* \theta_\alpha h,h\right>,
 \end{equation}
 where 
 $k$ is the smallest nonnegative integer such that
  $\theta_{\alpha_0}\neq 0$ for some $\alpha_0\in \FF_n^+$ with $|\alpha_0|=k$.
 Moreover, the entropy of the multi-Toeplitz  kernel $K_\theta$ 
 satisfies the inequality
 \begin{equation}\label{edeca}
 e(K_\theta)\geq \ln \det \left[\sum_{|\alpha|=k}\theta_{\alpha}^*\theta_{\alpha}\right].
 \end{equation}
  \item[(ii)]
  If $\cE$ is an arbitrary Hilbert space and $M_\theta\in R_n^\infty\bar\otimes B(\cE)$ is 
  a nonzero multi-analytic operator, then inequality \eqref{infint-gen} 
  remains true. If, in addition, $\dim \cE<\infty$, then  inequality \eqref{edeca} holds. 
  \end{enumerate}
  In particular, if $\theta\in F^2(H_n) $, then 
  $ e(K_\theta)>-\infty$ and 
  $e(K_\theta)\geq \ln|\theta(0)|^2$. 
 \end{theorem}
 \begin{proof}
  Consider the Fourier representation 
  $$
  \theta h := \sum_{\alpha\in \FF_n^+} e_{\tilde\alpha}\otimes \theta_\alpha h, \quad h\in \cE,
  $$
  and assume that $\theta(0)\neq 0$ \ (recall that $\theta(0):= \theta_{g_0}\in B(\cE)$).
  Since $\dim \cE<\infty$,  there is an invertible operator
  $X:\ker \theta(0)\to \ker \theta(0)^*$. The operator 
  $\theta^\epsilon(0)\in B(\cE)$ defined by
  $$
  \theta^\epsilon(0):=\theta(0)+ \epsilon X P_{\ker \theta(0)}
  $$
  is invertible for any
  $\epsilon>0$.
 Define $\theta^\epsilon:\cE\to F^2(H_n)\otimes \cE$ by
 \begin{equation}\label{t^e}
 \theta^\epsilon h:= 1\otimes \theta^\epsilon(0)h+ \sum_{|\alpha|\geq 1} 
  e_{\tilde\alpha}\otimes\theta_\alpha h,
 \quad h\in \cE.
 \end{equation}
 Since  $N:= \sum\limits_{\alpha\in \FF_n^+} 
  {\theta^\epsilon _\alpha}^*  \theta^\epsilon_\alpha\in B(\cE)$ is invertible,  
 we can define the bounded operator $\psi: \cE\to F^2(H_n)\otimes \cE$ by setting
 $$
 \psi h:= \sum_{\alpha\in \FF_n^+} e_{\tilde\alpha}
  \otimes \theta^\epsilon_\alpha N^{-1/2}h, 
 \quad h\in \cE.
 $$
 Applying Theorem \ref{in-out} part (i), we find an outer operator 
 $\varphi: \cE\to F^2(H_n)\otimes \cE_3$  such that $K_\psi=K_\varphi$.
 Since 
 $$
 K_\psi(g_0,g_0)=\sum\limits_{\alpha\in \FF_n^+} \psi_\alpha^* \psi_\alpha= I_\cE,
 $$
  we can apply 
  Theorem 4.1 from \cite{Po-structure} to the normalized multi-Toeplitz
   kernel $K_\psi$ and deduce
  the equality
  \begin{equation}\label{inf}
  \inf_{h_{g_0}=h,~ h_\sigma\in \cE} \sum \left< K_\psi(\omega,\sigma)h_\sigma, 
  h_\omega\right>=
  \left< \varphi(0) ^*\varphi(0)h,h\right>, \quad h\in \cE,
  \end{equation}
  where the sum  is taken over all finitely supported sequences 
  $\{h_\sigma\}_{\sigma\in \FF_n^+}\subset \cE$. 
 Since $K_\psi=K_\varphi$, we can apply Theorem \ref{outer21} part (i) when $k=0$, to infer that
 $$
  \left< \varphi(0) ^*\varphi(0)h,h\right>\geq \left< \psi(0) ^*\psi(0)h,h\right>,\quad  
  h\in \cE.
  $$
 This inequality together with  equation \eqref{inf} imply
 \begin{equation}\label{dec}
 \Delta_{K_\psi}\geq \psi(0) ^*\psi(0).
 \end{equation}
 Since $K_\psi(\sigma, \omega)= N^{-1/2} K_\theta(\sigma, \omega)N^{-1/2}$, relation 
 \eqref{dec} implies 
 \begin{equation}\label{decc}
 \Delta_{K_{\theta^\epsilon}}\geq \theta^\epsilon (0) ^*\theta^\epsilon(0)
 \end{equation}
 for any $\epsilon >0$.
 Note that
 $$
 \left<  \Delta_{K_{\theta^\epsilon}}h,h\right> =
 \inf_{p\in \cP(\cE),~ p(0)=0} \|(I\otimes \epsilon
  X P_{\ker \theta(0)} + M_\theta) (h-p)\|^2.
 $$
 Taking  $\epsilon\to 0$ in  inequality \eqref{decc}, we get 
 $$
 \Delta_{K_{\theta}}\geq \theta (0) ^*\theta(0).
 $$
  
  Now, assume that $\theta(0)=0$ and let 
 $k$ is the smallest nonnegative integer such that
  $\theta_{\alpha_0}\neq 0$ for some $\alpha_0\in \FF_n^+$ with $|\alpha_0|=k$.
  Then we have
  $$
  \theta=\sum_{|\alpha|\geq k} e_{\tilde\alpha}\otimes \theta_\alpha=\sum_{|\beta|=k} 
  (R_\beta\otimes I)\Lambda_\beta
  $$
  for some operators  $\Lambda_\beta:\cE\to F^2(H_n)\otimes \cE$. Note that
   $\Lambda_\beta(0)= \theta_\beta$, \ $\beta\in \FF_n^+$.
   Since $\{ R_\beta\}_{|\beta|=k}$ are isometries with orthogonal ranges, we deduce
   $$
   \|M_\theta(h-p)\|^2= \sum_{|\beta|=k}\|M_{\Lambda_\beta} (h-p)\|^2
   $$
   for any $h\in \cE$ and $p\in \cP(\cE)$ with $ p(0)=0$.
   Therefore, applying the first part of the proof to each operator 
   $\Lambda_\beta$,   $\beta\in \FF_n^+$ with $|\beta|=k$, we have
   \begin{equation*}\begin{split}
    \inf_{p\in \cP(\cE),~ p(0)=0}\|M_\theta(h-p)\|^2 &\geq 
    \sum_{|\beta|=k}  \inf_{p\in \cP(\cE),~ p(0)=0} \|M_{\Lambda_\beta} (h-p)\|^2\\
    &\geq
    \sum_{|\beta|=k} \left < \Lambda_\beta(0)^* \Lambda_\beta(0) h, h \right>
   \end{split}
   \end{equation*}
   for any $h\in \cE$.
   Since $\Lambda_{\alpha_0}(0)=\theta_{\alpha_0}\neq 0$ and 
   $$
   \left< \Delta_{K_\theta} h,h\right>=\inf_{p\in \cP(\cE),~ p(0)=0}\|M_\theta(h-p)\|^2,
   $$
   we deduce that $\Delta_{K_\theta}\neq 0$ and relation \eqref{infint-gen} is proved.
   The inequality \eqref{edeca} is now obvious.
   
 To prove  part (ii) of the theorem, let $\cE$ be an arbitrary Hilbert space 
 and let $M_\theta\in R_n^\infty\bar\otimes B(\cE)$ be a nonzero multi-analytic 
 operator such that $\theta(0)\neq 0$.
 Using  the Szeg\"o type result of Theorem \ref{Sze} when $T= M_\theta^* M_\theta $, we obtain
 \begin{equation}\label{Sz-eq}
 \inf_{p\in \cP(\cE), ~p(0)=0} \left<  M_\theta^* M_\theta (h-p), h-p\right>= 
 \left< \varphi(0)^* \varphi(0) h, h\right>,
 \end{equation}
 where $M_\varphi$ is the maximal outer  factor of $T$.
 Moreover, due to Theorem \ref{sub-fac}, we have 
 $ M_\theta^* M_\theta=M_\varphi^*M_\varphi$.
 Now, Theorem \ref{outer21} part (ii) implies 
 \begin{equation}\label{mac}
 \left< \varphi(0)^* \varphi(0) h, h\right>\geq \left< \theta(0)^* \theta(0) h, h\right>
 \end{equation}
 for any $h\in \cE$.
 Combining relations \eqref{Sz-eq} and  \eqref{mac}, we get
 $$\Delta_{ M_\theta^* M_\theta}\geq \theta(0)^* \theta(0) .
 $$
 When $\theta(0)=0$, the proof is the same as that  of  part  (i) of the theorem.
 Noting that $\Delta_{K_\theta}=\Delta_{ M_\theta^* M_\theta}$ and 
  $e(K_\theta)= e( M_\theta^* M_\theta)$, one can easily complete the proof.
 \end{proof}

  Let us remark that 
 if $\cE$ is an arbitrary Hilbert space  and
  $\theta\in  B(\cE, F^2(H_n)\otimes \cE)$ is such that  the operator 
  $ K_\theta (g_0,g_0)$  has closed range, then
  the results of Theorem {\rm \ref{integr-gen}} part {\rm (i)} remain true.

 \begin{remark}\label{Szo-case}
 If $f\in F^2(H_1)$, then 
 \begin{equation} \label{cl}
 e(K_f)= \frac{1}{2\pi} \int_{-\pi}^{\pi} \ln |f(e^{it})|^2 ~dt.
 \end{equation}
 \end{remark}
 \begin{proof}
 Under the  canonical identification of the full Fock space $ F^2(H_1)$ with 
 the Hardy space $H^2(\DD)$, we have
 \begin{equation*}
 \begin{split}
 \inf \|M_f(1-p)\|^2 &= \inf \|f(1-p)\|^2\\
 &= \inf \left\{\frac{1}{2\pi} \int_{-\pi}^{\pi} |1-p(e^{it})|^2 |f(e^{it})|^2 ~dt \right\}\\
 &=\exp\left[ \frac{1}{2\pi} \int_{-\pi}^{\pi} \ln |f(e^{it})|^2 ~dt\right],
 \end{split}
 \end{equation*}
 where the infimum is taken over all polynomials $p\in  F^2(H_1)$  with $p(0)=0$.
  The latter equality is due to Szeg\" o's theorem \cite{H}. According 
  to  relation \eqref{DEK} and 
  \eqref{entK},
  the equality \eqref{cl} follows.
 \end{proof} 
 
 A characterization for the outer operators in 
 $B(\cE, F^2(H_n)\otimes \cE)$, when $\cE$ is finite dimensional, 
 is proved  in what follows.

 \begin{theorem}\label{outer2-gen} Let $k$ be a fixed nonnegative integer, 
 $\cE$ be a finite dimensional Hilbert space, and 
 let $\theta\in  B(\cE, F^2(H_n)\otimes \cE)$. 
    Then the following statements  are equivalent:
  \begin{enumerate}
  \item[(i)]
   $\theta $ is  an outer operator;
   \item[(ii)]  $\theta(0)$ is invertible and if
    $\psi\in B(\cE, F^2(H_n)\otimes \cE)$ is  
     such that $K_\psi=K_\theta$, then 
    $$\|\PP_k M_\psi p\|\leq  \|\PP_k M_\theta p\|, \quad p\in \cP_k(\cE);
    $$
    \item[(iii)] $\theta(0)$ is invertible and
 \begin{equation}
 \label{infouter}
 \inf_{h_\sigma \in \cE, ~ h_e=h} \sum \left< K_\theta (\omega, \sigma)h_\sigma, 
 h_\omega \right>
 = \left< \theta(0)^* \theta(0)h, h\right>, \quad h\in \cE,
 \end{equation}
 where the sum is taken over all finitely supported sequences
  $\{h_\sigma\}_{\sigma\in \FF_n^+}\subset \cE$ such that $h_e=h$;
    \item[(iv)] $\theta(0)$ is invertible and 
    the entropy of the  positive definite multi-Toeplitz  kernel
 $K_\theta $ satisfies the equation
 $$
 e(K_\theta)= \ln \det [\theta(0)^* \theta(0)].
 $$
 \end{enumerate}
 \end{theorem}
 
 \begin{proof}
 Assume that $\theta$ is an  outer operator. First,  we show that the operator 
 $\theta(0):=P_\cE\theta\in B(\cE)$  is invertible. Suppose that there exists 
 $y\in \cE$, $y\neq 0$, such that $y$ is orthogonal to the range of $\theta(0)$.
 Note that, for any $h\in \cE$, $\alpha\in \FF_n^+$, we have
 $$\left< (S_\alpha\otimes I_\cE) \theta h, y\right>=0, \quad \text{ if } \ |\alpha|\geq 1,
 $$
 and 
 $$
 \left< \theta h,y\right> =\left<\theta(0)h, y\right>=0.
 $$
 Hence, $y$ is orthogonal to the linear  span $\bigvee_{\sigma\in \FF_n^+}
  (S_\sigma\otimes I_\cE) \theta\cE$, which contradicts that $\theta$ is outer. Therefore,
  the range of $\theta(0)$ is dense in $\cE$. Since $\dim \cE<\infty$, it is clear that 
  $\theta(0)$ must
   be invertible. Now, the implication (i)$\implies$ (ii) follows from Theorem
    \ref{outer21} part (i).
 Assume (ii) holds.    According to the same theorem, the operator  $\theta$ has a factorization
 \begin{equation}\label{te}
 \theta=M_\chi \psi_0
 \end{equation}
 with $\psi_0\in B(\cE, F^2(H_n)\otimes \cE_1)$ outer and $M_\chi=I_{F^2(H_n)}\otimes V$,
 where $V\in B(\cE_1, \cE)$ is an isometry. This shows that $\dim\cE_1\leq \dim \cE$.
 Since $\psi_0$ is outer, the range of $\psi_0(0)$ is dense in $\cE_1$. On the
  other hand, since 
 $\cE$ is finite dimensional, $V$ is an isometry, and $\theta(0)= V\psi_0(0)$ 
 is invertible,
 it follows that $\dim\cE_1= \dim \cE$. Using relation \eqref {te},  we deduce that $\theta$ is outer.
 Therefore (ii)$\Leftrightarrow$ (i).
 
 Assume now that  (i) holds. Since the operator
 $X:= \sum\limits_{\alpha\in \FF_n^+} \theta_\alpha^* \theta_\alpha\in B(\cE)$ 
 is invertible,
 $\theta X^{-1/2}\in B(\cE, F^2(H_n)\otimes \cE)$ is outer and $K_{\theta X^{-1/2}} (g_0, g_0)=I$. 
 Consequently, applying Theorem 4.1 from \cite{Po-structure}, we get
 $$
 \inf_{h_{g_0}=h,~ h_\sigma\in \cE} \sum\left< K_{\theta X^{-1/2}} 
 (\omega, \sigma) h_\sigma, h_\omega\right>=
 \left< \theta(0)^* \theta(0) X^{-1/2}h, X^{-1/2}h\right>, \quad h\in \cE,
 $$
 where the sum is  taken over all finitely supported sequences 
  $\{h_\sigma\}_{\sigma\in \FF_n^+}\subset \cE$.
 Hence, (iii) follows.
 Assume now that (iii) holds. Let $\varphi:\cE\to F^2(H_n)\otimes \cE_1$ be 
 the maximal outer factor of the multi-Toeplitz  kernel $K_{\theta X^{-1/2}}$. According to 
 Theorem 3.3 from \cite{Po-structure}, we have $K_{\theta X^{-1/2}}=K_\varphi$.
 Using  relation \eqref{infouter} and Theorem 4.1 from \cite{Po-structure},   
   we obtain
 \begin{equation*}
 \begin{split}
 \left< \theta(0)^* \theta(0) X^{-1/2}h, X^{-1/2}h\right>&=
 \inf_{h_{g_0}=h,~ h_\sigma\in \cE} \sum\left< K_{\theta X^{-1/2}} 
 (\omega, \sigma) h_\sigma, h_\omega\right>\\
 &= \left< \varphi(0)^* \varphi(0) h, h\right>
 \end{split}
 \end{equation*}
 for any $h\in \cE$.
 Now, as in the proof of Theorem \ref{outer21} (when $k=0$), we obtain the factorization
 $$
 \theta X^{-1/2}= M_\chi \varphi,
 $$
 where $M_\chi=I_{F^2(H_n)} \otimes W$ and $W\in B(\cE_1, \cE)$ is an isometry.
 Since $\theta$ and $X$ are invertible and $\theta(0)X^{-1/2}= W\varphi(0)$,
  we must have $W\cE_1=\cE$. Thus,  $M_\chi$ is a  unitary operator. 
  Since $\varphi$ is outer,  it is clear that
   $\theta= M_\chi \varphi X^{1/2}$ is outer, i.e., (i) holds.
   Therefore (iii)$\Leftrightarrow$ (i).
   
   To prove the equivalence  (iii) $\Leftrightarrow$ (iv) is enough to show that if
   $\theta(0)$ is invertible, then
   $\Delta_{K_\theta}=\theta(0)^* \theta(0)$ if and only if 
   \begin{equation}\label{det}
   \det\Delta_{K_\theta}=\ \det [\theta(0)^* \theta(0)].
   \end{equation}
   Indeed,  according to Theorem \ref{integr-gen}, we have
   $$
   \Delta_{K_\theta}\geq \theta(0)^* \theta(0).
   $$
 Hence, it follows that there is a contraction $C\in B(\cE)$ such that 
 \begin{equation} \label{DCCD}
 \theta(0)^* \theta(0)=\Delta_{K_\theta}^{1/2} C^* C \Delta_{K_\theta}^{1/2},
 \end{equation}
 whence we deduce that
 \begin{equation}\label{det2}
  \det [\theta(0)^* \theta(0)]= \det \Delta_{K_\theta}  \det(C^* C).
 \end{equation}
 Since $\theta(0)$ is  an invertible operator,
  one can see that relations  \eqref{det}  and \eqref{det2}
  imply  that $\det (C^* C)=1$. Hence, and using the fact that $C$ is a contraction, 
  we get  $C^* C=I$.
  Thus, relation \eqref{DCCD} becomes  $\theta(0)^* \theta(0)=\Delta_{K_\theta}$, 
  which proves our claim.
   The proof is complete.
 \end{proof}

 \begin{corollary}\label{outer2} 
 If $f\in F^2(H_n)$, then $f$ is outer if and only if $f(0)\neq 0$ and 
 $$
 e(K_f)=\ln |f(0)|^2.
 $$
 \end{corollary}

 Note that in the particular case when $n=1$ and  $\cE=\CC$,  Theorem \ref{outer2-gen} and 
 Remark \ref{Szo-case} imply 
   the following classical result.  A function $f\in H^2(\DD)$ 
 is outer if and only if  $f(0)\neq 0$ and 
 $$
 \frac {1} {2\pi} \int_{-\pi}^{\pi} \ln |f(e^{it})| ~dt = \ln |f(0)|.
 $$

 \bigskip
 \subsection{Extreme points of the unit ball of  
 $F_n^\infty$}
 \label{extreme points}

 It is well-known  \cite{H} that a function $f\in H^\infty(\DD)$ is an extreme point of 
 the unit ball of $H^\infty(\DD)$ if and only if
 $$
 \frac{1}{2\pi} \int_{-\pi}^\pi \ln(1-|f(e^{it})|^2)~dt=-\infty.
 $$
 In what follows, we present some results concerning the extreme points
  of the unit ball of the 
 noncommutative  analytic Toeplitz algebra $F_n^\infty$ and 
 $F_n^\infty \bar\otimes B(\cH)$, where $\dim\cH<\infty$.
 In particular, we prove that if $\varphi\in F_n^\infty$, $\|\varphi\|\leq 1$, 
 and the entropy $E(\varphi)=-\infty$, then $\varphi$ is an extreme point 
 of the unit ball of 
 $F_n^\infty$. For the converse, a weaker form is provided.

 Let $[V_1'~\cdots ~V_n']$, \ $V_i'\in B(\cK')$,  be isometries with orthogonal
  ranges on a Hilbert space $\cK'$ and let $B:F^2(H_n)\otimes \cE\to \cK'$ 
  be a contractive generalized multiplier with respect to $\{S_i\otimes I_\cE\}_{i=1}^n$ and 
  $\{V_i'\}_{i=1}^n$, i.e., 
  $$
  B(S_i\otimes I_\cE)=V_i'B, \quad i=1,\ldots n.
  $$
 Note  that
  \begin{equation*}\begin{split}
  (S_i^*\otimes I_\cE)(I_{F^2(H_n)\otimes \cE} -B^*B) (S_j\otimes I_\cE)&=
  \delta_{ij} I_{F^2(H_n)\otimes \cE}- B^* {V_i'}^* V_j B\\
  &= \delta_{ij} (I_{F^2(H_n)\otimes \cE}-B^*B)
  \end{split}
  \end{equation*}
  for any $i,j=1,\ldots, n$.
 Therefore, $D_B^2$ is a multi-Toeplitz operator on $F^2(H_n)\otimes \cE$.
 If $\dim \cE<\infty$, we define the prediction entropy
  of the generalized multiplier  $B$ by setting
 \begin{equation}\label{EB}
 E(B):= e(D_B^2).
 \end{equation}

  In   particular, if $\cK':= F^2(H_n)\otimes\cK$ and $V_i':= S_i \otimes I_\cK$, 
  ~$i=1,\ldots, n$, then
  the generalized multipliers $B$  are  multi-analytic operators on Fock spaces. 
Let
 $\Theta: F^2(H_n)\otimes \cE\to F^2(H_n)\otimes \cK $ be  a  multi-analytic
  operator,
 i.e., 
  $\Theta\in R_n^\infty\overline{\otimes} B(\cE, \cK)$.
  If $\|\Theta\|\leq 1$ and   $\dim \cE<\infty$,  then the
   prediction entropy of
   the multi-analytic operator $\Theta$   satisfies the equation
   $$
   E(\Theta)= \ln \det \Delta(\Theta),
   $$
   where 
 \begin{equation}\label{DTheta}
 \left< \Delta(\Theta) x,x\right>:=\inf \{ \left< (I-\Theta^*\Theta) (x-p), x-p\right>:
 \ p\in F^2(H_n)\otimes \cE,\  
 p(0)=0\}
 \end{equation}
 for any $x\in \cE$.
 Using Szeg\" o's theorem in the particular case
  when $n=1$ and  $\cE= \cK=\CC$, we have 
 $$
 E(f)= \frac{1} {2\pi} \int_{-\pi}^{\pi} \ln (1-|f(e^{it})|^2)~dt,
 $$
  which is the classical 
 definition of the entropy
 of  $f\in H^\infty(\DD)$ with $\|f\|\leq 1$.

  Since $R_n^\infty=U^* F_n^\infty U$, where $U$ is the flipping operator, 
  all the results of this section are true
  for both algebras $F_n^\infty$ and $R_n^\infty$.

 \begin{theorem}\label{extreme}
 Let $M_\theta\in R_n^\infty\bar{\otimes} B(\cE)$  be  such that
     $\|M_\theta\|\leq 1$.
 \begin{enumerate}
 \item[(i)]
  If    $\Delta(M_\theta)=0$, then 
   $M_\theta $ is an extreme point 
  of the unit ball of $R_n^\infty \bar\otimes B(\cE)$.
   \item[(ii)] Assume that there is a multi-analytic operator
   $M_\varphi\in R_n^\infty \bar\otimes B(\cE)$, $ M_\varphi\neq 0$, 
   such that its range is orthogonal to the range of $M_\theta$.
   If $\dim \cE<\infty$ and  the entropy $E(M_\theta)>-\infty$, then  $M_\theta$ is not
   an extreme point 
  of the unit ball of $R_n^\infty\bar{\otimes} B(\cE)$.  
   \end{enumerate}
 \end{theorem}
 
 \begin{proof}
 Assume that $\Delta(M_\theta)=0$, and let $G\in R_n^\infty\bar{\otimes} B(\cE)$ be such that 
 $$
 \|M_\theta+ G\|\leq 1 \quad \text{ and } \quad \|M_\theta- G\|\leq 1.
 $$
 Hence, we get 
 $$
 M_\theta^* M_\theta +G^* G\leq I,
 $$ 
 whence
 \begin{equation*}
 \begin{split}
 \inf_{p\in \cP(\cE), ~p(0)=0}\left< G^*G (h-p), h-p\right> &\leq
 \inf_{p\in \cP(\cE), ~p(0)=0}\left< (I-M_\theta^* M_\theta) (h-p), h-p\right> \\
 &= \left<\Delta(M_\theta)h,h\right>=0,
 \end{split}
 \end{equation*}
 for any $h\in \cE$.
 Hence, and using Theorem \ref{integr-gen}, we obtain
 $$
 0=\inf_{p\in \cP(\cE),~ p(0)=0}\left< G^*G (h-p), h-p\right>\geq
  \left< G(0)^* G(0)h, h\right>
 $$
 for any $h\in \cE$. This implies $G(0)=0$.
  Now, assume that $G\neq 0$ and let $k$ be the smallest nonnegative integer such
  $G_{\alpha_0}\neq 0$ for some 
    $\alpha_0\in \FF_n^+$ with  $|\alpha_0|=k$. Using again Theorem \ref{integr-gen},
     we get
    $$
 0=\inf_{p\in \cP(\cE), ~p(0)=0}\left< G^*G (h-p), h-p\right>\geq
  \left< G_{\alpha_0}^* G_{\alpha_0}h, h\right>
 $$
 for any $h\in \cE$, which implies $G_{\alpha_0}=0$, a contradiction.
  Therefore, $G=0$ and consequently
 $M_\theta $ is an extreme point 
  of the unit ball of $R_n^\infty \bar\otimes B(\cE)$.
 
 We prove now  the second part of the theorem.
 Let
 $M_\theta\in R_n^\infty \bar{\otimes} B(\cE)$ and assume that $\dim \cE<\infty$
  and  $\|M_\theta\|\leq 1$.
  Since  the entropy of $M_\theta$ satisfies  
  the equation $E(M_\theta)=e(I-M_\theta^* M_\theta)$, we can apply Corollary \ref{entro-fact}
   to the multi-Toeplitz operator $I-M_\theta^* M_\theta$ in order to find an  outer  
    multi-analytic operator
   $M_\psi\in R_n^\infty \bar{\otimes} B(\cE)$ such that
   \begin{equation} \label{MMMM}
   M_\theta^* M_\theta + M_\psi^* M_\psi\leq I.
   \end{equation}
 Define the operators  
 $$
 X:= M_\theta+ M_\varphi M_\psi \quad \text { and  }\quad
  Y:= M_\theta- M_\varphi M_\psi.
  $$
 Since  the operators 
  $M_\theta$ and $ M_\varphi $
    have orthogonal ranges, we can assume that $\|M_\varphi\|\leq 1$.
     Using  relation \eqref{MMMM}, we obtain
    $$
    X^* X\leq  M_\theta^* M_\theta + M_\psi^*M_\varphi^*M_\varphi  M_\psi \leq 1.
    $$
     Similarly,  we get  $\|Y\|\leq 1$.
    Since 
    $
    M_\theta= {\frac {1} {2}}(X+ Y)
    $
    and $X\neq Y$, 
    we deduce that  $M_\theta$ is not
   an extreme point 
  of the unit ball of $R_n^\infty\bar{\otimes} B(\cE)$.
 This completes the proof. 
 \end{proof}
 
 Let us remark that  
 if $\cE=\CC$, then $\Delta(M_\theta)\neq 0$ if and only if the
  entropy $E(M_\theta)>-\infty$.

 \begin{corollary}\label{implic} If $\varphi\in F_n^\infty$ 
 (resp.~$\cA_n$, the noncommutative disc algebra),  $\|\phi\|\leq 1$, 
   and $E(M_\varphi)=-\infty$, 
 then  $\varphi$ is an extreme point 
  of the unit ball of $F_n^\infty$   (resp. $\cA_n)$.
 \end{corollary}
 
 Whether or not  the converse of this corollary  is  true remains
   an open problem. 
 For the time being, according to Theorem \ref{extreme}, 
 we have a weaker form, namely, 
 if $\varphi\in F_n^\infty$, $\|\varphi\|\leq 1$, and there is $\psi\in F_n^\infty$
  such that   $\varphi $   and $\psi$ have orthogonal ranges, then $\varphi$ is 
   an extreme point 
  of the unit ball of $F_n^\infty$  if and only if the entropy $E(\varphi)=-\infty$.
  
  In \cite{ArPo}, we prove that $H^\infty(\DD)$ can be completely isometrically
   embedded
   into the noncommutative analytic Toeplitz algebra $F_n^\infty$ by the mapping
   $f\mapsto f(S_1)$, \ $f\in H^\infty(\DD)$.
   Surprisingly, under this embedding, the extreme points of the unit ball
    of  $H^\infty(\DD)$ remain so in the unit ball of $F_n^\infty$.

 \begin{theorem}\label{H-F}
 A function $f\in H^\infty(\DD)$ is an extreme point of the unit ball of $H^\infty(\DD)$
 if and only  if  $f(S_1)  $ is an extreme point of the unit ball 
 of $F_n^\infty$.
 \end{theorem}

\begin{proof}
If $f\in H^\infty(\DD)$ is an extreme point of the unit ball of $H^\infty(\DD)$,
 then according to \cite{H},  
 $$
 E(f)= \frac{1} {2\pi} \int_{-\pi}^{\pi} \ln (1-|f(e^{it})|^2)~dt=-\infty.
 $$
Accorging to relation \eqref{DTheta} and the remark that follows, it is clear that
  $E(f(S_1))\leq E(f)$. Therefore,   we have $E(f(S_1))=-\infty$.
   From Corollary
\ref{implic}, it follows that $f(S_1)$ is an extreme point 
of the unit ball of $F_n^\infty$. 

Now, assume that $f\in H^\infty(\DD)$  is not an extreme point
 of the unit ball of $H^\infty(\DD)$. Then
there exist $\varphi, \psi\in (H^\infty(\DD))_1$, \  $\varphi \neq \psi$ such that
$f=\frac{\varphi + \psi} {2}$.
Since the map $f\mapsto f(S_1)$ is a complete isometry of $H^\infty(\DD)$  into 
$F_n^\infty$, and 
   $f(S_1)=\frac{\varphi(S_1) + \psi(S_1)} {2}$,
 we deduce that $f(S_1)$
  is not an extreme point 
of the unit ball of $F_n^\infty$. This completes the proof.
\end{proof}

This result shows that there are extreme  points of the unit ball 
$(F_n^\infty)_1$ 
which are not inner operators. Indeed,  if $f\in H^\infty(\DD)$ is 
an extreme point
 of the unit ball of $H^\infty(\DD)$ such that $f$  
 is not an inner function (see \cite{H}), then, according to Theorem \ref{H-F}, 
 $f(S_1)$  and $S_i f(S_1)$, $ i=2,\ldots, n$, are  also
  extreme points
 of the unit ball of $F_n^\infty$ which are  not  inner operators.
 For some other examples,  take   $\varphi(S_2,\ldots, S_n)\in F_n^\infty$  to be 
  inner  operator  such 
 that $\varphi(0)=0$, and note that
   $ \varphi(S_2,\ldots, S_n)f(S_1)$ has the required property.

\clearpage

 \section{Noncommutative commutant lifting theorem: 
 geometric structure and  maximal entropy solution}

 We obtain several geometric 
 characterizations for the  multivariable central 
 intertwining lifting, a maximum principle, and a permanence principle for the 
 noncommutative commutant lifting theorem.
 Under certain natural conditions, we find explicit forms for 
 the maximal entropy
 solution of this multivariable commutant lifting theorem, and concrete
  formulas for its entropy.

 \subsection{Multivariable intertwining liftings and geometric structure}\label{geom}  
In this section, we obtain some results concerning the geometric
 structure of the
 intertwining liftings in the
noncommutative commutant lifting theorem  
(see \cite{Po-isometric} and  \cite{Po-intert}).
It is shown that there is a one-to-one correspondence between the set of all 
intertwining liftings with tolerance $t>0$ and certain families of  contractions
$\{C_k\}_{k=1}^\infty$ and $\{\Lambda_\alpha\}_{\alpha\in \FF_n^+}$ (see Theorem \ref{main2}
and Theorem \ref{geo}). 
The geometric structure of the multivariable  intertwining liftings will
 play an important role in our investigation.
       
Let us recall from \cite{Po-models}, \cite{Po-isometric}, and  \cite{Po-charact} a few results concerning 
the noncommutative dilation theory for sequences of operators (see
\cite{SzF-book} for the classical case $n=1$).
     A sequence of operators $\cT:=(T_1,\ldots, T_n)$, $T_i\in B(\cH)$,
        is called
 row contraction if 
 $$
 T_1T_1^*+\cdots +T_nT_n^*\le I_\cH.
 $$
We say that a sequence of isometries   $\cV:=(V_1,\ldots, V_n)$,\ $V_i\in B(\cK)$,  
  is a  minimal isometric dilation (m.i.d.) of $\cT$  on a Hilbert space 
$\cK\supseteq \cH$ if
 the following properties are satisfied:
\begin{enumerate}
\item[(i)] $V_i^* V_j=0$   ~for all $i\neq j$, \ $i,j\in\{1,\ldots, n\}$;
\item[(ii)]   $V^*_j|\cH=T^*_j$  ~for all  $j=1,\ldots,n$;
\item[(iii)] $\cK=\bigvee\limits_{\alpha\in\FF^+_n}V_\alpha\cH$.
\end{enumerate} 
If $\cV$ satisfies only the condition (i) and $P_\cH V_i=T_i P_\cH, \ i=1,\dots, n$, 
then $\cV$ is called isometric lifting of $\cT$.
The minimal isometric dilation of $\cT$ is an
 isometric lifting and is uniquely determined up to an isomorphism 
 \cite{Po-isometric}. 

Let us consider a canonical  realization of it on Fock spaces. 
 For convenience of notation, we will sometimes identify the $n$-tuple 
 $(T_1,\ldots, T_n)$ 
 with the row operator  $[T_1~\cdots~T_n]$.
Define  the operator $D_\cT:\oplus_{j=1}^n\cH\to \oplus_{j=1}^n\cH$ by  setting
$D_\cT:=(I_{\oplus_{j=1}^n\cH}-\cT^*\cT)^{1/2}$, and set
 $\cD:=\overline{D_\cT(\oplus_{j=1}^n\cH)}$,
 where $\otimes_{j=1}^n \cH$ denotes the direct sum of $n$ copies of $\cH$.
Let $D_i:\cH\to 1\otimes \cD\subset F^2(H_n)\otimes \cD$ be defined by
$$
D_ih:= 1\otimes D_\cT(\underbrace{0,\ldots,0}_{{i-1}\mbox{ \scriptsize
times}},h,0,\ldots).
$$
Consider the Hilbert space $\cK:=\cH\oplus [F^2(H_n)\otimes \cD]$ and 
define $V_i:\cK\to\cK$ by
\begin{equation}\label{dil}
V_i(h\oplus \xi):= T_ih \oplus [D_ih +(S_i\otimes 
I_\cD)\xi]
\end{equation}
for any $h\in \cH, ~\xi\in F^2(H_n)\otimes\cD$.
Note that 
\begin{equation}\label{dilmatr}
V_i=\left[\begin{matrix} T_i& 0\\ 
D_i& S_i\otimes I_\cD
\end{matrix}\right]
\end{equation}
with respect to the decomposition $\cK=\cH\oplus [F^2(H_n)\otimes \cD]$.
In \cite{Po-isometric},  we proved  that  
 $\cV:=[V_1~\cdots~ V_n]$  is the   minimal isometric dilation of $\cT$.
Let $\cH_0:=\cH$ and 
$$
\cH_k:= \cH_{k-1}\bigvee(\bigvee_{|\alpha|=1}V_\alpha \cH_{k-1})\quad 
\text{ if } k\geq 2. 
$$
Note that
 ~$\cK=\bigvee_{k=0}^\infty \cH_k$, \ 
 $\cH_k\subset \cH_{k+1}$, and all subspaces $\cH_k$ are invariant 
 under each ~$V_i^*$, \ $ i=1,\dots,n$.
On the other hand, we have $\cH_1=\cH\oplus \cD$ and 
$$
\cH_k=\cH\oplus \bigoplus_{|\alpha|\leq k-1} e_\alpha\otimes \cD \quad 
\text{ if }
 k\geq 2. 
$$
 Denote $\cV_0:=\cT$ and  
$\cV_k:=[V_{1,k}~\cdots~ V_{n,k}]$ if $k\geq 1$, 
where $V_{i,k}:=P_{\cH_k} V_i|_{\cH_k}$,  for any  $i=1,\dots,n$, and $P_{\cH_k}$ 
is the orthogonal projection from $\cK$ onto $\cH_k$.
Note that  the operators 
$V_{i,k}, ~i=1,\dots,n$, are partial isometries with orthogonal final spaces
 and initial space $\cH_{k-1}$.
It is easy to see that $\cV$ is also the minimal isometric dilation 
of  $\cV_k$,\ 
 $V_i^*|_{\cH_k}= V_{i,k}^*$, and 
 $$V_i^*=\text {\rm SOT}-\lim_{k\to\infty}
V_{i,k}^*P_{\cH_k}
$$
 for any $i=1,\ldots, n$. 
On the other hand, let us mention that   $\cV_{k+1}$ is
 the one-step dilation of
 $\cV_k$, i.e., 
$\cH_{k+1}= \cH_k\oplus D_{\cV_k}(\oplus_{j=1}^n \cH_k)$, and, for each
 $i=1,\dots,n$,
$$
V_{i,k+1}(x\oplus y)=V_{i,k}x\oplus D_{\cV_k}
(\underbrace{0,\ldots,0}_{{i-1}\mbox{ \scriptsize times}},x,0,\ldots)
$$
for any $x\in \cH_k$ and $y\in D_{\cV_k}(\oplus_{j=1}^n \cH_k)$.

Let   $\cT':=[T_1'~\cdots~ T_n']$, ~$T_i'\in B(\cH')$,  be  another  row contraction
and   let   $\cV':=[V_1'~\cdots~ V_n']$  
 be the   minimal isometric dilation of $\cT'$ on the  Hilbert space
  $\cK':= \cH'\oplus [F^2(H_n)\otimes \cD_{\cT'}]$.
Let  $A\in B(\cH, \cH')$ be an operator  
 satisfying $AT_i=T_i' A$, for any $ i=1,\dots ,n$.
 
  An  intertwining lifting of $A$ is an  
  operator $B\in B(\cK,\cK')$ satisfying 
$BV_i=V_i' B$, for any $ i=1,\dots ,n,$ and $P_{\cH'} B=AP_\cH$.
Let $\cV'_k:=[V'_{1,k}~\cdots ~ V'_{n,k}]$, $k\geq 0$, 
$V'_{i,k}:= P_{\cH_k} V_i' |\cH_k$,
 be the $k$-step dilation  of $\cT'$. An operator $A_k\in B(\cH_k, \cH_k')$ 
 is a $k$-step intertwining lifting of $A_0:=A$, if 
 $$V_{i,k}' A_k= A_k V_{i,k}\quad i=1,\ldots, n,
 $$
 and $P_{\cH'} A_k = A P_\cH|\cH_k$.  A sequence $\{A_k\}_{k=0}^\infty$ of 
 $k$-step intertwining  liftings of $A$ is compatible if 
 $$
 P_{\cH_k'} A_{k+1}= A_k P_{\cH_k}|\cH_{k+1}, \quad k=0,1,\ldots.
 $$
  Note that the sequence 
 $\{A_k\}_{k=0}^\infty$ 
 is compatible if $A_{k+1}$ is a one-step intertwining lifting of $A_k$ for
  any $k=0,1,\ldots$.
  According to  \cite{Po-intert} (see \cite{FF-book} for the classical case $n=1$), 
  there is a one-to-one correspodence between the set of all multivariable
  intertwining liftings $B$ of $A$
  such that $\|B\|\leq 1$ and the set of all compatible sequences $\{A_k\}_{k=0}^\infty$,
  $\|A_k\|\leq 1$,  of 
 $k$-step intertwining  liftings of $A$. Moreover the correspondence is given by
 $$
 B=\text{\rm SOT}-\lim_{k\to \infty} A_k P_{\cH_k}
 \quad \text{
   and  } \quad
 A_k= P_{\cH_k} B |\cH_k, \quad  k=0,1,\ldots.
 $$
 The
 noncommutative commutant lifting theorem (see \cite{Po-isometric},  \cite{Po-multi},
 \cite{Po-intert}) states
  that there always exists an  intertwining lifting $B$ of $A$ 
   with $\|B\|=\|A\|$ (see \cite{SzF-book}, \cite{FF-book}
    for the classical case).

Define the subspaces
 \begin{equation}
 \label{nmd}
 \begin{split}
 \cN:&=
 \left\{\sum_{i=1}^n D_{A} T_i h_i\oplus D_\cT(\oplus_{i=1}^nh_i):\
  \oplus_{i=1}^nh_i\in \oplus_{i=1}^n\cH \right\}^{-},\\
 \cM:&= [\cD_{A}\oplus \cD_\cT]\ominus \cN.
 \end{split}
\end{equation}
 and  the operator $W:\cN\to \cD_{\cT'}\oplus(\oplus_{i=1}^n \cD_A)$ 
 by setting
 \begin{equation}\label{W}
 W\left(\sum_{i=1}^n D_A T_i h_i\oplus D_\cT(\oplus_{i=1}^nh_i\right):=D_{\cT'}
 (\oplus_{i=1}^n Ah_i)\oplus (\oplus_{i=1}^n D_Ah_i).
 \end{equation}
 Since $AT_i=T_i' A,\ i=1,\dots ,n$, we have
 \begin{equation*}\begin{split}
 \left\|\sum_{i=1}^n D_{A}T_ih_i\right\|^2&+ \|D_{\cT}(\oplus_{i=1}^n h_i)\|^2\\
 &=  \left\|\sum_{i=1}^n T_ih_i\right\|^2- \left\|\sum_{i=1}^n AT_ih_i\right\|^2+ 
 \|D_{\cT}(\oplus_{i=1}^n h_i)\|^2\\
 &=\|\oplus_{i=1}^n h_i\|^2-\left\|\sum_{i=1}^n AT_ih_i\right\|^2\\
 &=\sum_{i=1}^n( \|D_{A}h_i\|^2 + \|Ah_i\|^2)-\left\|\sum_{i=1}^n T_i' Ah_i\right\|^2\\
 &=\sum_{i=1}^n \|D_{A}h_i\|^2+ \|D_{\cT'}(\oplus_{i=1}^n Ah_i)\|^2.
 \end{split}
\end{equation*}
 Due to  these  calculations,  it is clear that $W$ is an isometry.
  Consider  the subspaces
 \begin{equation*}\begin{split}
 \cN':&=\left\{\sum_{i=1}^n D_{\cT'}(\oplus_{i=1}^nA h_i)\oplus
  (\oplus_{i=1}^nD_Ah_i):\ 
 \oplus_{i=1}^nh_i\in \oplus_{i=1}^n\cH\right\}^{-},\\
 \cM':&= [\cD_{\cT'}\oplus (\oplus_{i=1}^n\cD_A)]\ominus \cN',
 \end{split}
\end{equation*}
 and note that $W:\cN\to \cN'$ is a unitary operator.
 As in \cite{Po-intert}, one can prove that
 \begin{equation}\label{PG}
 \overline{P_{\cM'}(\cD_{\cT'}\oplus \{0\})}=\cM'.
 \end{equation}
 Consider the orthogonal projections:
  \begin{equation*}\begin{split}
 P_{\cD'}&:\cD_{\cT'}\oplus (\oplus_{j=1}^n \cD_A)\to \cD_{\cT'},\qquad 
 P_{\cD'}[d\oplus (\oplus_{j=1}^n d_j)]:= d,\\
 P_{i}&:\cD_{\cT'}\oplus (\oplus_{j=1}^n \cD_A)\to \cD_{A},\qquad 
 P_{i}[d\oplus (\oplus_{j=1}^n d_j)]:= d_i
 \end{split}
\end{equation*}
 for any $i=1,\ldots, n$.
 
 Taking into account the results from \cite{Po-intert} 
 (see \cite{FF-book} for the case $n=1$),  
  we can prove the following.
 
 \begin{theorem}\label{onestep2} Let $\cT:=[T_1~\cdots~ T_n]$, $T_i\in B(\cH)$,  
  and $\cT'=[T_1'~\cdots~ T_n']$, $T_i'\in B(\cH')$,   be row contractions, 
  and let $\cV$ and $\cV'$ 
  be the corresponding minimal isometric dilations, respectively.
 Let $A\in B(\cH,\cH')$ be a contraction such that $AT_i=T_i'A$ for any $i=1,\ldots, n$.
 Then
 there is a one-to-one correspondence between the set of 
 all contractive one-step intertwining liftings $A_1$ of $A$ and the set of all
 contractions  $C:\cM\to \cM'$, and
  any contractive one-step intertwining lifting 
 $A_1:\cH\oplus \cD_\cT\to \cH'\oplus \cD_{\cT'}$ of $A$ is given by
 \begin{equation}\label{one2}
  A_1=
   \left[\begin{matrix}
 I_\cH&0\\
 0& P_{\cD'}(WP_\cN+C P_\cM)
 \end{matrix}\right]
  \left[\begin{matrix}
 AP_\cH\\
  D_A\oplus I_{\cD_\cT}
 \end{matrix}\right], 
\end{equation}
 where $C:\cM\to \cM'$ is a  contraction.

  Moreover, 
 there exists a unitary operator 
 $\Omega:\cD_{A_1}\to (\oplus_{j=1}^n \cD_A)\oplus \cD_C$
 satisfying
  \begin{equation}\label{omega1}
   \Omega D_{A_1}=
   \left[\begin{matrix}
   \left[\begin{matrix} P_1\\\vdots\\ P_n\end{matrix}
    \right](WP_\cN+CP_\cM)\\
   D_C P_\cM 
   \end{matrix}\right] (D_A\oplus I_{\cD}).
 \end{equation}
 \end{theorem}
 \begin{proof}

 Let 
 $A_1:\cH\oplus \cD_\cT\to \cH'\oplus \cD_{\cT'}$ be a contractive  one-step 
 intertwining 
 lifting of $A$.  Since  $\|A_1\|\leq 1$, $ A_1^*|\cH'= A^*$, and taking into account the geometric structure of 
 $2\times 2$  matrix contractions, we deduce that
 \begin{equation}
 \label{Aone}
 A_1=\left[ \begin{matrix} A&0\\
 XD_A& Y\end{matrix}\right],
 \end{equation}
 where $[\begin{matrix} X& Y \end{matrix}]: \cD_A\oplus \cD_T\to \cD_{\cT'}$ is a contraction.
 One the other hand, we must have
 $$
 \left[ \begin{matrix} A&0\\
 XD_A& Y\end{matrix}\right]\left[ \begin{matrix} T_i&0\\
 D_i& 0\end{matrix}\right]= 
 \left[ \begin{matrix} T_i&0\\
 D_i& 0\end{matrix}\right]\left[ \begin{matrix} A&0\\
 XD_A& Y\end{matrix}\right] 
 $$
 for any $i=1,\ldots, n$.
  This  equality holds if and only if 
  $$
  AT_i=T_i'A \ \text{   and  }\  XD_AT_i +YD_i= D_i A
  $$
 for any $i=1,\ldots, n$. The latter relation is equivalent to 
 $$
 [\begin{matrix} X& Y \end{matrix}]\left[\begin{matrix} D_A T_i\\ D_i 
 \end{matrix}\right]=
 D_iA, \quad i=1,\ldots, n,
 $$
 which shows that $[\begin{matrix} X& Y \end{matrix}]|\cN=F_0$, where the operator
  $F_0:\cN\to \cD_{\cT'}$ is defined by
 \begin{equation}\label{Fzero}
 F_0\left[\sum_{i=1}^n D_A T_i h_i\oplus D_\cT(\oplus_{i=1}^nh_i)\right]:=
 D_{\cT'}(\oplus_{i=1}^nAh_i).
 \end{equation} 
  It is well known (see \cite{FF-book}) that the set of all contractions
  satisfying the equation $[\begin{matrix} X& Y \end{matrix}]|\cN=F_0$ is given
 by 
 $$
 [\begin{matrix} X& Y \end{matrix}]=F_0P_\cN + D_{F_0^*} F_1 P_\cM,
 $$
 where $F_1:\cM\to \cD_{F_0^*}$ is an arbitrary contraction. 
 Then any contractive one-step intertwining lifting $A_1$ of $A$ is given by
 \begin{equation}\label{one}
  A_1=
   \left[\begin{matrix}
 I_\cH&0\\
 0& F_0P_\cN+D_{F_0^*} F_1 P_\cM
 \end{matrix}\right]
  \left[\begin{matrix}
 AP_\cH\\
  D_A\oplus I_{\cD_\cT}
 \end{matrix}\right], 
\end{equation}
 where $F_1:\cM\to \cD_{F_0^*}$ is an arbitrary contraction.
 Moreover, there is a one-to-one correspondence between the set of 
 all contractive one-step intertwining liftings $A_1$ of $A$ and
  the set of all
 contractions $F_1:\cM\to \cD_{F_0^*}$.

  Since 
  \begin{equation}\label{f0}
  F_0=P_{\cD'}W,
  \end{equation}
  we have
  \begin{equation*}\begin{split}
  \|D_{F_0^*}h\|^2&= \|h\|^2-\|F_0^*\|^2=
   \|h\|^2-\|W^* P_{\cN'}P_{\cD'}^* h\|^2\\
  &= \|h\|^2-\|P_{\cN'}P_{\cD'}^* h\|^2=\|P_{\cM'}P_{\cD'}^* h\|^2
  \end{split}
  \end{equation*}
  for any $h\in\cD'$.
  Hence, the relation $RD_{F_0^*}=P_{\cM'} P_{\cD'}^*$ 
  defines an isometry from $\cD_{F_0^*}$
  into $\cM'$. Using relation \eqref{PG}, we infer that $R$ is onto.
  Therefore, the operator  $R^*:\cM'\to \cD_{F_0^*}$ is   unitary  and 
\begin{equation}
\label{R*}
D_{F_0^*}R^*=P_{\cD'}|\cM'.
  \end{equation}
 Now, using  relations \eqref{one}, \eqref{f0}, and \eqref{R*}, we get
 \begin{equation*}\begin{split}
 A_1&=
   \left[\begin{matrix}
 I_\cH&0\\
 0& P_{\cD'}(WP_\cN+D_{F_0^*}R^* R F_1 P_\cM)
 \end{matrix}\right]
  \left[\begin{matrix}
 AP_\cH\\
  D_A\oplus I_{\cD_\cT}
 \end{matrix}\right]\\
 &=
   \left[\begin{matrix}
 I_\cH&0\\
 0& P_{\cD'}(WP_\cN+(P_{\cD'} |\cM') RF_1 P_\cM)
 \end{matrix}\right]
  \left[\begin{matrix}
 AP_\cH\\
  D_A\oplus I_{\cD_\cT}
 \end{matrix}\right],
 \end{split}
 \end{equation*} 
 which proves relation \eqref{one2}, where $C:= RF_1:\cM\to \cM'$.
  Now, define the operator 
 $$
 V:\cH\oplus \cD\to \cH'\oplus \left[\cD'\oplus\left(
 \operatornamewithlimits{\oplus}\limits_{i=1}^n \cD_A\right)\right]\oplus \cD_C\simeq
 (\cH'\oplus \cD')\oplus \left(\operatornamewithlimits{\oplus}\limits_{i=1}^n \cD_A
 \oplus \cD_C\right),
 $$
 by setting
 $$
 V:=\left[\begin{matrix}I&0\\0& WP_\cN +CP_\cM\\
 0& D_C P_\cM
 \end{matrix} 
 \right]
 \left[\begin{matrix}AP_\cH\\ D_A\oplus I_{\cD}
 \end{matrix} 
 \right].
 $$
 Note that  $V$ is  a product of two isometries. According to  relation \eqref{one2},
  we have
   $A_1=P_{\cH'\oplus \cD'}V$.
  For each $x\in \cH\oplus \cD$, we get
  \begin{equation*} \begin{split}
  \|D_{A_1}x\|^2&= \|x\|^2-\|A_1x\|^2 =\|Vx\|^2-\|A_1x\|^2\\
  &=\left\|P_{\left(\operatornamewithlimits{\oplus}\limits_{i=1}^n
   \cD_A\right)\oplus \cD_C} Vx\right\|^2.
  \end{split}
  \end{equation*}
 Hence, it is clear that $\Omega$ defines an isometry from $\cD_{A_1}$ to
 $(\operatornamewithlimits{\oplus}\limits_{i=1}^n
   \cD_A)\oplus \cD_C$. 
  On the other hand,
  for any $h_i\in \cH$, we have
  \begin{equation*}\begin{split}
  \Omega D_{A_1}\left(\sum_{i=1}^n T_i h_i\oplus D_{\cT}(\oplus_{i=1}^n h_i)
  \right)
  &= \left[\begin{matrix} P_1\\\vdots\\ P_n\end{matrix} 
 \right] W\left(\sum_{i=1}^n D_AT_i h_i\oplus D_{\cT}(\oplus_{i=1}^n h_i)\right)
 \oplus 0\\
 &=\operatornamewithlimits{\oplus}\limits_{i=1}^n
   D_Ah_i\oplus 0. 
  \end{split}
  \end{equation*}
Hence, the range of $\Omega$ contains  the subspace 
 $\operatornamewithlimits{\oplus}\limits_{i=1}^n
   D_Ah_i\oplus \{0\}$.
   Notice that, if $x\in \cM$, then  $\Omega D_{A_1} 
   y=z\oplus D_Cx$ for some
    $z\in \operatornamewithlimits{\oplus}\limits_{i=1}^n
   \cD_A$. It is clear that
    $\{0\}\oplus\cD_C\cM$ is also in the range of $\Omega$.
    Therefore, $\Omega$ is  unitary.
 The proof is complete. 
 \end{proof}
  
  A repeated application of Theorem \ref{onestep2} provides compatible 
  contractive sequences $\{A_k\}_{k=0}^\infty$ of $k$-step intertwining liftings of $A$.
  Setting $B:=\text{\rm SOT}-\lim\limits_{k\to \infty} A_k P_{\cH_k}$, we get 
  an intertwining lifting $B$ of $A$ satisfying $\|A\|=\|B\|$, which proves 
  the noncommutative 
  commutant lifting theorem.
   For details, see \cite{Po-intert}
  (resp. \cite{FF-book} for the classical case $n=1$).
  
  We should mention that, as in the classical case, 
  the general setting of the noncommutative commutant lifting theorem can be reduced to 
  the case when $\cT:=[T_1~\cdots ~T_n]$ is a row 
  isometry (see Corollary 2.2 in \cite{Po-nehari}).

 Due to this reason, we will assume 
  from now on  that $\cT:=[T_1~\cdots ~T_n]$ 
 is an isometry. 
  Let $\cV_k':=[V_{1,k}'~\cdots~ V_{n,k}']$, $k\geq 0$, 
  $V_{i,k}':= P_{\cH_k} V_i'|\cH_k$,
  be the $k$-step dilation of $ \cT':=[T_1'~\cdots ~T_n']$.
  Let $A\in B(\cH, \cH')$, $t>0$,  and assume that $\|A\|\leq t$.
  Define the defect operator
  $D_{A,t}:= (t^2I-A^*A)^{1/2}$ and the subspace $\cD_{A,t}:= \overline{D_{A,t} \cH}$.
  In the particular case when $t=1$, we use the classical notation $D_{A,1}:=D_A$ and 
  $\cD_{A,1}:=\cD_A$.

  We say  that $A_k:\cH\to \cH_k'$ is
   a  $k$-step intertwining lifting  of $A_0:=A$ with tolerance $t>0$ if
    $\|A_k\|\leq t$,
   \ $V_{i,k}' A_k = A_k T_i$,  for any $i=1,\ldots, n$, and 
  $P_{\cH'}A_k= A$.
  Note that
  $\cD_k':= \cD_{\cV_k'}$ can be identified with 
  $ \operatornamewithlimits{\oplus}\limits_{|\alpha|=k}
    e_\alpha\otimes \cD '$,
  where $\cD ':= \cD_{\cT'}$.
  For each $k\geq 0$,  define the following subspaces:
  \begin{equation}
  \label{nmk}
  \begin{split}
 \cN_k:&=\left\{\sum_{i=1}^n D_{A_k,t} T_i h_i :\ 
 \oplus_{i=1}^nh_i\in \oplus_{i=1}^n\cH\right\}^{-}\subseteq \cD_{A_k,t},\\
 \cM_k:&= \cD_{A_k,t} \ominus \cN_k\subseteq \cD_{A_k,t},
 \end{split}
\end{equation}
 and 
 
 \begin{equation*}\begin{split}
 \cN_k':&=\left\{\sum_{i=1}^n D_k'(\oplus_{i=1}^nA_k h_i)\oplus
  (\oplus_{i=1}^nD_{A_k,t}h_i):\ 
 \oplus_{i=1}^nh_i\in \oplus_{i=1}^n\cH\right\}^{-},\\
 \cM_k':&= [\cD_k'\oplus (\oplus_{i=1}^n\cD_{A_k,t})]\ominus \cN_k',
 \end{split}
\end{equation*}
 where $D_k':= D_{\cV_k'}$. Note that $\cN_0= \cN$, $\cM_0= \cM$, $\cN_0'= \cN'$, and 
  $\cM_0'= \cM'$.
 Consider the orthogonal projections:
  \begin{equation*}\begin{split}
 P_{\cD_k'}&:\cD_k'\oplus (\oplus_{j=1}^n \cD_{A_k,t})\to \cD_k',\qquad 
 P_{\cD_k'}[d\oplus (\oplus_{j=1}^n d_j)]:= d,\\
 P_{i,k}&:\cD_k'\oplus (\oplus_{j=1}^n \cD_{A_k,t})\to \cD_{A_k,t},\qquad 
 P_{i,k}[d\oplus (\oplus_{j=1}^n d_j)]:= d_i
 \end{split}
\end{equation*}
 for any  $i=1,\ldots, n$. Note that $P_{\cD_0'}= P_{\cD'}$  and
  $P_{i,0}=P_i$ for any
  $i=1,\ldots, n$.
 Define the unitary operator 
 $W_k:\cN_k\to  \cN'_k$    by
 \begin{equation}\label{wk}
 W_k\left(\sum_{i=1}^n D_{A_k,t} T_i h_i\right) :=D_k'
 (\oplus_{i=1}^n A_kh_i)\oplus (\oplus_{i=1}^n D_{A_k,t}h_i).
 \end{equation}
Note that if  $k=0$ and $t=1$,   then we have $W_0=W$ 
(see  relation \eqref{W} in the particular case  when $\cT$ is a row isometry).
 
 Now we can prove the following result.
  
 \begin{theorem}\label{main2} Let $\cT:=[T_1~\cdots ~T_n]$, $T_i\in B(\cH)$,  
 be a row  isometry and let  $\cT':=[T_1'~\cdots ~T_n']$, $T_i'\in B(\cH')$, be a row 
 contraction
 with minimal isometric dilation $\cV'$.
 Let $A\in B(\cH,\cH')$,  $t>0$, be such that $\|A\|\leq t$ and  $AT_i=T_i'A$ 
 for any $i=1,\ldots, n$. For each $k=0,1,\ldots$,  
   let $A_k$ be  a $k$-step intertwining lifting  of $A$ with $\|A_k\|\leq t$.
 Then any   one-step intertwining lifting $A_{k+1}$ of $A_k$  such that $\|A_{k+1}\|\leq t$
 is given by
 \begin{equation}\label{one4}
  A_{k+1}=
   \left[\begin{matrix}
 A_k\\
  P_{\cD_k'}(W_kP_{\cN_k}+C_k  P_{\cM_k})D_{A_k,t}
 \end{matrix}\right]:\cH\to \cH_k'\oplus \cD_k',
\end{equation}
 where $C_k:\cM_k\to \cM_k'$ is a   contraction.
 Moreover, there is a one-to-one correspondence between the set of 
 all   one-step intertwining liftings $A_{k+1}$ of $A_k$ with $\|A_{k+1}\|\leq t$ and 
 the set of all
 contractions  $C_k:\cM_k\to \cM_k'$, and  
 there exists a unitary operator 
 $\Omega_k:\cD_{A_{k+1}, t}\to \left(\operatornamewithlimits{\oplus}
 \limits_{j=1}^n \cD_{A_k,t}\right)\oplus \cD_{C_k}$
 satisfying
  \begin{equation}\label{omegak}
   \Omega_k D_{A_{k+1},t}=
   \left[\begin{matrix}
   \left[\begin{matrix} P_{1,k}\\\vdots\\ P_{n,k}\end{matrix} \right]
   (W_kP_{\cN_k}+C_kP_{\cM_k})D_{A_k,t}\\
   D_{C_k} P_{\cM_k} D_{A_k,t}
   \end{matrix}\right].
 \end{equation}
 \end{theorem}
 
 \begin{proof}
 Let $A\in B(\cH,\cH')$ be  such that  $\|A\|\leq t$ and  $AT_i=T_i'A$ for any $i=1,\ldots, n$.
  Note that since $\cT$ is a row isometry and $\|A\|\leq t$,
  Theorem \ref{onestep2} remains true in a slightly adapted version. More precisely, taking into account
  that $\cD_\cT=\{0\}$ and replacing $D_A$ with $D_{A,t}$,
   we
  deduce that
    any  one-step intertwining lifting $A_1$ of $A$ such that $\|A_1\|\leq t$  is given by
 \begin{equation}\label{one3}
  A_1=
   \left[\begin{matrix}
 A\\
  P_{\cD'}(WP_\cN+C P_\cM)D_{A, t}
 \end{matrix}\right]:\cH\to \cH'\oplus \cD_{\cT'}.
\end{equation}
 where $C:\cM\to \cM'$ is a  contraction.
 Moreover, there is a one-to-one correspondence between the set of 
 all  one-step intertwining liftings $A_1$ of $A$ such that $A_1\|\leq t$ and the set of all
 contractions  $C:\cM\to \cM'$, and
 there exists a unitary operator 
 $\Omega:\cD_{A_1,t}\to (\oplus_{j=1}^n \cD_{A,t})\oplus \cD_C$
 satisfying
  \begin{equation*} 
   \Omega D_{A_1, t}=
   \left[\begin{matrix}
   \left[\begin{matrix} P_1\\\vdots\\ P_n\end{matrix} \right](WP_\cN+CP_\cM)D_{A,t}\\
   D_C P_\cM D_{A,t}
   \end{matrix}\right].
 \end{equation*}
   Since $A_{k+1}$ is a one-step intertwining lifting of $A_k$, we can 
    apply  the first part of this proof  to the $k$-step intertwining 
    lifting $A_k$ of $A$,  and complete the proof.
 \end{proof}
 
 According to Theorem \ref{main2}, there is a one-to-one correspondence between 
 the set of all  
 intertwining liftings of $A$ with tolerance $t$, and the set of all contractions 
 $C_k:\cM_k\to \cM_k'$, \ $k=1,2,\ldots$, given by
 $$
 B=\text{\rm SOT}-\lim_{k\to \infty} A_k P_{\cH_k},
 $$
  where  the sequence $\{A_k\}$ is  defined by  
 \eqref{one4}.
 We remark that if $t=1$,  then   this
 correspondence is, up to  unitary operators, exactly the one obtained in \cite{Po-intert}
  (see \cite{FF-book} for the case $n=1$) between 
 the set of all contractive
 intertwining liftings of $A$ and the set of all generalized choice sequences.

 Now, we show  that there is a one-to-one correspondence between the set of all 
intertwining liftings with tolerance $t>0$ and certain families of  contractions
  $\{\Lambda_\alpha\}_{\alpha\in \FF_n^+}$.  
 
  \begin{theorem}\label{geo}
 Let  $\cT:=[T_1~\cdots~ T_n]$, \ $T_i\in B(\cH)$, be a row isometry and 
 let  $\cT':=[T_1'~\cdots~ T_n']$, \ $T_i'\in B(\cH')$,  be a row contraction
 with its minimal isometric dilation  $\cV':=[V_1'~\cdots~ V_n']$, \ $V_i'\in B(\cK')$,
 on the Hilbert space
  $\cK':= \cH'\oplus [F^2(H_n)\otimes \cD']$.
 Let $A\in B(\cH, \cH')$ be such that $\|A\|\leq t$  and 
 $$
 AT_i=T_i'A, \quad i=1,\ldots, n.
 $$
  Then $B$ is an intertwining lifting of $A$ with respect to $\cV'$ and tolerance $t>0$
  if and only if  $B$ has a matrix decomposition 
  \begin{equation}\label{decomp}
  B=\left[\begin{matrix}A\\
  \Lambda D_{A, t}
  \end{matrix}
  \right]: \cH\to \cH'\oplus [F^2(H_n)\otimes \cD'],
  \end{equation}
 where $\Lambda:\cD_{A, t}\to F^2(H_n)\otimes \cD'$ is a contraction with
  $\Lambda h=\sum_{\alpha\in \FF_n^+} e_\alpha \otimes \Lambda_\alpha h$,
   $h\in \cD_{A, t}$, and  such that the operators $\Lambda_\alpha\in B(\cD_{A,t}, \cD')$ are given by
   \begin{equation}
   \begin{split}\label{one'}
   \Lambda_{g_0}&= P_{\cD'} W P_\cN+ Y_{g_0} P_\cM\\
   \Lambda_{g_j\alpha}&=\Lambda_\alpha P_jWP_\cN+ Y_{g_j\alpha} P_\cM,
   \end{split}
   \end{equation}
$j=1,\ldots, n$,  for some contractions  $Y_\alpha\in B(\cM, \cD')$, $\alpha\in \FF_n^+$, where 
 the subspaces $\cN$ and $\cM$  are defined by relation \eqref{nmk} (when $k=0$).
 Moreover, the families of contractions $\{\Lambda_\alpha\}_{\alpha\in \FF_n^+}$ and 
 $\{Y_\alpha\}_{\alpha\in \FF_n^+}$ uniquely determine each other.
 \end{theorem}
 
 \begin{proof}
 Since $B$ is an intertwining lifting of $A$ with respect lo $\cV'$ we must have
 $$
  B=\left[\begin{matrix}A\\
   X
  \end{matrix}
  \right]: \cH\to \cH'\oplus [F^2(H_n)\otimes \cD'].
  $$
  It is easy to see that $\|B\|\leq t$  if and only if $\|Xh\| \leq \|D_{A,t}h\|$, 
  and therefore, 
  if and only if there is a contraction 
  $\Lambda:\cD_{A, t}\to F^2(H_n)\otimes \cD'$   such that $X=\Lambda D_{A,t}$. Note that 
  $\Lambda$ is uniquely determined by $B$.
  Since $B T_i= V_i' B$ and  $A T_i= T_i' A$ for any $i=1,\ldots, n$,   we
  can use relation \eqref{dilmatr}  and get
  \begin{equation}
  \label{LDA}
  \Lambda D_{A, t} T_i= D_i'A+ (S_i\otimes I_{\cD'})\Lambda D_{A,t}
  \end{equation}
 for any $i=1,\ldots, n$.
 Setting 
 $$
 \Lambda h=\sum_{\alpha\in \FF_n^+} e_\alpha \otimes \Lambda_\alpha h, \quad 
   h\in \cD_{A, t},
   $$
      we infer that relation \eqref{LDA} is equivalent to
   $$
   \sum_{\beta\in \FF_n^+} e_\beta \otimes \Lambda_\beta D_{A,t} T_i h
   = 1\otimes D_i'Ah+ 
   \sum_{\alpha\in \FF_n^+} e_{g_i\alpha} \otimes \Lambda_\alpha D_{A,t} h
   $$
   for any  $h\in \cD_{A, t}$ and $i=1,\ldots, n$.
 It is clear that the latter equality is equivalent to the  following equations: 
 \begin{equation}
   \begin{split}\label{one''}
   \Lambda_{g_0}D_{A,t} T_i&=  D_i' A\\
   \Lambda_{g_j\alpha}D_{A,t} T_i&= \delta_{ij} \Lambda_\alpha D_{A,t}
   \end{split}
   \end{equation}
  for any   $i, j\in \{1,\ldots, n\}$ and $\alpha\in \FF_n^+$.
 We recall that
 \begin{equation}
 \label{PiPm}
 P_\cM D_{A, t} T_i h=0\quad \text{ and } \quad
   P_j W P_\cN D_{A,t} T_ih= \delta_{i j} D_{A,t} h
 \end{equation}
 for any $h\in \cH$ and $i, j\in \{1,\ldots, n\}$.
 On the other hand, relation \eqref{W} implies 
 \begin{equation}
 \label{alt}
 P_{\cD'} W\left(\sum_{i=1}^n D_{A, t} T_i h_i\right)= D'(\oplus_{i=1}^n A h_i).
 \end{equation}
  Taking into account  relations \eqref{PiPm},  \eqref{alt},  and   \eqref{one''}, we deduce 
  \begin{equation*}
   \begin{split}
   P_{\cD'} W P_\cN D_{A,t} T_i&= P_{\cD'} W  D_{A,t} T_i= D_i'A\\
   &=\Lambda_{g_0} D_{A, t} T_i =\Lambda_{g_0} P_\cN D_{A, t} T_i
   \end{split}
  \end{equation*}
  for any $i=1,\ldots, n$, 
  which implies 
  $$
  \Lambda_{g_0}= P_{\cD'} W P_\cN+ Y_{g_0} P_\cM,
  $$
  where $Y_{g_0}\in B(\cM, \cD')$.
  Using again relations \eqref{PiPm}, \eqref{alt},  and   \eqref{one''},
  we obtain 
  \begin{equation*}
   \begin{split}
   \Lambda_\alpha P_j W P_\cN D_{A,t} T_i&=
   \delta_{ij} \Lambda_\alpha D_{A,t}= \Lambda_{g_i\alpha} D_{A, t} T_i\\
   &= \Lambda_{g_i\alpha}P_\cN D_{A, t} T_i
   \end{split}
  \end{equation*}
  for any $i,j\in\{1,\ldots, n\}$, 
   which implies
   $$
    \Lambda_{g_j\alpha}=\Lambda_\alpha P_jWP_\cN+ Y_{g_j\alpha} P_\cM,\quad j=1,\ldots, n,
    $$
    where $Y_{g_j \alpha}\in B(\cM, \cD')$.
 Now, it is easy to see that  
 the families of contractions $\{\Lambda_\alpha\}_{\alpha\in \FF_n^+}$ and 
 $\{Y_\alpha\}_{\alpha\in \FF_n^+}$ uniquely determine each other.
 The proof is complete.
 \end{proof}

 \bigskip
 
 \subsection{Central lifting in several variables and geometric characterizations}
 \label{central}

The main results of this section (see Theorem \ref{main3} 
 and  Theorem \ref{restr=0})
 show that the intertwining lifting corresponding to
   the parameters $C_k=0$, $k=1,2,\ldots$, 
   (resp.~$\Lambda_\alpha=0$, $\alpha\in \FF_n^+$)  coincides with the central 
   intertwining lifting  with tolerance $t$,
   for which we have an explicit form.
 At the end of this section, we show that, under certain conditions,
 there is only one intertwining lifting $B$ of $A$ such that $\|B\|=\|A\|$, namely the 
  central 
   intertwining lifting.
   The geometric structure of the  central 
   intertwining lifting will play a very important role in our investigation.

 Our first  result shows that the  intertwining
  lifting of $A$ corresponding to the parameters $C_k=0$ for any  
  $k=1,2,\ldots$, coincides with the  central intertwining lifting $B_c$ of 
  $A$ (as defined  in \cite{Po-central} when $t=1$).

 \begin{theorem}\label{main3}
 Let $A\in B(\cH,\cH')$, $t>0$, be  such that $\|A\|\leq t$ 
  and  $AT_i=T_i'A$ for any $i=1,\ldots, n$, 
 where $\cT:=[T_1~\cdots ~T_n]$, $T_i\in B(\cH)$,  is an isometry and 
 $\cT':=[T_1'~\cdots ~T_n']$, $T_i'\in B(\cH')$,  is a row contraction. Let
 $\cK':=\cH'\oplus [F^2(H_n)\otimes \cD']$ and 
  $\cV':=[V_1'~\cdots ~V_n']$, $V_i'\in B(\cK')$, 
   be the minimal isometric dilation of $\cT'$.  
   
  Then the  intertwining
  lifting of $A$ with tolerance $t$ corresponding to the parameters $C_k=0$,    
  $k=1,2,\ldots$,  is given by
  $B_c:\cH\to \cH'\oplus [F^2(H_n)\otimes \cD']$ and 
  \begin{equation}\label{Bc}
  B_ch=Ah\oplus \sum_{\sigma\in \FF_n^+} e_\sigma 
  \otimes (P_{\cD'}W P_\cN)E_{\tilde\sigma} D_{A,t} h, \quad  h\in \cH,
  \end{equation}
  where $E_{g_0}:=I_{\cD_{A,t}}$ and $E_i:=  P_{i} W P_\cN$ for any  $i=1,\ldots, n$.
  In particular, we have $B_c^*|\cH'= A^*$ and $\|B_c\|\leq t$.
 \end{theorem}

 \begin{proof}
 Let $\{A_k\}_{k=0}^\infty$ be the set of
  the $k$-step   intertwining liftings of $A$ with $\|A_k\|\leq t$, obtained by 
  setting $C_k=0$ for
   any $k=0,1,\ldots$. Then, for each $k=0,1,\ldots$,  the operator 
    $A_{k+1}$ is the one-step   
   intertwining lifting of $A_k$ when $C_k=0$, and $A_0=A$.
   According to Theorem \ref{main2}, we have
   $$
   A_1h=Ah\oplus 1\otimes (P_{\cD'} WP_\cN)D_{A,t}h, \quad h\in \cH.
   $$ 
  By induction, we assume that
  \begin{equation}\label{Bk}
  A_kh=Ah\oplus \sum_{|\sigma|\leq k-1} e_\sigma 
  \otimes (P_{\cD'}W P_\cN)E_{\tilde\sigma} D_{A,t} h, \quad  h\in \cH.
  \end{equation}
 Note that, for each $h\in \cH$, we have
 \begin{equation*}
 \begin{split}
 t^2&\|h\|^2-\|A_kh\|^2=t^2\|h\|^2 -\|Ah\|^2-
 \sum_{|\sigma|\leq k-1}  
   \|(P_{\cD'}W P_\cN)E_{\tilde\sigma} D_{A,t} h\|^2\\
   &=
   \|D_{A,t}h\|^2- \sum_{|\sigma|\leq k-1}
   \left(\| W P_\cN E_{\tilde\sigma} D_{A,t} h\|^2-\sum_{i=1}^n\|(P_iWP_\cN)
   E_{\tilde\sigma} D_{A,t}h\|^2\right)\\
   &=
   \|D_{A,t}h\|^2-\sum_{|\sigma|\leq k-1} 
   \|  P_\cN E_{\tilde\sigma} D_{A,t} h\|^2+\sum_{i=1}^n \sum_{|\sigma|\leq k-1} 
   \|E_{\tilde{\sigma g_i}} D_{A,t}h\|^2\\
   &=
   \|D_{A,t}h\|^2-\|P_\cN D_{A,t}h\|^2-
   \sum_{1\leq |\sigma|\leq k-1} 
   \|  P_\cN E_{\tilde\sigma} D_{A,t} h\|^2 +
   \sum_{1\leq|\omega|\leq k} 
   \| E_{\tilde\omega} D_{A,t} h\|^2\\
   &=
   \|P_\cM D_{A,t}h\|^2+\sum_{1\leq |\sigma|\leq k-1} 
   \| P_\cM E_{\tilde\sigma} D_{A,t} h\|^2+\sum_{|\omega|=k} \|E_{\tilde\omega} D_{A,t} h\|^2.
 \end{split}
 \end{equation*}
 Hence, we infer that $\|A_k\|\leq t$ and 
 $$
 \|D_{A_k,t}h\|^2=
 \sum_{|\sigma|\leq k-1} 
   \| P_\cM E_{\tilde\sigma} D_{A,t} h\|^2+\sum_{|\omega|=k}
    \|E_{\tilde\omega} D_{A,t} h\|^2.
    $$
 Now, it is clear that there is an isometry
 $$
 Z_k: \cD_{A_k,t}\to  \left[ \operatornamewithlimits{\oplus}
 \limits_{|\sigma|\leq k-1} e_\sigma\otimes 
 \cM\right]\oplus
 \left[ \operatornamewithlimits{\oplus}\limits_{|\omega|=k} e_\omega\otimes 
 \cD_{A,t}\right]
 $$
 satisfying the equation
 \begin{equation}
 \label{Zk}
 Z_kD_{A_k,t}h= 
 \left(\sum_{ |\sigma|\leq k-1} e_\sigma\otimes P_\cM E_{\tilde\sigma} D_{A,t}h\right)\oplus
 \left(\sum_{|\omega|=k} e_\omega\otimes E_{\tilde\omega} D_{A,t}h\right) 
 \end{equation}
 for any $h\in \cH$.
 
 We show now that $Z_k$ is a unitary operator.
 First, note that
 $$
 P_jWP_\cN D_{A,t}T_ih=\delta_{i,j} D_{A,t}h \quad \text{ and } \quad
  P_\cM D_{A,t} T_ih=0
  $$
   for any $h\in \cH$
  and 
  $i,j\in \{1,2,\ldots, n\}$. 
  Hence, we infer that if $|\sigma|=|\omega|$, then
  \begin{equation}\label{EDT}
  E_{\tilde\sigma} D_{A,t}T_\omega h=\delta_{\omega, \sigma} D_{A,t}h, \quad h\in \cH.
 \end{equation}
 Hence, and using \eqref{Zk}, we infer that if $|\omega|=k$, then
 $$
 Z_kD_{A_k,t}\left(
 \sum_{|\omega|=k} T_\omega h_\omega\right)=\sum_{|\omega|=k} 
 e_\omega\otimes D_{A,t}h_\omega.
 $$
 On the other hand, 
 $$
 Z_kD_{A_k,t}\left(
 \sum_{|\omega|=k-1} T_\omega h_\omega\right)=\left(\sum_{|\omega|=k} 
 e_\omega\otimes d_\omega\right) \oplus \left(\sum_{|\omega|=k-1} 
 e_\omega\otimes P_\cM D_{A,t} h_\omega\right)
 $$
 for some $d_\omega\in \cD_{A,t}$.
 The last two equations show that the closed  linear span of the subspaces 
 $Z_kD_{A_k,t}\left(
 \bigvee_{|\omega|=k-1} T_\omega \cH\right)$ and 
 $Z_kD_{A_k,t}\left(
 \bigvee_{|\omega|=k} T_\omega \cH\right)$ 
 coincides with
 $$
 \left[ \operatornamewithlimits{\oplus}\limits_{|\omega|=k} e_\omega\otimes 
 \cD_{A,t}\right]\oplus \left[ \operatornamewithlimits{\oplus}\limits_{|\sigma|=k-1} e_\omega\otimes 
 \cM\right].
 $$
 Continuing this process, one can show that $Z_k$ is surjective,
  and consequently a unitary operator.
 Since $\cN_k= \bigvee_{i=1}^n D_{A_k,t} T_i \cH$, it is easy to see that 
 \begin{equation}
 \label{ZNM} \begin{split}
 Z_k\cN_k
 &=
 \left[ \operatornamewithlimits{\oplus}\limits_{|\omega|=k} e_\omega\otimes 
 \cD_{A,t}\right]\oplus \left[ \operatornamewithlimits{\oplus}\limits_{1\leq 
 |\sigma|\leq k-1} e_\omega\otimes 
 \cM\right]\quad \text{ and }\\
 Z_k\cM_k&=1\otimes \cM.
 \end{split}
 \end{equation}
 According to Theorem \ref{main2} (when $C_k=0$),
 to complete our proof by induction, we need to show that
 \begin{equation}\label{formula}
 P_{\cD_k'} W_k P_{\cN_k} D_{A_k,t}h = \sum_{|\sigma|=k} e_\sigma\otimes 
 (P_{\cD'} W P_{\cN})E_{\tilde\sigma} D_{A,t}h, \quad h\in \cH,
 \end{equation}
 where $\cD_k'$ is identified with 
 $\operatornamewithlimits{\oplus}\limits_{|\omega|=k} e_\omega\otimes 
 \cD'$ and $\cD':= \cD_{\cT'}$.
 Using relations \eqref{Zk} and \eqref{ZNM}, we obtain
 
 \begin{equation*}\begin{split}
 &P_{\cD_k'} W_k P_{\cN_k} D_{A_k,t}h =P_{\cD_k'} W_k P_{\cN_k}Z_k^* Z_k D_{A_k,t}h\\
 &\overset {\eqref {Zk}}=
 P_{\cD_k'} W_k P_{\cN_k}Z_k^*
 \left(\sum_{|\omega|=k} e_\omega \otimes E_{\tilde\omega} D_{A,t}h\oplus 
 \sum_{|\sigma|\leq k-1} e_\sigma \otimes P_\cM E_{\tilde\sigma} D_{A,t}h\right)\\
 &\overset {\eqref {ZNM}} =
 P_{\cD_k'} W_k Z_k^*
 \left(\sum_{|\omega|=k} e_\omega \otimes E_{\tilde\omega} D_{A,t}h\oplus 
 \sum_{1\leq|\sigma|\leq k-1} e_\sigma \otimes P_\cM E_{\tilde\sigma} D_{A,t}h\right)\\
 &=
 \sum_{j=1}^n P_{\cD_k'} W_k Z_k^*
 \left(\sum_{|\alpha|=k-1} e_{g_j\alpha}\otimes E_{\tilde\alpha} E_{g_j} D_{A,t}h\oplus 
 \sum_{0\leq|\beta|\leq k-2} e_{g_j\beta} \otimes P_\cM E_{\tilde\beta} E_{g_j}
  D_{A,t}h\right).
 \end{split}
 \end{equation*}
 
 Now, using relations \eqref{EDT}, \eqref{Zk}, we deduce that
 \begin{equation*}\begin{split}
 &P_{\cD_k'} W_k Z_k^*
 \left(\sum_{|\alpha|=k-1} e_{g_j\alpha}\otimes E_{\tilde\alpha } D_{A,t}y\oplus 
 \sum_{0\leq|\beta|\leq k-2} e_{g_j\beta} \otimes P_\cM E_{\tilde\beta} 
  D_{A,t}y\right)\\
  &\overset{\eqref {EDT}} =P_{\cD_k'} W_k Z_k^*\sum_{i=1}^n 
  \left(\sum_{|\alpha|=k-1} e_{g_i\alpha}\otimes E_{\tilde{g_i\alpha}}  D_{A,t}T_jy\oplus 
 \sum_{0\leq|\beta|\leq k-2} e_{g_i\beta} \otimes P_\cM E_{\tilde{g_i\beta}} 
  D_{A,t} T_jy\right)\\
  &\overset {\eqref {Zk}} =
  P_{\cD_k'} W_k P_{\cN_k} Z_k^* Z_k D_{A_k,t} T_jy=
  P_{\cD_k'} W_k P_{\cN_k}  D_{A_k,t} T_jy\\
  &\overset{\eqref {wk}}=
  D_k'(\underbrace{0,\ldots,0}_{{j-1}\mbox{ \scriptsize
times}}, A_k y,0,\ldots)= 
S_j P_{\operatornamewithlimits{\oplus}\limits_{|\omega|=k} (e_\omega\otimes 
 \cD')} A_ky\\
 &=
\sum_{|\alpha|=k-1} e_{g_j}
 e_\alpha \otimes (P_{\cD'} W P_\cN) E_{\tilde\alpha} D_{A,t}y.
 \end{split}
 \end{equation*}
 Here, we  used the fact that $D_k'=
 \operatornamewithlimits{\oplus}\limits_{j=1}^n
 P_{\operatornamewithlimits{\oplus}\limits_{|\omega|=k} (e_\omega\otimes 
 \cD')}$
 and $\cD_k'$ is identified with
 $$
 \operatornamewithlimits{\oplus}\limits_{j=1}^n 
 \operatornamewithlimits{\oplus}\limits_{|\omega|=k}
  (e_{g_j\omega}\otimes 
 \cD').
 $$
 Now, since $D_{A,t}$ has dense range in $\cD_{A,t}$, we infer that 
 \begin{equation*} 
 \begin{split}
 P_{\cD_k'} W_k Z_k^*
 \left(\sum_{|\alpha|=k-1} e_{g_j\alpha}\otimes E_{\tilde \alpha } x\right. &\oplus 
 \left.\sum_{0\leq|\beta|\leq k-2} e_{g_j\beta} \otimes P_\cM E_{\tilde\beta} x \right)\\
  &= 
 \sum_{|\alpha|=k-1} e_{g_j}
 e_\alpha \otimes (P_{\cD'} W P_\cN) E_{\tilde \alpha} x
 \end{split}
 \end{equation*}
 for any $x\in \cD_{A,t}$ and $j=1,2,\ldots, n$. 
 Summing up the results obtained so far, we  deduce
 \begin{equation*}
 \begin{split}
 P_{\cD_k'} W_k P_{\cN_k} D_{A_k,t}h 
 &=
 \sum_{j=1}^n P_{\cD_k'} W_k Z_k^*
 \left(\sum_{|\alpha|=k-1} e_{g_j\alpha} 
 \otimes E_{\tilde\alpha} (E_{g_j} D_{A,t}h)\right.\\ 
 &\qquad\left.\oplus 
 \sum_{0\leq|\beta|\leq k-2} e_{g_j\beta} \otimes P_\cM E_{\tilde\beta} (E_{g_j} D_{A,t}h )
 \right)\\
 &=
 \sum_{j=1}^n \sum_{|\alpha|=k-1} e_{g_j\alpha}
  \otimes (P_{\cD'} W P_\cN) E_{\tilde\alpha} (E_{g_j} D_{A,t}h )\\
  &= \sum_{|\beta|=k} e_{\beta}
  \otimes (P_{\cD'} W P_\cN) E_{\tilde\beta} D_{A,t}h ).
 \end{split}
 \end{equation*}
 Therefore, relation \eqref{formula} is proved. Since $\|A_k\|\leq k$ for any $k=1,2,\ldots$, we get
 $\|B_c\|\leq t$,  and the   proof is complete.
 \end{proof}

 Note that  if $t=\|A\|$, then Theorem \ref{main3} implies 
 that the central intertwining lifting $B_c$ of $A$
 satisfies $\|B_c\|=\|A\|$.
 
 The following result  provides another characterization for the 
 central intertwining lifting in the noncommutative commutant lifting theorem.
 We employ the notation from Theorem \ref{geo}.

 \begin{theorem}\label{restr=0}
 Let  $\cT:=[T_1~\cdots~ T_n]$, \ $T_i\in B(\cH)$, be a row isometry and 
 let  $\cT':=[T_1'~\cdots~ T_n']$, \ $T_i'\in B(\cH')$,  be a row contraction
 with its minimal isometric dilation  $\cV':=[V_1'~\cdots~ V_n']$, \ $V_i'\in B(\cK')$,
 on the Hilbert space
  $\cK':= \cH'\oplus [F^2(H_n)\otimes \cD']$.
 Let $A\in B(\cH, \cH')$ be such that $\|A\|\leq t$  and 
 $$
 AT_i=T_i'A, \quad i=1,\ldots, n.
 $$
 Let $B$ be an intertwining liftig of $A$ with tolerance $t>0$, and let 
 $B=\left[\begin{matrix}A\\
  \Lambda D_{A,t}
  \end{matrix}
  \right]$
  be the corresponding decomposition.
  Then the following statements are equivalent:
  \begin{enumerate}
  \item[(i)] $B$ is the central intertwining lifting of $A$ with tolerance $t>0$;
  \item[(ii)] $\Lambda|\cM=0$;
  \item[(iii)] $\Lambda_\alpha|\cM=0$ for any $\alpha\in \FF_n^+$;
  \item[(iv)] $Y_\alpha|\cM=0$ for any $\alpha\in \FF_n^+$.
  \end{enumerate}
 In particular, if $\cM=\{0\}$, then 
  the central intertwining lifting of $A$ with tolerance $t>0$ is the unique 
  intertwining liftig $B$  of $A$ such that $\|B\|\leq t$.
 \end{theorem}
 
 \begin{proof}
 Note that, using Theorem \ref{geo}, Theorem \ref{main3}, and  relations \eqref{decomp}, \eqref{one'}, and \eqref{Bc},
   the  result follows.
 \end{proof}

 Let  $\cT':=[T_1'~\cdots~ T_n']$, \ $T_i'\in B(\cH')$,  be a row contraction and let
   $\cV':=[V_1'~\cdots~ V_n']$, \ $V_i'\in B(\cK')$,
 be its minimal isometric dilation on the Hilbert space
  $\cK':= \cH'\oplus [F^2(H_n)\otimes \cD']$.
 If  $\cW':=[W_1'~\cdots~ W_n']$, \ $W_i'\in B(\cG')$, is an arbitrary
  isometric lifting of $\cT'$ on a Hilbert space 
 $\cG'\supseteq \cH'$, then there exists a unique isometry
 $\Phi : \cK'\to \cG'$ such that 
 $\Phi V_i'= W_i'\Phi$,  for any  $i=1,\ldots, n$, and $\Phi|\cH'= I_{\cH'}$.
 Moreover, using the geometric structure of the minimal isometric dilation 
  (see \cite{Po-isometric}), one can deduce  that $\Phi$ is given by
  $$
  \Phi\left(h\oplus \sum_{\alpha\in \FF_n^+} e_\alpha \otimes D_{\cT'} h_\alpha\right)
  = h+\sum_{\alpha\in \FF_n^+}W_\alpha'[W_1'-T_1'~\cdots ~ W_n'-T_n']h_\alpha
  $$
  for any $h\oplus \sum\limits_{\alpha\in \FF_n^+} e_\alpha \otimes 
   h_\alpha$ in ~$\cH'\oplus [F^2(H_n)\otimes \cD']$.
 Indeed, we have
 \begin{equation*}
 \begin{split}
 \Phi\left(0\oplus \sum_{\alpha\in \FF_n^+} e_\alpha \otimes D_{\cT'} h_\alpha\right)
 &=\Phi\left(0\oplus \sum_{\alpha\in \FF_n^+} V_\alpha'(1 \otimes D_{\cT'} h_\alpha)\right)\\
 &=
 \sum_{\alpha\in \FF_n^+}\Phi V_\alpha'[V_1'-T_1'~\cdots ~ V_n'-T_n']h_\alpha\\
 &=
 \sum_{\alpha\in \FF_n^+}W_\alpha' [W_1'-T_1'~\cdots ~ W_n'-T_n']h_\alpha.
 \end{split}
 \end{equation*}
Let $\cT:= [T_1~\cdots ~T_n]$, $ T_i\in B(\cH)$, be an isometry and let
 $A\in B(\cH, \cH')$ be an 
 operator such that $\|A\|\leq t$ and  $AT_i= T_i' A$  for any $i=1,\ldots, n$.
  If $B$ is an intertwining 
lifting of $A$  
with respect to $\cV'$ such that $\|B\|\leq t$, then $\Phi B$ is an intertwining lifting of $A$ with 
respect to $\cW'$. Since $\Phi$ is  unique, we call the operator 
$\tilde{B_c}:= \Phi B_c$ the central intertwining 
lifting of $A$ with respect to $\cW'$ and tolerance $t$.
Using these remarks and Theorem \ref{main3}, we can easily obtain the following form 
for $\tilde{ B_c}$.

\begin{proposition}
 Let $\cT:= [T_1~\cdots ~T_n]$, \ $T_i\in B(\cH)$, be  an isometry and 
  $\cT':= [T_1'~\cdots ~T_n']$, \ $T_i'\in B(\cH')$, be a row contraction. 
 Let  $\cW':= [W_1'~\cdots ~W_n']$, \ $W_i'\in B(\cG')$, be  an isometric lifting of
 of $\cT'$, and $A\in B(\cH, \cH')$ be such that $\|A\|\leq t$ and 
 $$
 AT_i =T_i'A, \quad i=1,\ldots, n.
 $$
Then the central intertwining lifting $\tilde{B_c}$ of $A$ with 
respect to $\cW'$  and tolerance $t>0$ is given by
$$
\tilde{B_c}= A+\sum_{\alpha\in \FF_n^+} W_\alpha'(P'UP_\cN) E_{\tilde \sigma}' D_{A,t},
$$
where $E_{g_0}:= I_{\cD_{A,t}}$, \  $E_i':= P_i'U P_\cN$, \ $i=1,\ldots, n$, 
$P' ($resp.  $P_i')$ 
 is the orthogonal projection of 
~$\cG'\oplus (\oplus_{j=1}^n \cD_{A,t})$~ onto ~$\cG'  ($resp. the i-th component 
of ~$\oplus_{j=1}^n \cD_{A,t}),$ and the operator ~$U:\cN\to \cG' \oplus (\oplus_{j=1}^n \cD_{A,t})$ is the isometry
defined by
$$
U\left( \sum_{i=1}^n D_{A,t} T_i h_i\right):=\left(
\sum_{i=1}^n(W_i'-T_i') Ah_i\right)\oplus \left(
\oplus_{i=1}^n D_{A,t} h_i \right).
$$
\end{proposition}

At the end of this section, we show that, under certain conditions,
 there is only one intertwining 
 lifting $B$ of $A$ such that $\|B\|=\|A\|$. 
 We say that an operator $A\in B(\cH, \cH')$ attains its norm if there is a 
 vector $h\in \cH$ of norm one such that $\|Ah\|=\|A\|$.

\begin{proposition}\label{M0}
Let $\cT':= [T_1'~\cdots ~T_n']$, \ $T_i'\in B(\cH')$, be a row contraction,
 and let $A: F^2(H_n)\to \cH'$ be a contraction which attains 
 its norm  $\|A\|=1$.
 If 
 $$
 AS_i= T_i'A, \quad i=1,\ldots, n,
 $$
 then 
 \begin{equation}
 \label{MZERO}
 \cM:= \cD_A \ominus \bigvee_{i=1}^n D_A S_i F^2(H_n)=\{0\}. 
 \end{equation}
 In particular, if  $\cV':= [V_1'~\cdots ~V_n']$ is the minimal isometric dilation of $\cT'$,
  then there is only one intertwining lifting $B$ of $A$
   with respect to $\cV'$ such that $\|B\|=\|A\|$.
\end{proposition}

\begin{proof}
First, note that $A$ attains its norm at a vector $f\in F^2(H_n)$ if and only if $f\in \ker D_A$.
We need to prove that $A$ attains its norm at    $f\in F^2(H_n)$ such that
$f(0)\neq 0$.
To this end, 
assume that $A$ attains its norm at a vector  $g\in F^2(H_n)$ such that
$g(0)=0$. 
 Let $g:= \sum\limits_{\alpha\in \FF_n^+} a_\alpha e_\alpha$  and
let $k$ be the least positive integer such that there is  $\alpha_0\in \FF_n^+$
with $|\alpha_0|=k$ and  $a_{\alpha_0}\neq 0$.
Then $g$ can be written as $g= \sum\limits_{|\alpha|=k} S_\alpha \varphi_\alpha$, 
where $\varphi_\alpha\in F^2(H_n)$. Since $\|g\|=1$ and $ \| Ag\|=\|A\|$, we have
\begin{equation*}
\begin{split}
\|A\|^2 \sum_{|\alpha|=k} \|\varphi_\alpha\|^2&=
\|A\|^2 \|g\|^2= \|Ag\|^2\\
&= \left\|A\left(\sum_{|\alpha|=k} S_\alpha \varphi_\alpha\right)\right\|^2=
\left\| \sum_{|\alpha|=k} T_\alpha' A \varphi_\alpha\right\|^2 \\
&\leq   \sum_{|\alpha|=k}\|A\varphi_\alpha\|^2\leq \|A\|^2 
\sum_{|\alpha|=k}\|\varphi_\alpha\|^2.
\end{split}
\end{equation*}
Therefore, we must have 
$$
\|A \varphi_\alpha\|= \|A\| \|\varphi_\alpha\|,\quad \text{ if } \  |\alpha|=k.
$$
Setting $f:= \frac{\varphi_{\alpha_0}} {\|\varphi_{\alpha_0}\|} $, we have 
$\|Af\|=\|A\|$, $\|f\|=1$, and $f(0)\neq 0$.
If we assume  that there is $\psi\in \cM$, $\psi\neq 0$,
then $\left< \psi D_AS_i h\right>=0$ for any $h\in F^2(H_n)$ and $i=1,\ldots, n$.
This implies
$D_A \psi=a_0\in \CC$, $ a_0\neq 0$. Since $f\in \ker D_A$, we have
$$
0= \left< f, D_A \psi \right> =\left< f, a_0\right>.
$$
On the other hand, $f(0)\neq 0$ and $a_0\neq 0$ imply $\left< f, a_0\right>\neq 0$, 
which is a contradiction. Therefore, $\cM=\{0\}$. 
The uniqueness of the intertwining lifting follows from Theorem \ref{onestep2} and 
Theorem \ref{main2}, noticing that $C_k=0$ for any $k=1,2,\ldots$.
The proof is complete.
\end{proof}

It is easy to see that one can obtain  a version of Proposition \ref{M0} if $\|A\|=t>0$. 
In this case,  we have to replace $D_A$ by $D_{A,t}$.

 \bigskip

\subsection{A maximum principle  for
 the noncommutative commutant lifting theorem}\label{max} 

 In this section, we prove a maximum principle  for
 the noncommutative commutant lifting theorem, which also provides
  a new characterization for the central intertwining lifting 
  (see Theorem \ref{maxi} and Theorem \ref{maxprin2}).
  This is a key result used in Section \ref{max-entro-sol} to find the maximal 
    entropy solution  for the noncommutative commutant lifting theorem.

Let $T_1,\ldots, T_n\in B(\cH)$  be isometries with orthogonal ranges and let 
$\cL:=\bigcap_{i=1}^n \ker T_i^*$ be the wandering subspace in
 the Wold type decomposition \cite{Po-isometric}.   If  $X\in B(\cH, \cH')$  
satisfies $\|X\|\leq t$  ~
for some
$t>0$, 
 we define 
the positive operator $\Delta(X)\in B(\cL)$ by setting
\begin{equation}\label{DX}
\left< \Delta(X)\ell, \ell\right>:= 
\inf\{\|D_{X,t}(\ell-T_1h_1-\cdots-T_nh_n)\|^2:\ h_1,\ldots, h_n\in \cH\}
\end{equation}
for any $\ell\in\cL$. Following the classical case (see \cite{FFGK-book}),
 we call the operator  $\Delta(X)$
 the Schur complement of $D_{X,t}^2$ with respect to $T_1,\ldots, T_n$.
 Note that if $\|X\|<t$, then $\Delta(X)$ is indeed 
 the Schur complement of $D_{X,t}^2$ with respect to the  orthogonal decomposition
 $$
 \cH=\cL\oplus[T_1\cH\oplus\cdots \oplus T_n\cH].
 $$
  It is easy to see that 
 $$
 \left< \Delta(X)\ell, \ell\right>:= 
\inf\{\left< D_{X,t}^2h,h\right>:\ h\in\cH \text{ and }  P_\cL h=\ell\}.
$$
On the other hand,  one can prove that
\begin{equation}\label{Dform}
\Delta(X)= P_\cL D_{X,t} P_{\cM_X} D_{X,t}|\cL,
\end{equation}
where 
$$
\cM_X:= \cD_{X,t}\ominus \bigvee_{i=1}^n D_{X,t} T_i\cH
$$
 and $P_\cL$, $P_{\cM_X}$ are   the orthogonal
 projections onto $\cL$ and $\cM_X$, respectively.
 Indeed, we have 
 \begin{equation*} \begin{split}
 \left< \Delta(X)\ell, \ell\right>&= 
\inf\{\|D_{X,t}(\ell-T_1h_1-\cdots-T_nh_n)\|^2:\ h_1,\ldots, h_n\in \cH\}\\
&=\inf\{\|D_{X,t} \ell-k\|^2:\ k\in \bigvee_{i=1}^n D_{X,t} T_i\cH\}\\
&= \|P_{\cM_X}D_{X,t}\ell\|^2=\left< (P_\cL D_{X,t} P_{\cM_X} D_{X,t}|\cL)\ell, 
\ell\right>
 \end{split}
 \end{equation*}
 for any $\ell\in \cL$.
 Hence, we  obtain relation  \eqref{Dform}.
 
We can prove now  a maximum principle for  the noncommutative commutant lifting theorem.
This  result also provides  a new characterization for the central
  intertwining lifting.

\begin{theorem}\label{maxi}
 Let $A\in B(\cH,\cH')$ be  an operator   satisfying $\|A\|\leq t$ and  $AT_i=T_i'A$ for any
  $i=1,\ldots, n$, 
 where $\cT:=[T_1~\cdots ~T_n]$, $T_i\in B(\cH)$,  is an isometry and 
 $\cT':=[T_1'~\cdots ~T_n']$, $T_i'\in B(\cH')$,  is a row contraction.
 Let $\cV' $ be the minimal isometric dilation of $\cT'$ and let 
  $B_c$ be  the central intertwining lifting of $A$ with respect to $\cV'$. Then
 \begin{equation}\label{DDD}
 \Delta (A)=\Delta(B_c)\geq \Delta(B)
 \end{equation}
for any   intertwining lifting $B$ of $A$ with tolerance $t$.
Moreover, $\Delta(B_c)= \Delta(B)$ if and only if $B_c=B$.
\end{theorem}

\begin{proof}

According to Theorem \ref{main2}, for each $k=0,1,\ldots$, there exists a unitary operator 
 $\Omega_k:\cD_{A_{k+1},t}\to \left(\operatornamewithlimits{\oplus}
 \limits_{j=1}^n \cD_{A_k,t}\right)\oplus \cD_{C_k}$ such that
 $$
 P_{i,k} W_k P_{\cN_k} \left(\sum_{j=1}^n D_{A_k,t} T_j h_j\right)= D_{A_k,t} h_i,
  \quad i=1,\ldots, n.
 $$
 Hence, and using relation \eqref{omegak}, we have
 $$
 \Omega_kD_{A_{k+1},t} \left(\sum_{j=1}^n T_j h_j\right)= D_{A_k,t}h_1\oplus D_{A_k,t}h_2\oplus
 \cdots \oplus D_{A_k,t}h_n.
 $$
 Therefore, we have  $\Omega_k \cN_{k+1}= \oplus_{j=1}^n \cD_{A_k,t} \oplus \{0\}
 $
 and $\Omega_k \cM_{k+1}=\{0\}\oplus \cD_{C_k}$.
 Using again Theorem \ref{main2}, we infer that
 $$
 \|\Omega_k P_{\cM_{k+1}} D_{A_{k+1},t} x\|=\|D_{C_k} P_{\cM_{k}} D_{A_{k},t} x\|
$$
for any $x\in \cH$.
This implies
\begin{equation} \label{ine1}
\|P_{\cM_{k}} D_{A_{k},t} x\|\geq \|D_{C_k} P_{\cM_{k}} D_{A_{k},t} x\|
=\|P_{\cM_{k+1}} D_{A_{k+1},t} x\|,\quad x\in \cH.
\end{equation}
Moreover, the equality holds for any $x\in \cH$ if and only if $C_k=0$ for any 
$k=0,1,\ldots$.
Hence, we deduce that
\begin{equation}\label{PGD}
\|P_\cM D_{A,t}x\|\geq \|P_{\cM_1} D_{A_1,t}x\|\geq\cdots \geq \|P_{\cM_k} D_{A_k,t}x\|.
\end{equation}
Now, let $x\in \cH$ and note that
\begin{equation}\label{dax}
\|D_{A_k,t} x\|^2=t^2\| x\|^2-\|A_kx\|^2=\|D_{B,t}x\|^2+\|(I-P_k) Bx\|^2.
\end{equation}
Hence, $\|D_{A_k,t} x\|^2\geq \|D_{B,t}x\|^2$, which implies
$$
\|P_{\cM_B}D_{B,t}h\|\leq \left\|D_{B,t}h-D_{B,t}\left(\sum_{j=1}^n T_j h_j\right)
\right\|\leq 
\left\|D_{A_k,t}h-D_{A_k,t}\left(\sum_{j=1}^n T_j h_j\right)\right\|
$$
for any $h_1,\ldots, h_n\in \cH$.
Taking the infimum over $h_1,\ldots, h_n\in \cH$, we obtain
\begin{equation}\label{ineg}
\|P_{\cM_B}D_{B,t}h\|\leq \|P_{\cM_k} D_{A_k,t} h\|,\quad h\in \cH.
\end{equation}
Combining \eqref{PGD} with \eqref{ineg}, we infer that 
$$
\|P_{\cM_B}D_{B,t}h\|\leq \|P_\cM D_{A,t}h\|, \quad  h\in \cH,
$$ 
and 
$\|P_{\cM_B}D_{B,t}h\|= \|P_\cM D_{A,t}h\|$ for any $h\in \cH$ if and only if
 $C_k=0$ for any 
$k=0,1,\ldots$, which means that $B=B_c$, the central 
intertwining lifting of $A$ with tolerance $t$. 
Taking into account  relation \eqref{Dform} and the fact that 
$$
P_{\cM_k} D_{A_k,t} y
=P_{\cM_B} D_{B,t} y=0
$$
for any $y\in \bigvee_{j=1}^n T_j \cH$,
 we complete
 the proof of the theorem.
\end{proof}

Let us remark that if $\{A_k\}_{k=0}^\infty$ is the sequence 
of $k$-step intertwining liftings of $A$ 
corresponding to $B$, then
\begin{equation}\label{DB}
\Delta(B)=\text {\rm SOT}-\lim_{k\to\infty} \Delta(A_k).
\end{equation}
Indeed,  according to the proof of Theorem \ref{maxi}, the sequence
 $\{D_{A_k,t} P_{\cM_k} D_{A_k,t}\}_{k=0}^\infty$ of positive operators is a decreasing.
 Therefore, 
 $$
 Q:= \text {\rm SOT}-\lim_{k\to \infty} D_{A_k,t} P_{\cM_k} D_{A_k,t}
 $$
 exists.
Using  relation \eqref{ineg}, we infer that
\begin{equation}\label{DQD}
D_{B,t}P_{\cM_B}D_{B,t}\leq Q\leq D_{A_k,t} P_{\cM_k}D_{A_k,t}
\end{equation}
for any $k=0,1,\ldots$.
On the other hand, using  relation \eqref{dax},
 we obtain 
\begin{equation*}
\begin{split}
\left< Qx,x\right>  &\leq \|P_{\cM_k}D_{A_k,t}x\|^2\leq
\left\|D_{A_k,t} x -\sum_{j=1}^n D_{A_k,t} T_j h_j\right\|^2\\
&= \left\|D_{B,t}\left(x-\sum_{j=1}^n  T_j h_j\right)\right\|^2+ 
\left\|(I-P_k)B\left(x-\sum_{j=1}^n  T_j h_j\right)\right\|^2
\end{split}
\end{equation*}
for any $k=0,1,\ldots$, and $h_1,\ldots, h_n\in \cH$.
Since $P_k\to I$ strongly, as $k\to \infty$, 
we obtain
$$
\left< Qx,x\right>\leq \left\|D_{B,t}\left(x-\sum_{j=1}^n  T_j h_j\right)\right\|^2.
$$
Taking the infimum over $h_1,\ldots, h_n\in \cH$ we get $Q\leq D_{B,t}P_{\cM_B} D_{B,t}$,
 which together with relation \eqref{DQD} imply $Q=D_{B,t}P_{\cM_B}D_{B,t}$.
Now, it is clear that relation   \eqref{DB} holds.

\smallskip

 Does Theorem \ref{maxi} remain true if $\cV'$ is an arbitrary isometric lifting of $\cT'$ $?$
 The answer is given in what follows.
 
 \begin{theorem}\label{maxprin2}
 Let  $\cT:=[T_1~\cdots~ T_n]$, \ $T_i\in B(\cH)$, be a row isometry and 
   $\cT':=[T_1'~\cdots~ T_n']$, \ $T_i'\in B(\cH')$,  be a row contraction.
 Let    $\cW':=[W_1'~\cdots~ W_n']$, \ $V_i'\in B(\cG')$, be an isometric lifting of $\cT'$,
  and 
 let $A\in B(\cH, \cH')$  be such that  $\|A\|\leq t$  and
 $$
 AT_i=T_i'A, \quad i=1,\ldots, n.
 $$
 If 
   $B$ is an intertwining lifting of $A$ with respect to $\cW'$ and tolerance $t>0$,
   then
   \begin{equation}\label{DBtild}
   \Delta(B)\leq \Delta(\tilde{B_c})= \Delta(A),
   \end{equation}
   where  $\tilde{B_c}$ is the central intertwining lifting of $A$ with tolerance $t>0$.
   Moreover, 
   $\Delta(B) = \Delta(A)$ if and only if $\Gamma|\cM=0$, where 
  $\Gamma:\cD_{A, t}\to \cG'\ominus \cH'$ is the unique contraction   
  in the matrix decomposition
   $
  B=\left[\begin{matrix}A\\
  \Gamma D_{A, t}
  \end{matrix}
  \right]. 
   $
 \end{theorem}

\begin{proof}
Since  $B$ is an intertwining lifting of $A$ with tolerance $t>0$, we have 
the  matrix representation 
$
  B=\left[\begin{matrix}A\\
  \Gamma D_{A, t}
  \end{matrix}
  \right], 
   $
 where 
  $\Gamma:\cD_{A, t}\to \cG'\ominus \cH'$  is a contraction.
Since $\|D_{A,t}h\|= \|D_{\Gamma} D_{A,t} h\|$,  $h\in \cH$, there is a unitary operator
 $U_\Gamma:\cD_{B,t} \to \cD_\Gamma $ such that 
 \begin{equation}
 \label{UDDD}
 U_\Gamma D_{B,t}=D_{\Gamma} D_{A, t}.
 \end{equation}
Note that, 
 for any $\ell\in \cL:= \bigcap_{i=1}^n \ker T_i^*$, we have
\begin{equation*}
\begin{split}
\left< \Delta(B)\ell, \ell\right>&=
\inf \left\{ \|D_{B,t}(\ell-T_1h_1-\cdots -T_nh_n)\|^2:\ h_i\in \cH\right\}\\
&=
\inf \left\{ \|D_{\Gamma} D_{A, t}(\ell-T_1h_1-\cdots -T_nh_n)\|^2:\ h_i\in \cH\right\}\\
&\leq 
\inf \left\{ \|D_{A,t}(\ell-T_1h_1-\cdots -T_nh_n)\|^2:\ h_i\in \cH\right\}\\
&=\left< \Delta(A)\ell, \ell\right>.
\end{split}
\end{equation*}
Therefore, we obtain
\begin{equation}\label{DBDA}
\Delta(B)\leq \Delta(A).
\end{equation}
On the other hand, since  $D_{{\tilde B_c}, t}=D_{ B_c, t}$, Theorem \ref{maxi} implies 
$$
  \Delta(\tilde{B_c})=\Delta(B_c)= \Delta(A).
  $$
  
Let $
  B=\left[\begin{matrix}A\\
  \Gamma D_{A, t}
  \end{matrix}
  \right] 
   $
   be an intertwining lifting of $A$ with tolerance $t>0$, and assume that $\Gamma|\cM=0$.
   Then $D_\Gamma^2|\cM=I_\cM$, which implies $D_\Gamma|\cM=I_\cM$ and $\cM$ is
   a reducing subspace for 
   $D_\Gamma$.
   Hence, and taking into account that 
   $$
   \cN=\bigvee_{i=1}^n D_{A, t} T_i\cH,  \quad D_\Gamma \cN\subset \cN, \quad  
  \text{  and } \ D_\Gamma|\cM=I_\cM,
  $$
   we obtain
   \begin{equation*}
\begin{split}
\left< \Delta(B)\ell, \ell\right>&=
\inf \left\{ \|D_{\Gamma} D_{A, t}(\ell-T_1h_1-\cdots -T_nh_n)\|^2:\ h_i\in \cH\right\}\\
&=
\inf \left\{ \|D_{\Gamma} D_{A, t}\ell-D_\Gamma k\|:\ k\in \cN \right\}\\
&\geq
\inf \left\{ \|D_{\Gamma} D_{A, t}\ell- x\|:\ x\in \cN \right\}\\
&=\|P_\cM D_\Gamma D_{A, t}\ell\|^2= \| D_\Gamma P_\cM D_{A, t}\ell\|^2\\
&=\|   P_\cM D_{A, t}\ell\|^2= 
\left< \Delta(A)\ell, \ell\right>.
\end{split}
\end{equation*}
Therefore, $\Delta(B)\geq \Delta(A)$, which together with relation \eqref{DBDA} imply
$\Delta(B)= \Delta(A)$.

Now, let $
  B=\left[\begin{matrix}A\\
  \Gamma D_{A, t}
  \end{matrix}
  \right] 
   $
   be an intertwining lifting of $A$ with tolerance $t>0$, and assume that
   $\Delta(B)= \Delta(A)$.
Since  $\cN=\bigvee_{i=1}^n D_{A, t} T_i\cH$ and $\cD_{A,t}=\cN\oplus \cM$,
  we have 
 \begin{equation*}
\begin{split}
\left< \Delta(B)\ell, \ell\right>&=
\inf \left\{ \|D_{\Gamma} D_{A, t}(\ell-T_1h_1-\cdots -T_nh_n)\|^2:\ h_i\in \cH\right\}\\
&=  \| D_\Gamma P_\cM D_{A, t}\ell\|^2\\
&\leq \|   P_\cM D_{A, t}\ell\|^2= 
\left< \Delta(A)\ell, \ell\right>
\end{split}
\end{equation*}
for any  $\ell\in \cL$.
Since $\Delta(B)= \Delta(A)$, we infer that 
$$
\| D_\Gamma P_\cM D_{A, t}\ell\|
= \|   P_\cM D_{A, t}\ell\|, \quad \ell\in \cL.
$$
Hence, we get $\Gamma|\cM=0$. The proof is complete.
\end{proof}

\bigskip
\subsection{A permanence principle
 for the central intertwining lifting}
 \label{permanence}

 In  this section, we present a  permanence principle 
 for the central intertwining lifting. This generalizes the permanence principle
  for the Carath\' eodory
 interpolation problem  in \cite{EGL}
  (case $n=1$ ) to our multivariable setting.
 Applications of this principle  will be considered in  the next sections.

 Let  $[T_1'~\cdots ~T_n']$, \ $T_i'\in B(\cH')$,  be a row contraction and let
   $[V_1'~\cdots ~V_n']$, \ $V_i'\in B(\cK')$,
 be its minimal isometric dilation.
 Let $[T_1~\cdots ~T_n]$, \ $T_i\in B(\cH)$,  be a row  isometry and let 
 $\cM'\subseteq \cK'$ be an invariant subspace under  each 
 ${V_i'}^*$, $i=1,\ldots, n$,  such that $\cH'\subseteq \cM'$.
 Let $A\in B(\cH, \cH')$ be such that $\|A\|\leq t$ and 
 $AT_i= T_i'A$ for any $ i=1,\ldots, n$.
 An operator $Y:\cH\to \cM'$ is called partial intertwining lifting of $A$ if 
 \begin{equation}\label{partint}
 P_{\cM'} Y=A\quad \text{ \rm and }\quad YT_i=(P_{\cM'} V_i'|\cM')Y, 
 \quad i=1,\ldots, n.
 \end{equation}
 
 For example, a $k$-step intertwining lifting $A_k$ is a partial 
 intertwining lifting of $A$.

 \begin{theorem}\label{perman}
 Let $[T_1~\cdots ~T_n]$, 
 \ $T_i\in B(\cH)$, be a row isometry and let
  $[T_1'~\cdots ~T_n']$, 
 \ $T_i'\in B(\cH')$, 
  be a row contraction with minimal isometric dilation $\cV'$.  
If  $A: \cH\to \cH'$  is    
  such that $\|A\|\leq t$ and 
 $$
 AT_i= T_i'A, \quad i=1,\ldots, n,
 $$
 and $Y:\cH\to \cM'$ is a contractive partial intertwining lifting of $A$, then
 \begin{equation}\label{DEDE}
 \Delta(A)\geq \Delta(Y).
 \end{equation}
 The equality holds if and only if  $Y=P_{\cM'} B_c$, where $B_c$ 
 is the central intertwining lifting of $A$ with respect to $\cV'$ and tolerance $t$.
 Moreover, the central intertwining lifting of $P_{\cM'} B_c$ is precisely
  the central
 intertwining lifting of A. 
 \end{theorem}
 \begin{proof}
 First note that $[V_1'~\cdots ~V_n']$ is also the minimal isometric
  dilation of the row contraction 
 $[P_{\cM'}V_1'|\cM'~\cdots~P_{\cM'} V_n'|\cM']$. This follows from
  the fact that $\cH'\subset \cM'$, ${V_i'}^*\cM'\subset \cM'$ for any 
   $ i=1,\ldots, n$, 
   and the minimality condition 
 $\cK'=\bigvee_{\alpha\in \FF_n^+} V_\alpha' \cH'$.
 On the other hand, any intertwining lifting $\hat Y$ of $Y$ is also
  an intertwining lifting  of $A$. Indeed, we have 
  $V_i'\hat Y=\hat Y T_i$, $i=1,\ldots, n$, and 
  $$
  P_{\cH'}\hat Y= P_{\cH'}P_{\cM'}\hat Y= P_{\cH'} Y= A.
  $$
  Now, let $Y_c$ be the central intertwining of $Y$. Applying Theorem \ref{maxi},
  we get $\Delta (Y)=\Delta  (Y_c)$. Since $Y_c$ is an intertwining
   lifting of $A$, we can use again Theorem \ref{maxi} to obtain
   $$
   \Delta (A)\geq \Delta ( Y_c).
   $$
   Therefore, we infer  relation \eqref{DEDE}. According to Theorem \ref{maxi}, 
   the equality 
   $\Delta (A)= \Delta (Y_c)$ holds if and only if $Y_c=B_c$, where 
   $B_c$ is the 
   central intertwining of $A$. Hence, $Y=P_{\cM'} Y_c=P_{\cM'} B_c$.

 To prove the second part of the theorem, let $Y=P_{\cM'} B_c$. 
 Since $\cM'$ is an invariant subspace under each 
 $V_i^*$, $i=1,\ldots, n$, we have
 \begin{equation*}\begin{split}
 (P_{\cM'}V_i'|\cM')Y &= P_{\cM'}V_1'P_{\cM'}B_c=P_{\cM'}V_1'B_c\\
 &= P_{\cM'}B_c T_i=YT_i,
 \end{split}\end{equation*}
 for any $i=1,\ldots, n$. Hence, $Y$ is a partial intertwining lifting of $A$.
 Applying the first part of the theorem to the partial intertwining 
 lifting $Y$, we have 
 \begin{equation}\label{DDD1}
 \Delta(A)\geq \Delta(Y)=\Delta(Y_c),
 \end{equation}
 where $Y_c$ is the central intertwining lifting of $Y$.
 Since $B_c$ is an intertwining lifting of $Y$, Theorem \ref{maxi} implies 
 \begin{equation}\label{DDD2}
 \Delta(Y_c)\geq \Delta(B_c)=\Delta(A).
 \end{equation}
 Combining relation \eqref{DDD1} with \eqref{DDD2}, we get $ \Delta(Y_c)=\Delta(B_c)$. 
 Since $Y_c$ and $B_c$ are intertwining liftings of $A$, the uniqueness part
 of Theorem \ref{maxi} implies $Y_c=B_c$.
 The proof is complete.
 \end{proof}

\bigskip

  \subsection{Quasi  outer spectral factorizations}
  \label{Quasi}

 In this section,  we obtain  explicit formulas for the  quasi outer  spectral
   factor  of  the defect operator  $t^2I- B_c^* B_c$ of the central intertwining lifting $B_c$
   (see Theorem \ref{quasi}). This leads, 
   in the next section,  to concrete formulas for 
  the entropy of $B_c$   as well as to
   a maximum principle
   and  a characterization of the   central intertwining lifting
  $\tilde{B_c}$ with respect to non-minimal isometric liftings.

  Let  $\cT:=[T_1~\cdots~ T_n]$, \ $T_i\in B(\cH)$,  be a row isometry,
    $\cT':=[T_1'~\cdots ~T_n']$, \ $T_i'\in B(\cH')$,  be a row contraction,
     and 
 let 
 $A:\cH\to \cH'$  be    such that  $\|A\|\leq t$ and 
 $$
 A T_i= T_i'A, \quad i=1,\ldots, n.
 $$
 Define the operator $X_A\in B(\oplus_{i=1}^n \cH, \cH)$ by
 \begin{equation}\label{XA}
 X^A:=D_{A,t}^2 [ T_1 ~\cdots ~ T_n ]
  [T_i^* D_{A,t}^2T_j]_{n\times n},
 \end{equation}
 and let $X^A=[X^A_1 ~\cdots ~X^A_n]$  be its matrix representation
   with $X^A_i\in B(\cH)$, ~$i=1,\ldots, n$.
 As is  \cite{Po-nehari},  one can prove that if $\|A\|<t$,
  then the central intertwining  lifting of $A$  satisfies the equation
  \begin{equation}
  \label{BC*}
  B_ch=Ah\oplus \sum_{\sigma\in \FF_n^+} e_\sigma \otimes D_{\cT'}
  (\oplus_{i=1}^n A) (X^A_\sigma )^*h
  \end{equation}
 for any $h\in \cH$.

 \begin{lemma}\label{lem1}
 If $\|A\|<t$, then the following statements hold:
 \begin{enumerate}
 \item[(i)]
 $\cM:= \cD_{A,t}\ominus \cN=D_{A,t}^{-1} \cL,$ where 
 $\cL:=\bigcap_{i=1}^n \ker T_i^*$ and $\cN$ is defined by \eqref{nmk} $($when $k=0);$
 \item[(ii)] $B_c D_{A,t}^{-2}\ell= AD_{A,t}^{-2}\ell$ for any $\ell\in \cL;$
 \item[(iii)] $E_i= D_{A,t} (X^A_i)^* D_{A,t}^{-1}$,  where $E_i:= P_i WP_\cN$
 for 
 any $ i=1,\ldots, n$, and the operator $W$ is defined by \eqref{wk} $($when $k=0)$.
 \end{enumerate}
 \end{lemma}
 
 \begin{proof}
 To prove (i), note that $h\in \cM$ if and only if $h\perp \bigvee_{i=1}^n D_{A,t} T_i \cH$,
  which is equivalent to $D_{A,t}h\in \bigcap_{i=1}^n \ker T_i^*$. 
  Hence, $h\in D_{A,t}^{-1} \cL$. 
  The statement (ii) follows from relation \eqref {BC*} and the fact that, for any 
  $\ell\in \bigcap_{i=1}^n \ker T_i^*$, 
  $$
  (X^A)^* D_{A,t}^{-2} \ell= 
  [T_j^* D_{A,t}^2T_i]_{n\times n} \left[\begin{matrix} T_1^*  \\ \vdots\\ T_n^* \end{matrix}\right]
   D_{A,t}^2 (D_{A,t}^{-2}\ell)=0.
   $$
   The proof of the  similarity relation (iii)  is exactly the same as that from
    \cite{Po-nehari} for the case $t=1$. We shall omit it.
 \end{proof}
 
 We mention that  part (i) can be used together with Theorem \ref{main3} to give 
 another proof for  part (ii) of Lemma \ref{lem1}.
 We recall that $\cT:=[T_1~\cdots~ T_n]$, ~$T_i\in B(\cH)$,  is a $C_0$-row contraction if 
 $$
 \lim_{k\to\infty} \sum_{|\alpha|=k} \|T_\alpha^* h \|^2=0
 $$
 for any $h\in \cH$.
 If $\cT= [S_1\otimes I_\cM~\cdots ~S_n\otimes I_\cM]$ for
  some Hilbert space $\cM$, then $\cT$ is called
  orthogonal shift of multiplicity ~$\dim \cM$.

 \begin{lemma}\label{lem2}
  Let  $\cT:=[T_1~\cdots~ T_n]$, \ $T_i\in B(\cH)$,  be a row isometry,
    $\cT':=[T_1'~\cdots~ T_n']$, \ $T_i'\in B(\cH')$,  be a row contraction, and 
 let 
 $A:\cH\to \cH'$  be  such that  $\|A\|\leq t$ and 
 $$
 A T_i= T_i'A, \quad i=1,\ldots, n.
 $$
 Then there exists a unique positive operator $Q\in B(\cD_{A,t})$ such that 
 \begin{equation}\label{Q}
 \|QD_{A,t} h\|^2=\lim_{k\to \infty} \sum_{|\sigma|=k} \|E_\sigma D_{A,t} h\|^2, \quad h\in \cH.
 \end{equation}
 Moreover, if $\cT$ is an orthogonal shift and $\|A\|<t$, then $Q=0$ and 
 $[X^A_1 ~\cdots ~X^A_n]$ is a $C_0$-row contraction.
 \end{lemma}
 
 \begin{proof}
 Note that if $E_i:= P_i WP_\cN$, \ $i=1,\ldots, n$,
 then $[E_1^*~\cdots ~ E_n^*]$ is a row contraction from $\oplus_{i=1}^n \cD_{A,t}$
  to $\cD_{A,t}$.
 Since $\{\sum_{|\alpha|=k} E_\alpha^* E_\alpha \}_{k=1}^\infty$ is a decreasing 
 sequence of contractions,  we can define  the  operator
 $$
 Q:=\left(\text{\rm SOT}-\lim_{k\to \infty} \sum_{|\alpha|=k} E_\alpha^*
  E_\alpha\right)^{1/2}.
  $$
   Hence, relation \eqref{Q}
 holds.
 Now, assume that $\cT$ is an orthogonal shift and $\|A\|<1$.
 Since, for any $i=1,\ldots, n$, 
 \begin{equation}\label{PWP}
 P_j WP_\cN D_{A,t} T_i h= \delta_{ij}D_{A,t}h, \quad h\in \cH,
 \end{equation}
we deduce
\begin{equation}\label{EDA}
E_{\tilde{\omega}} D_{A,t} T_\sigma h= \delta_{\sigma, \omega} D_{A,t} h,
\end{equation} 
for any  $\sigma, \omega\in \FF_N^+$.
According to Lemma 3.3. from \cite{Po-nehari}, the closed linear
 span $\bigvee_{\sigma\in \FF_n^+} D_{A,t} T_\sigma D_{A,t}^{-2}\cL $ is  equal to 
 $\cH$, where
 $\cL:=\bigcap_{i=1}^n \ker T_i^*$. Taking into account that $\cD_{A,t}= \cN\oplus \cM$ and
  $\cM= D_{A,t}^{-1} \cL$, we have
 $$
 \sum_{j=1}^n \sum_{|\alpha|=k} \|E_j E_\alpha D_{A,t} T_\sigma D_{A,t}^{-2}\ell\|^2=
  \sum_{j=1}^n \|E_j D_{A,t}^{-1} \ell\|^2=0.
 $$
 This shows that $[E_1^*~\cdots ~ E_n^*]$ is a $C_0$-row contraction.
 Using Lemma \ref {lem1} part (iii), we conclude that $X^A$ is a a $C_0$-row contraction.
 \end{proof}

  \begin{lemma}\label{lem3}
 Let $\cY_Q:=\overline{\text{\rm range} ~Q}$ and define the operators 
 $V_i\in B(\cY_Q)$, 
 \ $i=1,\ldots, n$, by setting
 \begin{equation}\label{Vi}
 V_iQD_{A,t} h:= QD_{A,t}T_ih, \qquad h\in \cH.
 \end{equation}
 Then  $V_1,\ldots, V_n$ are isometries with orthogonal ranges and 
 $$
 V_1V_1^*+\cdots +V_n V_n^*=I_{\cY_Q}.
 $$
 \end{lemma}

 \begin{proof}
 Using Lemma \ref{lem2} and relation \eqref{EDA}, we have
 \begin{equation*}
 \begin{split}
 \|QD_{A,t} T_ih\|^2&= 
 \lim_{k\to \infty} \sum_{j=1}^n \sum_{|\alpha|=k-1} \| E_\alpha E_j D_{A,t} T_i h\|^2\\
 &=  \lim_{k\to \infty}  \sum_{|\alpha|=k-1} \| E_\alpha   D_{A,t}h\|^2= 
 \|QD_{A,t} h\|^2
 \end{split}
 \end{equation*}
 for any $h\in \cH$ and $i=1,\ldots, n$.
 Similar computations show that 
 $$
 \left< QD_{A,t} T_i h, QD_{A,t} T_j h\right>=
 \text{\rm SOT}-\lim_{k\to \infty} \left<\sum_{|\alpha|=k} E_\alpha^*
  E_\alpha D_{A,t} T_ih, D_{A,t} T_j h\right>=0
  $$
  for any $i\neq j$. Therefore, $\{V_i\}_{i=1}^n $ are isometries with orthogonal ranges.
  On the other hand, since $Q\geq 0$, $\cN= \bigvee_{i=1}^n D_{A,t} T_i\cH$, and 
  the range of $Q$ is included in $\cN$, we infer that
  $$
  \bigvee_{i=1}^n QD_{A,t} T_i\cH= \overline{ Q\cD_{A,t}}.
  $$
  Now, it is clear that  
   $\bigvee_{i=1}^n QD_{A,t} T_i\cH $ is dense in $\cY_Q$ and 
 $$
 V_1V_1^*+\cdots +V_n V_n^*=I_{\cY_Q}.
 $$
 \end{proof}
 
 In the next theorem, we prove that the defect operator 
  $t^2I-B_c^* B_c $ admits a quasi outer spectral factor, i.e.,
  $ t^2I-B_c^* B_c =Z^* Z$ for some operator $Z$ with dense range. 
 
 \begin{theorem}\label{quasi}
  Let  $\cT:=[T_1~\cdots~ T_n]$, \ $T_i\in B(\cH)$,  be a row isometry,
    $\cT':=[T_1'~\cdots~ T_n']$, \ $T_i'\in B(\cH')$,  be a row contraction, and 
 let 
 $A:\cH\to \cH'$  be    such that  $\|A\|\leq t$ and 
 $$
 AT_i= T_i'A, \quad i=1,\ldots, n.
 $$
 Let $Z:\cH\to \cY_Q\oplus [F^2(H_n)\otimes \cM]$ be defined by
 $$
 Zh:= QD_{A,t}h\oplus \sum_{\sigma\in \FF_n^+} e_\sigma 
 \otimes P_\cM E_{\tilde \sigma} D_{A,t} h,
 $$
 where $\cM:= \cD_{A,t} \ominus \cN$.
 Then the operator $Z$ has dense range and 
 $$
 t^2I-B_c^* B_c =Z^* Z,
 $$
 where $B_c$ is the central intertwining lifting of $A$ with respect to the minimal isometric dilation of 
 $\cT'$ and tolerance $t$.
 Moreover, for any $i=1,\ldots, n$, we have
 \begin{equation}\label{ZTi}
 ZT_i=\left[\begin{matrix}
 V_i& 0\\
 0& S_i\otimes I_\cM
 \end{matrix}\right] Z, 
 \quad  i=1,\ldots, n,
 \end{equation}
 where the isometries $V_i\in B(\cY_Q)$ are defined  by relation \eqref{Vi}.
 \end{theorem}
  
 \begin{proof}
 As in the proof of Theorem \ref{main3}, for any $h\in \cH$, we have
 \begin{equation*}\begin{split}
 \|D_{B_c, t}h\|^2&= \|D_{A,t}h\|^2- \|P_\cN D_{A,t}h\|^2 \\
 & {}\qquad +\lim_{k\to \infty} 
 \left( \sum_{1\leq |\omega|\leq k} \| E_{\tilde \omega} D_{A,t} h\|^2-
 \sum_{1\leq |\sigma|\leq k-1} \| P_\cN E_{\tilde \sigma} D_{A,t} h\|^2
 \right)\\
 &=
 \|P_\cM D_{A,t}h\|^2 + \lim_{k\to \infty}
  \sum_{|\alpha|=k} \| E_{\tilde \alpha} D_{A,t} h\|^2 + \sum_{|\sigma|\geq 1}
   \|P_\cM E_{\tilde \sigma}D_{A,t}h\|^2.
 \end{split}
 \end{equation*}
 Hence, and using Lemma \ref{lem2}, we get
 \begin{equation*}\begin{split}
 \|D_{B_c,t}h\|^2&=\|QD_{A,t}h\|^2+ \sum_{\sigma \in \FF_n^+}
   \|P_\cM E_{\tilde \sigma}D_{A,t}h\|^2\\
   &= \|Zh\|^2
 \end{split}
 \end{equation*}
 for any $h\in \cH$. Therefore,  we have $t^2I-B_c^* B_c= Z^* Z$.
 
 Now, we prove the equality \eqref{ZTi}. 
 Using Lemma \ref{lem3} and relation \eqref{EDA}, we obtain
 \begin{equation*}\begin{split}
 \left[\begin{matrix}
 V_i& 0\\
 0& S_i\otimes I_\cM
 \end{matrix}\right] Zh&= 
 V_iQD_{A,t}h\oplus \sum_{\sigma\in \FF_n^+} e_{g_i\sigma} 
 \otimes P_\cM E_{\tilde \sigma} D_{A,t} h\\
 &=
 QD_{A,t} T_ih\oplus \sum_{\sigma\in \FF_n^+} e_{g_i\sigma} 
 \otimes P_\cM E_{\tilde \sigma}E_i D_{A,t} T_i h\\
 &=
 QD_{A,t}T_ih\oplus \sum_{\sigma\in \FF_n^+}\sum_{j=1}^n  e_{g_j\sigma} 
 \otimes P_\cM E_{\tilde \sigma}E_j D_{A,t} T_i h\\
 &=
 QD_{A,t}T_ih\oplus \sum_{\omega\in \FF_n^+} e_{ \omega} 
 \otimes P_\cM E_{\tilde \omega} D_{A,t} T_i h\\
 &= ZT_ih
 \end{split}
 \end{equation*}
 for any $h\in \cH$ and $i=1,\ldots, n$.
 
 It remains to show that $Z$ has dense range.
 Note that $\cY_Q\subseteq \cN$ and let $f\oplus \varphi$ be 
 in $\cY_Q\oplus [F^2(H_n)\otimes \cM]$ such that  $f\oplus \varphi\perp Z\cH$. 
 Consider the representation 
 $$
 \varphi= \sum_{\alpha\in \FF_n^+} e_\alpha\otimes h_\alpha, \quad h_\alpha\in \cM.
 $$
 Using the definition of $Z$, we have
 \begin{equation*}
 \begin{split}
 0=\left< f\oplus \varphi, Zh \right>&=\left< f, QD_{A,t}h \right>+
 \sum_{\alpha\in \FF_n^+} \left< h_\alpha, P_\cM E_{\tilde \alpha} D_{A,t} h \right>\\
 &=
 \left<D_{A,t} Qf, h\right>+ \sum_{\alpha\in \FF_n^+} \left<D_{A,t}E_{\tilde \alpha}^* 
 h_\alpha, h\right>,
 \end{split}
 \end{equation*}
 for any $h\in \cH$.
Therefore, 
$$
D_{A,t}\left[ Qf+ \sum_{\alpha\in \FF_n^+} E_{\tilde \alpha}^*  h_\alpha\right]=0.
$$
Since $D_{A,t}$ is a one-to-one operator on $\cD_{A,t}=\cN\oplus \cM$ 
and  
 the range
 of $Q$ is included in $\cN$, it follows that 
 $$
 h_{g_0}+Qf+ \sum_{|\alpha|\geq 1} E_{\tilde \alpha}^*  h_\alpha=0.
 $$
 Since $h_{g_0}\in \cM$ and the other summands are in $\cN$, we get 
 $h_{g_0}=0$ and
 \begin{equation}\label{Qf}
 Qf+ \sum_{|\alpha|\geq 1} E_{\tilde \alpha}^*  h_\alpha=0.
 \end{equation}
 Now, note that, for any $x\in \cD_{A,t}$,   we have
 \begin{equation} \label{QP}
 \sum_{i=1}^n \|Q(P_iW P_\cN)x\|^2=
 \lim_{k\to \infty} \sum_{|\sigma|=k} \sum_{i=1}^n \|E_\sigma E_i x\|^2
 = \|Qx\|^2.
 \end{equation}
 Define the operators $W_i\in B(\cY_Q)$, \ $i=1,\ldots, n$,  by setting
 \begin{equation}\label{Wi}
 W_iQh:= QE_i h, \quad h\in \cH.
 \end{equation}
 According to  relation \eqref{QP}, we have
 \begin{equation}\label{W*}
 W_1^*W_1+\cdots + W_n^* W_n= I_{\cY_Q}.
 \end{equation}
 Using  relations \eqref{Wi} and \eqref{W*}, we obtain 
 $QW_i^* W_i= E_i^* Q E_i$,  for any  $i=1,\ldots, n$, and  
   $Q= \sum_{i=1}^n E_i^* Q W_i$, which
  together with relation \eqref{Qf} imply
  $$
  \sum_{i=1}^n E_i^* Q W_if+ \sum_{|\alpha|\geq 1} E_{\tilde \alpha}^*  h_\alpha=0.
  $$
  This equality can be written as
  \begin{equation}\label{EI*}
  \sum_{i=1}^n E_i^* \left(Q W_if+ h_{g_i}+  \sum_{|\sigma|\geq 2} E_{\tilde \sigma}^* 
   h_{g_i\sigma}\right)=0.
  \end{equation}
 Note that the subspace 
 $$
 \{E_1k\oplus E_2 k \oplus \cdots \oplus E_n k:\quad k\in \cD_{A,t}\}
 $$
  is dense 
 in $\oplus_{i=1}^n \cD_{A,t}$.
 Indeed, this is due to the equations 
 $$E_i D_{A,t} T_j h =\delta_{ij} D_{A,t} h, \quad  i,j=1,\ldots, n.
 $$
 Therefore, the operator $[E_1^*~\cdots ~ E_n^*]$ is injective and relation \eqref{EI*}
 implies
 $$
 Q W_if+ h_{g_i}+  \sum_{|\sigma|\geq 2} E_{\tilde \sigma}^* 
   h_{g_i\sigma} =0, \quad  i=1,\ldots, n.
   $$
 Since $h_{g_i}\in \cM$ and the other summands are in $\cN$, we obtain $h_{g_i}=0$,  for any 
 $i=1,\ldots, n$.
 Continuing this process, one can prove that $h_\alpha=0$ for any $\alpha\in \FF_n^+$.
 Therefore,
 relation \eqref{Qf} implies $Qf=0$. Since $f$ is in the range of $Q$, we get $f=0$.
 The proof is complete.
  \end{proof}

 \begin{corollary}\label{quasi2}
  Let  $\cT:=[T_1~\cdots~ T_n]$, \ $T_i\in B(\cH)$,  be a row isometry and 
  $\cL:= \bigcap_{i=1}^n \ker T_i^*$. 
   Let $\cT':=[T_1'~\cdots~ T_n']$, \ $T_i'\in B(\cH')$,  be a row contraction, and 
 let 
 $A:\cH\to \cH'$  be    such that  $\|A\|< t$ and 
 $$
 AT_i= T_i'A, \quad i=1,\ldots, n.
 $$
 Let $Z_\cL:\cH\to \cY_Q\oplus [F^2(H_n)\otimes \cL]$ be defined by
 $$
 Z_\cL h:= QD_{A,t}h\oplus \sum_{\sigma\in \FF_n^+} e_\sigma 
 \otimes N_A P_\cL (X^A_\sigma)^* h,
 $$
 where $ N_A:= (P_\cL D_{A,t}^{-2}|\cL)^{-1/2}$ and $X^A$ is defined by relation  \eqref{XA}.
 Then $Z_\cL$ has dense range and 
 $$
 t^2I-B_c^* B_c =Z_\cL^* Z_\cL,
 $$
 where $B_c$ is the central intertwining lifting of $A$ with tolerance $t$.
 Moreover, for any $i=1,\ldots, n$, we have
 \begin{equation}\label{ZLT}
 Z_\cL T_i=\left[\begin{matrix}
V_i& 0\\
 0&  S_i\otimes I_\cL
 \end{matrix}\right] Z_\cL,\quad  i=1,\ldots, n,
 \end{equation}
 where the isometries $V_i\in B(\cY_Q)$ 
  are defined  by relation \eqref{Vi}.
 \end{corollary}
  
 \begin{proof}
 According to Lemma \ref{lem1}, we have $\cM=\cD_{A,t}\ominus \cN= D_{A,t}^{-1} \cL$, where 
 $\cL= \bigcap_{i=1}^n \ker T_i^*$.
 Note that the operator $X\in B(\cL, \cH)$ defined by $X:= D_{A,t}^{-1}|\cL$ is injective and
  has range equal to $\cM$.
  Hence, one can easily see that
  \begin{equation}\label{PM}
  P_\cM= X (X^* X)^{-1} X^*= (D_{A,t}^{-1}|\cL)(P_\cL D_{A,t}^{-2} |\cL)^{-1} P_\cL D_{A,t}^{-1}.
  \end{equation}
  Using again  part (iii) of Lemma \ref{lem1}  and relation \eqref{PM},  we obtain
  \begin{equation*}\begin{split}
  \sum_{\sigma\in \FF_n^+} e_\sigma \otimes P_\cM E_{\tilde \sigma} D_{A,t}h&= 
  \sum_{\sigma\in \FF_n^+} e_\sigma \otimes P_\cM D_{A,t} (X_\sigma^A)^*h\\
  &=\sum_{\sigma\in \FF_n^+} e_\sigma \otimes  \Omega N_A P_\cL (X_\sigma^A)^*h\\
  &= (I_{F^2(H_n)}\otimes \Omega)
  \left( \sum_{\sigma\in \FF_n^+} e_\sigma \otimes   N_A P_\cL (X_\sigma^A)^*h\right),
  \end{split}
  \end{equation*}
 where  
 $$\Omega:= (D_{A,t}^{-1}|\cL)(P_\cL D_{A,t}^{-2} |\cL)^{-1/2} \ \text{ and } \  
 N_A:=(P_\cL D_{A,t}^{-2} |\cL)^{-1/2}.
 $$
  Note that $\Omega$ is a
  unitary operator from $\cL$ onto $\cM$. 
  Using Theorem \ref{quasi}, we deduce
  $$
  Z=\left[ \begin{matrix}
 I_{\cY_Q}& 0\\
 0&  I_{F^2(H_n)}\otimes \Omega 
 \end{matrix} \right] Z_\cL,
 $$ 
 which completes the proof.
 \end{proof}

 \smallskip

 \begin{corollary} \label{strict-form}
 Let  $[T_1'~\cdots~ T_n']$, 
 \ $T_i'\in B(\cH')$, 
  be a row contraction,
 let $A:F^2(H_n)\otimes \cE\to \cH'$  be  
  such that  $\|A\|<t$ and 
 $$
 A(S_i\otimes I_\cE)= T_i'A, \quad i=1,\ldots, n.
 $$
 If  $B_c$ is   the central intertwining 
  lifting of $A$ with $\|B_c\|\leq t$,
 then  the multi-Toeplitz operator $t^2I-B_c^* B_c$ admits 
  an outer spectral factorization
  $$
  t^2I-B_c^* B_c= M_\varphi^*M_\varphi,
  $$
  where $ M_\varphi\in R_n^\infty\bar{\otimes} B(\cE)$ is an outer multi-analytic
   operator operator with the symbol $\varphi$
  given by
 \begin{equation}\label{vp}
 \varphi ( h):= \sum_{\sigma\in \FF_n^+} e_\sigma \otimes
  \Delta(A)^{1/2} P_\cE (X_\sigma^A)^*h, \qquad h\in \cE,
 \end{equation}
 where $\Delta (A)= P_\cE(t^2I-A^* A)^{-1}|\cE$. 
 \end{corollary}
 
 \begin{proof}
 Using Lemma \ref{lem2} when $T_i= S_i\otimes I_\cE$, \ $i=1,\ldots, n$, we get $Q=0$.
  In this particular case,  Corollary \ref{quasi2}  implies
 $$
 Z_\cE (S_i\otimes I_\cE) =  (S_i\otimes I_\cE) Z_\cE, \ i=1,\ldots, n.
 $$
 Hence, $Z_\cE= M_\varphi$ for some $M_\varphi\in R_n ^\infty\bar\otimes B(\cE)$.
 Since $Z_\cE$ has dense range in $F^2(H_n)\otimes \cE$, the operator $M_\varphi$ is outer.
 Using again Corollary \ref{quasi2}, one can easily  complete the proof.
 \end{proof}

 We remark that  the multi-analytic operator $M_\varphi$ in Corollary \ref{strict-form}
 satisfies the equation 
 $$
 M_\varphi = (I_{F^2(H_n)} \otimes \Delta (A)^{1/2} P_\cE) 
 (I-R_1\otimes (X_1^A)^*-\cdots -
 R_n\otimes (X_n^A)^*)^{-1}.
 $$
 Using  relation \eqref{vp}, we can deduce this formula  as in the proof
  of Remark \ref{state-sol}.

 \bigskip
 
 \subsection{Noncommutative commutant lifting theorem and the maximal
  entropy solution}\label{max-entro-sol}
 
  In this section, we prove that  the 
 central intertwining lifting is the maximal entropy solution for 
 the noncommutative commutant lifting theorem when 
 $\cT:=[S_1\otimes I_\cE~\cdots ~ S_n\otimes I_\cE]$ with $\dim \cE<\infty$
 (see Theorem \ref{inter-entro}). 
 Based on several results of this paper, we are  led to concrete formulas for 
  the entropy of $B_c$ (see Theorem \ref{strict}) and,
   under a certain condition of stability, to
   a maximum principle
   and  a characterization (in terms of  entropy) of the   central intertwining lifting
  $\tilde{B_c}$ with respect to non-minimal isometric liftings
   (see Theorem \ref{maxprin3} and Corollary \ref{en-gen}).

 Let  $[T_1'~\cdots~ T_n']$, \ $T_i'\in B(\cH')$,  be a row contraction and let
   $[V_1'~\cdots~ V_n']$, \ $V_i'\in B(\cK')$,
 be its minimal isometric dilation.
 Let   $A:F^2(H_n)\otimes \cE\to \cH'$  be   such that  $\|A\|\leq t$ and 
 $$
 A(S_i\otimes I_\cE)= T_i'A, \quad i=1,\ldots, n,
 $$
 and 
  let $B:F^2(H_n)\otimes \cE\to \cK'$  be an  intertwining
  lifting of $A$ satisfying  $\|B\|\leq t$. Then     
  $B$ is a generalized  multiplier 
  and  it makes sense to define its  prediction entropy  
   as in   Section \ref{extreme points}, i.e.,
   $E(B):= e(D_{B,t}^2)$, where $D_{B,t}^2=t^2I-B^* B$.
 
 In what follows, we show that the central intertwining lifting of $A$  with
  tolerance $t>0$ is 
 the  maximal 
 entropy intertwining lifting of $A$. More precisely,
  we can prove the following result.

 \begin{theorem}\label{inter-entro} Let  $[T_1'~\cdots~ T_n']$, 
 \ $T_i'\in B(\cH')$, 
  be a row contraction  with minimal isometric dilation $\cV'$, and
 let $A:F^2(H_n)\otimes \cE\to \cH'$  be    such that  $\|A\|\leq t$ and 
 $$
 A(S_i\otimes I_\cE)= T_i'A, \quad i=1,\ldots, n.
 $$
  If $\dim \cE<\infty$ and $B_c$ is the central intertwining 
  lifting of $A$ with respect to $\cV'$ and tolerance $t>0$, then
  $$
  E(B_c)\geq  E(B),
  $$
  for any   intertwining lifting  $B$ of $A$.
   Moreover, if  the entropy $ E(B_c)>-\infty$,
   then 
   $$
    E(B_c)= E(B)
    $$
     if and only if $B_c=B$.
 \end{theorem}
 
 \begin{proof} Let
 $B:F^2(H_n)\otimes \cE\to \cK'$ be   an  intertwining
  lifting of $A$ satisfying $\|B\|\leq t$.
  Using  relations \eqref{EB}, \eqref{ET}, \eqref{DTo}, 
  and \eqref{Dform} in our setting, we deduce
  \begin{equation*}\begin{split}
  E(B)&=e(D_{B,t}^2)=\ln \det \Delta_{D_{B,t}^2}\\
  &=\ln \det [P_\cE D_{B,t}P_{\cM} D_{B,t}|\cE]\\
  &= \ln \det \Delta(B),
  \end{split}
  \end{equation*}
  where 
  $$
  \cM:= \overline{D_{B,t}(F^2(H_n)\otimes \cE)}\ominus 
  \bigvee_{i=1}^n D_{B,t}(S_i\otimes I_\cE)(F^2(H_n)\otimes \cE).
  $$
  Similarly, we obtain $E(B_c)=\ln \det \Delta(B_c)$, where 
   $B_c$ is the central intertwining 
  lifting of $A$.
  According to Theorem \ref{maxi}, we have  $\Delta(B_c)\geq \Delta (B)$.
  Therefore, we get 
  $$
  E(B_c)=\ln \det \Delta(B_c)\geq \ln \det \Delta(B)=E(B).
  $$
  Now, assume that $E(B_c)>-\infty$ and $E(B_c)=E(B)$. This implies
  $$
  \det \Delta(B_c)=\det \Delta(B)\neq 0.
  $$
  Since $\Delta(B_c)$ and $\Delta(B)$ are strictly positive operators with 
  $\Delta(B_c)\geq \Delta (B)$, we infer that $\det \Delta(B_c)=\det \Delta(B)$
  if and only if 
  $\Delta(B_c)=\Delta(B)$. Using the uniqueness part of Theorem \ref{maxi},
   we get $B=B_c$,
   which completes the
  proof.
 \end{proof}

 \begin{theorem} \label{strict}
 Let  $[T_1'~\cdots~ T_n']$, 
 \ $T_i'\in B(\cH')$, 
  be a row contraction with minimal isometric dilation $\cV'$ and
 let $A:F^2(H_n)\otimes \cE\to \cH'$  be  
  such that  $\|A\|<t$ and 
 $$
 A(S_i\otimes I_\cE)= T_i'A, \quad i=1,\ldots, n.
 $$
 If   $B_c$ is  the central intertwining 
  lifting of $A$ with  respect to $\cV'$   and tolerance $t$,
  then the multi-Toeplitz operator $t^2I-B_c^* B_c$ admits 
  a square outer spectral factorization
  $$
  t^2I-B_c^* B_c= M_\varphi^*M_\varphi,
  $$
  where $ M_\varphi\in R_n^\infty\bar{\otimes} B(\cE)$ is an outer operator
  given by
  $$
   M_\varphi:= [I\otimes (P_\cE \psi)^{1/2}] M_\psi^{-1},
  $$
  and the symbol $\psi$  of the multi-analytic operator 
  $ M_\psi$ is given by $\psi:= (t^2I-A^*A)^{-1}|\cE$.
  
  Moreover, if $\dim \cE<\infty$, then
  \begin{equation}
  \label{EBc}
  E(B_c)=-\ln\det [P_\cE(t^2I-A^* A)^{-1}|\cE].
  \end{equation}
 \end{theorem}
 
 \begin{proof}
 Using  the proof of Corollary \ref{quasi2} and Corollary \ref{strict-form}, we  deduce that 
 $$
 (I\otimes \Omega^*)Z= Z_\cE= M_\varphi,
 $$
 and $M_\varphi$ is the outer spectral factor of ~$t^2I-B_c^* B_c$~ defined by \eqref{vp}. 
 Recall that $ \cD_{A,t}\ominus \cN =\cM=D_{A,t}^{-1} \cE$ and 
  $$
  \Omega:= (D_{A,t}^{-1}|\cE)(P_\cE D_{A,t}^{-2} |\cE)^{-1/2}
  $$ 
   is a
  unitary operator from $\cE$ onto $\cM$. 
 Hence, we have $P_\cN D_{A,t}^{-1} |\cH=0$ and
 \begin{equation}
 \begin{split}
 M_\varphi D_{A,t}^{-2}h&=  (I\otimes \Omega^*)ZD_{A,t}^{-2}h= \Omega^* P_\cM D_{A,t}^{-1} h\\
 &=
 \Omega^* D_{A,t}^{-1}h= (P_\cE D_{A,t}^{-2} |\cE)^{1/2}h=1\otimes
  (P_\cE \psi)^{1/2}h
 \end{split}
 \end{equation}
 for any $h\in \cE$.
 Applying Theorem \ref{outer} to the strictly positive multi-Toeplitz operator $T:=t^2I-A^* A$,
 we infer that $M_\psi$ is an invertible multi-analytic operator.
 Now, note that
 $$
 M_\varphi M_\psi (e_\alpha\otimes h)= (S_\alpha \otimes I_\cE) M_\varphi D_{A,t}^{-2}=
 e_\alpha \otimes (P_\cE \psi)^{1/2}h
 $$
 for any $h\in \cE$ and $\alpha\in \FF_n^+$.
 Hence,
 $$
 M_\varphi= [I_{F^2(H_n)}\otimes (P_\cE \Psi)^{1/2}] M_\psi^{-1}.
 $$
 To complete the proof, note that relation \eqref{vp}
 implies
 \begin{equation}\label{vava}
 \varphi(0)^* \varphi(0)=\Delta(A)= P_\cE(t^2I-A^*A)^{-1}|\cE.
 \end{equation}
 Since $M_\varphi$ is outer, Theorem \ref{outer2-gen} implies
 \begin{equation}\label{Ee}
 E(B_c)= e(M_\varphi^* M_\varphi)= \ln \det [\varphi(0)^* \varphi(0)].
 \end{equation}
 Combining relations \eqref{vava} and \eqref{Ee}, we get \eqref{EBc}. The proof is complete.
 \end{proof}

If $B$  be  an arbitrary intertwining 
lifting of $A$ with $\|B\|<t$, we can obtain  a  result  similar to that of  
Theorem \ref{strict} by applying Theorem \ref{sub-fac} and Corollary \ref{entro-fact} to the 
multi-Toeplitz operator $T:= t^2I-B^* B$.

 \begin{corollary} \label{strict2}
 Let  $[T_1'~\cdots~ T_n']$, 
 \ $T_i'\in B(\cH')$, 
  be a row contraction with minimal isometric dilation $\cV'$ and 
 let $A:F^2(H_n)\otimes \cE\to \cH'$  be 
  such that  $\|A\|<t$ and 
 $$
 A(S_i\otimes I_\cE)= T_i'A, \quad i=1,\ldots, n.
 $$
 If   $B$ is    an  intertwining 
  lifting of $A$  with  respect to $\cV'$ and tolerance $t$,
  then  the multi-Toeplitz operator $t^2I-B^* B$ admits 
  a square outer spectral factorization
  $$
  t^2I-B^* B= M_\chi^*M_\chi,
  $$
  where $ M_\chi\in R_n^\infty\bar{\otimes} B(\cE)$ is an outer operator
  given by
  $$
   M_\chi:= [I\otimes (P_\cE g)^{1/2}] M_g^{-1},
  $$
  and the symbol $g$  of the multi-analytic operator 
  $ M_g$ is given by $g:= (t^2I-B^*B)^{-1}|\cE$.
  
  Moreover, if $\dim \cE<\infty$, then
  \begin{equation*}
  E(B)=-\ln\det [P_\cE(t^2I-B^* B)^{-1}|\cE].
  \end{equation*}
 \end{corollary}

\smallskip
In what follows, we show that, under a certain condition of stability, there 
is a maximal principle
for the noncommutative  commutant lifting theorem and a characterization of 
the central  intertwining lifting
$\tilde B_c$ of $A$ with respect to non-minimal isometric liftings. The results will
 be very 
useful  in the next sections.

 \begin{theorem}\label{maxprin3}
 Let  $\cT:=[T_1~\cdots~ T_n]$, \ $T_i\in B(\cH)$, be a row isometry and 
   $\cT':=[T_1'~\cdots~ T_n']$, \ $T_i'\in B(\cH')$,  be a row contraction.
 Let    $\cW':=[W_1'~\cdots~ W_n']$, \ $V_i'\in B(\cG')$, be an isometric lifting of $\cT'$,
  and 
 let $A\in B(\cH, \cH')$  be  such that  $\|A\|\leq t$ and 
 $$
 AT_i=T_i'A, \quad i=1,\ldots, n.
 $$
 If 
   $B$ is an intertwining lifting of $A$ with respect to $\cW'$ and tolerance $t>0$,
   then
   \begin{equation}\label{DBtild2}
   \Delta(B)\leq \Delta(\tilde{B_c})= \Delta(A),
   \end{equation}
   where  $\tilde{B_c}$ is the central intertwining lifting of $A$ with tolerance $t>0$.
   Assume that $[E_1^* ~\cdots ~ E_n^*]$ is a $C_0$-row contraction,
  where  $E_i:= P_i W P_\cN\in B(\cD_{A,t})$, for any $i=1,\ldots, n$.
   Then 
   $\Delta(B) = \Delta(A)$ if and only if $B=\tilde{B_c}$.
   
   In particular, if $\|A\|<t$ and $[T_1~ \cdots ~ T_n]$ is unitarily equivalent
    to an orthogonal 
   shift
   $[S_1\otimes I~ \cdots ~ S_n\otimes I]$, then $\Delta(B) = \Delta(A)$
   if and only if $B=\tilde{B_c}$.
 \end{theorem}

\begin{proof} The inequality \eqref{DBtild2} was proved in Theorem \ref{maxprin2}.
Assume now that $[E_1^* ~\cdots ~ E_n^*]$ is a $C_0$-row contraction
and 
 $\Delta(B) = \Delta(A)$.  Recall that the operator
$W_i'$ admits a reducing decomposition $W_i'= V_i'\oplus U_i'$, \ $i=1,\ldots, n$, 
 of \ $\cG'= \cK'\oplus \cE'$
such that $\cV':=[V_1'~\cdots ~ V_n']$ is the minimal isometric dilation 
of $\cT'$ on $\cK'$, which can be identified
with $\cH'\oplus [F^2(H_n)\otimes \cD']$.
Since  $P_{\cH'} B=A$,  the operator $B$ has the matrix representation
\begin{equation}\label{B=matrix}
B=\left[ 
\begin{matrix}
A\\\Lambda D_{A,t}\\
XD_{A,t}
\end{matrix}\right]: \cH\to \cH'\oplus [F^2(H_n)\otimes \cD']\oplus \cE',
\end{equation}
where  the operator
$$
\left[
\begin{matrix}
 \Lambda  \\
X 
\end{matrix}\right]:  \cD_{A,t}\to  [F^2(H_n)\otimes \cD']\oplus \cE'
$$
is a contraction.
According to Theorem \ref{maxprin2}, the equation $\Delta(B) = \Delta(A)$ implies 
$
\left[
\begin{matrix}
 \Lambda  \\
X 
\end{matrix}
\right]
 |\cM=0$, i.e., $\Lambda|\cM=0$ and
 $X|\cM=0$.
According to  Remark 
\ref{restr=0}, we deduce that $B_c=\left[ 
\begin{matrix}
A\\\Lambda D_{A,t} 
\end{matrix}\right]
$
is the central intertwining lifting of $A$  with respect to $\cV'$ and tolerance $t$.
Since 
$
\left[
\begin{matrix}
 \Lambda  \\
X 
\end{matrix}\right] 
$
is a contraction,
 there exists another contraction $G:\cD_\Lambda\to \cE'$ such that $X=GD_\Lambda$.
On the other hand, since 
$\Lambda|\cM=0$, we have $D_\Lambda|\cM=I_\cM$. From $X|\cM=0$, we deduce $G|\cM=0$.
The matrix representation \eqref{B=matrix} becomes 
\begin{equation}\label{B=matrix2}
B=\left[ 
\begin{matrix}
A\\\Lambda D_{A,t}\\
GD_\Lambda D_{A,t}
\end{matrix}\right], 
\end{equation}
where $\Lambda$ and $G$ are contraction  such that $\Lambda|\cM=0$ and $G|\cM=0$.
Since $V_i' B_c= B_c T_i$  ~for any $i=1,\ldots, n$, we have 
\begin{equation*}
\begin{split}
\left\|D_{B_c,t}\left( \sum_{i=1}^n T_i h_i\right) \right\|^2
&=
\left< \sum_{i,j=1}^n T_i^* D_{B_c,t}^2 T_i h_i, h_j\right>\\
&=
\sum_{i,j=1}^n \left< \delta_{ij}   D_{B_c,t}^2  h_i, h_j\right>= 
\sum_{i=1}^n \|D_{B_c,t}^2  h_i\|^2.
\end{split}
\end{equation*}
Hence, we deduce that there are some unique isometries $V_i\in B(\cD_{B_c,t})$, 
\ $i=1,\ldots, n$, 
such that 
\begin{equation}
\label{VDB}
V_i D_{B_c, t}= D_{B_c, t} V_i, \quad, i=1,\ldots, n,
\end{equation}
and 
$$
\left\|\sum_{i=1}^n V_i D_{B_c, t}h_i\right\|^2=  \sum_{i=1}^n \left\|D_{B_c, t}h_i\right\|^2.
$$
This proves that $V_1,\ldots, V_n$ are isometries with orthogonal ranges.
Since $[E_1^* ~\cdots ~ E_n^*]$ is a $C_0$-row contraction, we have $Q=0$ in Lemma
\ref{lem2}. Therefore, Theorem \ref{quasi} implies
\begin{equation}
\label{and}
t^2I-B_c^* B_c=Z^* Z\quad \text{ and } \quad ZT_i= (S_i\otimes I_\cM) Z
\end{equation}
for any $i=1,\ldots, n$, where $Z: \cH\to F^2(H_n)\otimes \cM$ has dense range and is defined as in Theorem
\ref{quasi}.
According to the factorization \eqref{and}, there exists an isometric operator
$U: \cD_{B_c, t}\to F^2(H_n)\otimes \cM$ such that
\begin{equation}\label{UDZ}
UD_{B_c,t}= Z.
\end{equation}
Since $Z$ has dense range, we infer that $U$ is a unitary operator.
Using relations \eqref{VDB}, \eqref{and}, and \eqref{UDZ}, we obtain
\begin{equation*}
\begin{split}
UV_iD_{B_c, t}&= UD_{B_c, t}T_i= ZT_i\\
&=
(S_i\otimes I_\cM) Z= (S_i\otimes I_\cM)UD_{B_c, t}.
\end{split}
\end{equation*}
Therefore, 
$$
UV_i=(S_i\otimes I_\cM)U, \quad i=1,\ldots, n.
$$
On the other hand, since 
$B_c=\left[ 
\begin{matrix}
A\\\Lambda D_{A,t} 
\end{matrix}\right]
$ 
 and $\Lambda$ is a contraction, there is a unitary operator
 $U_\Lambda: \cD_{B_c,t}\to \cD_\Lambda$ such that
 \begin{equation}
 \label{mo}
 U_\Lambda D_{B_c,t}= D_\Lambda D_{A,t}.
 \end{equation}
Now, define the isometries $W_i\in B(\cD_\Lambda)$ by setting
\begin{equation}
\label{Wii}
W_i:= U_\Lambda V_i U_\Lambda^*, \quad i=1,\ldots, n.
\end{equation}
Note that relations \eqref{mo}, \eqref{Wii}, and \eqref{VDB} imply

\begin{equation}
 \label{mo2}
 W_i D_\Lambda D_{A,t}= D_\Lambda D_{A,t}T_i, \quad i=1,\ldots, n.
 \end{equation}
Indeed, we have
$$
U_\Lambda V_i U_\Lambda^*D_\Lambda D_{A,t}= U_\Lambda V_i D_{B_c, t} =
U_\Lambda D_{B_c, t}T_i =D_\Lambda D_{A,t}T_i
$$
for any $i=1,\ldots, n$.
Now, we prove that the subspace 
$
\cM:=\cD_{A, t}\ominus \bigvee_{i=1}^n D_{A,t} T_i\cH
$
is equal to $\bigcap_{i=1}^n \ker W_i^*$.
According to the Wold type decomposition for sequences of isometries with orthogonal ranges
(see \cite{Po-isometric}) and using \eqref{Wii},  we have
\begin{equation*}
\begin{split}
\bigcap_{i=1}^n \ker W_i^*&=
\cD_\Lambda \ominus [W_1\cD_\Lambda\oplus\cdots \oplus W_1\cD_\Lambda]\\
&=
\cD_\Lambda \ominus [\oplus_{i=1}^n W_i D_\Lambda D_{A,t} \cH]\\
&=
\cD_\Lambda \ominus [\oplus_{i=1}^n  D_\Lambda D_{A,t} T_i \cH]\\
&= 
\cD_\Lambda \ominus D_\Lambda\cN.
\end{split}
\end{equation*}
Since 
$\Lambda|\cM=0$,  we have $D_\Lambda|\cM=I_\cM$ and $\cM$ is a reducing  subspace 
for $D_\Lambda$.
Hence, we deduce that
$$
\cD_\Lambda =\overline{ D_\Lambda\cM}\oplus \overline{ D_\Lambda\cN}=
\cM\oplus \overline{ D_\Lambda\cN}.
$$
The previous computations show that $\bigcap_{i=1}^n \ker W_i^*= \cM$.
Now, since $W_i'B= BT_i$ and $W_i'= V_i'\oplus U_i'$ \ for any $i=1,\ldots, n$,
 we can take into account  the matrix representation
 \eqref{B=matrix2} and relation \eqref{mo2} to deduce 
 $$
 U_i' G D_\Lambda D_{A,t} = G D_\Lambda D_{A, t} T_i= GW_i D_\Lambda D_{A, t}
 $$
 for any $i=1,\ldots, n$.
 Hence, we get  $U_i' G= GW_i$, \ $i=1,\ldots, n$. Using the fact that 
 $G|\cM=0$,  we obtain
 $$
 GW_\alpha h= U_\alpha' h=0
 $$
 for any $\alpha\in \FF_n^+$ and $h\in \cM$. Since
  $\bigvee_{\alpha\in \FF_n^+} W_\alpha \cM= \cD_\Lambda$
  and $G:\cD_\Lambda \to \cE'$, we  get $G=0$.
  Therefore,
  $B=\left[ 
\begin{matrix}
B_c\\ 0
\end{matrix}\right]= \tilde{B_c}
$ 
is the central intertwining lifting of $A$.
To prove that last part of the theorem, note that if $\|A\|<t$ and $[T_1~\cdots ~ T_n]$
is an orthogonal shift, then according to Lemma \ref{lem2}, the operator $[E_1^*~\cdots ~ E_n^*]$
is a $C_0$-row contraction. Therefore, we can apply the first part
 of the theorem to complete the proof. 
\end{proof}

 The following result is  now a consequence of the previous theorem.
  Since the proof is similar to that of Theorem
 \ref{inter-entro}, we shall omit it.

\begin{proposition}\label{en-gen}
Let $\cT':=[T_1'~\cdots~ T_n']$, \ $T_i'\in B(\cH')$,  be a row contraction,
$\cW':=[W_1'~\cdots~W_n']$ be an isometric lifting of $\cT'$, 
and let 
 $A\in B(F^2(H_n)\otimes \cE, \cH')$   be  such that  $\|A\|\leq t$  and 
 $$
 A(S_i\otimes I_\cE)=T_i'A, \quad i=1,\ldots, n.
 $$
 If 
   $B$ is an intertwining lifting of $A$ with respect to $\cW'$ and tolerance $t>0$,
   then 
   \begin{equation}
   \label{e-g}
   E(A)=E(\tilde{B_c}) \geq E(B),
   \end{equation}
where  $\tilde{B_c}$ is the central intertwining lifting of $A$ with tolerance $t>0$.
Assume that $[E_1^* ~\cdots ~ E_n^*]$ is a $C_0$-row contraction,
 where  $E_i:= P_i W P_\cN\in B(\cD_{A,t})$  for any $i=1,\ldots, n$.
   If  the entropy ~$E(\tilde{B_c})>-\infty$, then  
     $E(\tilde{B_c})= E(B)$  if and only if $B=\tilde{B_c}$.
   
   In particular, if $\|A\|<t$, then  $E(\tilde{B_c})= E(B)$ 
    if and only if $B=\tilde{B_c}$.
\end{proposition}

 We remark that if $\|A\|<t$, then the operator 
 $\Delta(A)\in B(\cE)$ is invertible and $E(A)>-\infty$.

\clearpage
\section{Maximal entropy interpolation problems in several variables}

  We obtain explicit forms for the maximal entropy
 solution (as well as its entropy) of the Sarason \cite{S}, 
 Carath\' eodory-Schur \cite{Ca},
 \cite{Sc}, 
 and Nevanlinna-Pick  \cite{N} type interpolation problems for the 
 noncommutative (resp.~commutative) analytic Toeplitz algebra $ F_n^\infty$
 (resp.~$W_n^\infty$) and their tensor products  with $B(\cH, \cK)$.
 Moreover, under certain conditions, we also find explicit forms for
  the corresponding classical optimization problems, 
  in our multivariable noncommutative (resp.~commutative)
  setting.
  In particular, we provide explicit forms for the maximal entropy
 solutions of  several    interpolation problems  
     on the unit  ball of $\CC^n$. Finnaly, we apply our permanence principle
     to the Nevanlinna-Pick interpolation problem on the unit ball.

\subsection{Maximal entropy solution for the  Sarason  interpolation problem 
 for   analytic  Toeplitz algebras}
\label{Sara}

Given $R\in R_n^\infty\bar\otimes B(\cE, \cE')$ with $\dim \cE<\infty$, and $\Phi\in
 R_n^\infty\bar\otimes B(\cE_1, \cE')$ an inner operator, we consider the following 
 Sarason  (see \cite{S}) type 
 left interpolation problem  with tolerance $t>0$
 for $R_n^\infty\bar\otimes B(\cE, \cE')$:
 
 {\it Find $\Psi \in R_n^\infty\bar\otimes B(\cE, \cE')$  with maximal entropy such that 
 \begin{equation} \label{sa1}
 \|\Psi\|\leq t \ \text{ and } \ \Psi= R+ \Phi G,
 \end{equation}
where $G\in R_n^\infty\bar\otimes B(\cE, \cE_1)$.}

In this section,  we find the  solution to this problem 
and prove a permanence principle for it
 (see Theorem \ref{sara1} and Theorem \ref{sara-p}). 
  Under  a certain   condition, we 
obtain explicit forms for the maximal entropy
solution of the Sarason interpolation problem for analytic Toeplitz algebras in our multivariable 
setting
(see Theorem \ref{sara} and Corollary \ref{sara-com}).
We  also find an explicit form for the unique solution of the
corresponding Sarason  optimization problem
(see Theorem \ref{opti} and Theorem \ref{opti2}) in our setting.

We recall  (from Section \ref{factorizations}) that the entropy of  a multi-analytic  operator
 $\Psi\in R_n^\infty\bar\otimes B(\cE, \cE')$  with $\|\Psi\|\leq t$  is defined by
$E(\Psi):= \ln \det \Delta(\Psi)$, where 
$$
\left< \Delta(\Psi)h,h\right>:= \inf \left< (t^2I-\Psi^* \Psi)(h-p), h-p\right>,
\quad h\in \cE,
$$
where the infimum is taken over all polynomials $p\in F^2(H_n)\otimes \cE$ with $p(0)=0$.

\begin{theorem}\label{sara1} Let  $R\in R_n^\infty\bar\otimes B(\cE, \cE')$ 
with $\dim \cE<\infty$, and $\Phi\in
 R_n^\infty\bar\otimes B(\cE_1, \cE')$ be  an inner operator.
 Let $A:F^2(H_n)\otimes \cE\to \cH'$ be defined by $A:=P_{\cH'} R$, where
 \begin{equation}
 \label{h'}
 \cH':=[F^2(H_n)\otimes \cE']\ominus \Phi[F^2(H_n)\otimes \cE_1].
 \end{equation}
 Then there is a   solution for the interpolation  problem \eqref{sa1}
  if and only if 
 $\|A\|\leq t$.
 Moreover, the central intertwining lifting $\Psi_{\text{\rm max}}$ of $A$ with 
 respect to
  $\{S_i\otimes I_\cE\}_{i=1}^n$
 and  $\{S_i\otimes I_{\cE'}\}_{i=1}^n$ is the maximal entropy solution for
 the problem \eqref{sa1}.
\end{theorem}

\begin{proof}
First, note that if $\Psi $ is a solution to the interpolation problem 
\eqref{sa1}, then  we have $P_{\cH'} \Psi =P_{\cH'}R= A$ and $\|A\|\leq t$.
Conversely, let $T_i\in B(\cH')$, \  $i=1,\ldots, n$,  be defined  by $T_i':= P_{\cH'}(S_i\otimes I_{\cE'})|\cH'$ and note that
$$
A(S_i\otimes I_\cE)= T_i' A, \quad i=1,\ldots, n.
$$
Since $[S_1\otimes I_{\cE'}~\cdots ~ S_n\otimes I_{\cE'}]$ is an isometric dilation of 
$[T_1'~\cdots ~ T_n']$, we can apply Theorem \ref{main3} and the remarks following Theorem
\ref{restr=0}   to this setting. Therefore, the central
 intertwining lifting $\tilde{B_c}$ of $A$  with tolerance $t$  
 satisfies the conditions $P_{\cH'} \tilde{B_c}=A$,
  $\|\tilde{B_c}\|\leq t$, and
 $$
 \tilde{B_c} (S_i\otimes I_\cE)= (S_i\otimes I_{\cE'})\tilde{B_c}, \quad i=1,\ldots, n.
 $$
 Hence, $\tilde{B_c}= \Psi_{\text{\rm max}}$ for some  multi-anlytic  operator
 $\Psi_{\text{\rm max}}\in R_n^\infty\bar\otimes B(\cE, \cE')$. Since 
 $P_{\cH'}(\Psi_{\text{\rm max}}- R)=0$, according to \cite{Po-central}, we have
  $\Psi_{\text{\rm max}}- R\in \Phi(R_n^\infty\bar\otimes B(\cE, \cE'))$.
  Hence, $\Psi_{\text{\rm max}}= R+\Phi G$ for some $G\in R_n^\infty\bar\otimes B(\cE, \cE_1)$.
  On the other hand, 
  since $\Psi_{\text{\rm max}}$ is the central intertwining lifting of $A$, 
  Theorem \ref{maxprin3} and Proposition \ref{en-gen} show that 
  it is
   also the maximal entropy
  solution to the interpolation problem \eqref{sa1}. The proof is complete.  
\end{proof}

Let $\cK'\subseteq \cE'$ be an invariant subspace under  each operator
 $S_i^*\otimes I_{\cE'}$,
\ $i=1,\ldots, n$, such that
$\cH'\subseteq \cK'$, where $\cH'$ is given by relation \eqref{h'}.
 According to Theorem 2.2 from \cite{Po-charact}, there exists an inner 
operator $\Psi\in R_n^\infty\bar\otimes B(\cE_2, \cE')$ such that 
\begin{equation}\label{Beur}
\cK'= [F^2(H_n)\otimes \cE']\ominus \Psi[F^2(H_n)\otimes \cE_2].
\end{equation}
Since $\cK'\supseteq \cH'$, we have  
$$
\Psi[F^2(H_n)\otimes \cE_2]\subset
 \Phi[F^2(H_n)\otimes \cE_1].
 $$
  Hence,  as in the proof of Theorem 3.7 
 from \cite{Po-central}, one can show  that there is 
 $\Phi_1\in R_n^\infty\bar\otimes B(\cE_2, \cE_1)$ such that $\Psi =\Phi \Phi_1$. Conversely,
 if $\Psi =\Phi \Phi_1$, then the subspace $\cK'$ defined by \eqref{Beur} defines 
 an invariant subspace under  each operator $S_i^*\otimes I_{\cE'}$,
\ $i=1,\ldots, n$.

 Let  
 $\Phi\in
 R_n^\infty\bar\otimes B(\cE_1, \cE')$ 
 and $\Phi_1\in R_n^\infty\bar\otimes B(\cE_2, \cE_1)$ be  inner operators,
 and let $\Psi_{\text{\rm max}} \in R_n^\infty\bar\otimes B(\cE, \cE')$ be
 the maximal entropy solution of the problem  \eqref{sa1}.
  Consider the following 
   interpolation problem:

 {\it Find $\Gamma \in R_n^\infty\bar\otimes B(\cE, \cE')$  with maximal entropy such that 
 \begin{equation} \label{sa2}
 \|\Gamma\|\leq t \ \text{ and } \ \Gamma= \Psi_{\text{\rm max}}+  \Phi \Phi_1 G,
 \end{equation}
where $G\in R_n^\infty\bar\otimes B(\cE, \cE_2)$}.

We can prove the following permanence 
principle for the interpolation problem \eqref{sa1}.

\begin{theorem}\label{sara-p} 
Let  $R\in R_n^\infty\bar\otimes B(\cE, \cE')$ 
with $\dim \cE<\infty$, and $\Phi\in
 R_n^\infty\bar\otimes B(\cE_1, \cE')$ be
   a pure inner operator.
 Let $\Psi_{\text{\rm max}}$ be the maximal entropy solution for the interpolation 
 problem \eqref{sa1}.
 If $\Phi_1\in R_n^\infty\bar\otimes B(\cE_2, \cE_1)$ 
  is an  inner operator, then 
 $\Psi_{\text{\rm max}}$  is also  the maximal entropy solution for the interpolation 
 problem  \eqref{sa2}.
 \end{theorem}

\begin{proof}
Let $\Psi_{\text{\rm max}}$ be the central intertwining lifting of the operator 
$A:=P_{\cH'} R$. 
Then according to   Theorem \ref{sara1}, $\Psi_{\text{\rm max}}$ is a
 maximal entropy solution
 for the interpolation problem \eqref{sa1}.
 Note that since $\Phi$ is pure, i.e., $\|P_{\cE'} \Phi h\|< \|h\|$, \ $h\in \cE_1$,
  the results from
 \cite{Po-charact} imply that 
$[S_1\otimes I_{\cE'}~\cdots ~ S_n\otimes I_{\cE'}]$ is  the minimal  isometric dilation of 
the row contraction
$[T_1'~\cdots ~ T_n']$.
Set 
$$
\cK'= [F^2(H_n)\otimes \cE']\ominus \Phi \Phi_1[F^2(H_n)\otimes \cE_2] 
$$
and let $Y_{\cK'}:= P_{\cK'} \Psi_{\text{\rm max}}$ be the corresponding
 central partial intertwining lifting of $A$. The permanence principle of Theorem 
 \ref{perman}
 implies that $\Psi_{\text{\rm max}}$ is also the central intertwining lifting of $Y_{\cK'}$.
 Applying Theorem \ref{sara1} to the operator  $Y_{\cK'}$ instead of $A$, we conclude that 
  $\Psi_{\text{\rm max}}$ is also  
  the maximal entropy solution of the interpolation 
 problem  \eqref{sa2}.
\end{proof}

We remark that the previous theorem will be used in Section \ref{ball} 
to obtain  a permanence  principle for the Nevanlinna-Pick  interpolation 
problem  on the unit ball of $\CC^n$.

Now, using the results of the previous sections, we  can obtain explicit forms
 for the maximal entropy
 solution of the Sarason type interpolation  problem 
 for $R_n^\infty\bar{\otimes} B(\cE, \cY)$.

\smallskip

\begin{theorem}\label{sara}
Let $\cH'\subset F^2(H_n)\otimes \cY$ be an invariant subspace under each 
operator
$S_i^*\otimes I_\cY$, \ $i=1,\ldots, n$,  and
 let $A: F^2(H_n)\otimes \cE\to \cH'$ 
 be   such that $\|A\|<t$ and 
 $$
 A(S_i\otimes I_\cE)= [P_{\cH'}(S_i\otimes I_\cY)|\cH'] A, \quad i=1,\ldots, n.
 $$
 Let $ \psi:\cE\to F^2(H_n)\otimes \cE$ and $ \theta:\cE\to F^2(H_n)\otimes \cY$
  be operators defined by
 \begin{equation}\label{psi-theta}
 \psi h:= (t^2I-A^*A)^{-1}(1\otimes h)\quad \text{ and } 
 \quad \theta h:=A(t^2I-A^*A)^{-1}(1\otimes h),
 \end{equation}
 for any $h\in \cE$.
 Then $M_\psi, M_\theta$ are multi-analytic operators, with $M_\psi$ 
 invertible and 
 $$
 M_\psi^{-1} (1\otimes h)= \sum_{\sigma\in \FF_n^+} e_\sigma \otimes
  \Delta(A) P_\cE (X_\sigma^A)^*h, \qquad h\in \cE.
 $$ 
 The central intertwining lifting $\tilde{B_c}$ of $A$ with respect 
 to $\{S_i\otimes I_\cE\}_{i=1}^n$ and $\{S_i\otimes I_\cY\}_{i=1}^n$
 is equal to 
 $
   M_\theta M_\psi^{-1}
 $
 and the multi-Toeplitz operator
 $t^2I-\tilde{B_c}^* \tilde{B_c}$ admits 
  a square outer spectral factorization
  $$
  t^2I-\tilde{B_c}^* \tilde{B_c}= M_\chi^*M_\chi,
  $$
  where $ M_\chi\in R_n^\infty\bar{\otimes} B(\cE)$ is an outer operator
  given by
  $$
   M_\chi:= [I\otimes (P_\cE \psi)^{1/2}] M_\psi^{-1}.
  $$
   Moreover, if $\dim \cE<\infty$, then $ M_\theta M_\psi^{-1}$ is the maximal entropy 
   lifting of $A$ and 
   $$
   E( M_\theta M_\psi^{-1})=-\ln \det [P_\cE(t^2I-A^* A)^{-1}|\cE].
   $$
\end{theorem}
\begin{proof}
Note that $[S_1\otimes I_\cY~\cdots ~ S_n\otimes I_\cY]$ is an isometric lifting of 
the row-contraction
$$
\cC:=[P_{\cH'}(S_1\otimes I_\cY)|\cH' ~\cdots ~ P_{\cH'}(S_n\otimes I_\cY)|\cH'].
$$
According to the remarks following Theorem  \ref{restr=0}, there is a unique isometry
$$
\Phi: \cK':=\cH'\oplus[F^2(H_n)\otimes \cD']\to F^2(H_n)\otimes \cY,
$$
 such that 
$(S_i\otimes I_\cY)\Phi= \Phi V_i'$,  for any $i=1,\ldots, n$, and $\Phi|\cH'= \cH'$, where 
$[V_1'~\cdots V_n']$ is the  minimal isometric dilation of $\cC$ on 
$\cK'$.
 The operator ${\tilde B_c}:= \Phi B_c$ is the central intertwining
 dilation of $A$ with respect to
 $[S_1\otimes I_\cY~\cdots ~ S_n\otimes I_\cY]$.
Since 
$$
(S_i\otimes I_\cY){\tilde B_c}= {\tilde B_c} (S_i\otimes I_\cE), \ i=1,\ldots, n,
$$
we infer that ${\tilde B_c}= M_f$ for some multi-analytic operator
$M_f\in R^\infty_n\bar\otimes B(\cE, \cY)$. On the other hand we have 
$\|M_f\|= \|B_c\|\leq t$.
According to Theorem \ref{strict}, $M_\psi $ is an invertible multi-analytic operator on 
$F^2(H_n)\otimes \cE$.
Using Lemma \ref{lem1}, we get
\begin{equation*}
\begin{split}
{\tilde B_c} \psi h&= \Phi B_c D_{A,t}^{-2} (1\otimes h)=\Phi AD_{A,t}^{-2}(1\otimes h) \\
&= A D_{A,t}^{-2} (1\otimes h)= \theta h
\end{split}
\end{equation*}
for any $h\in \cE$.
Therefore ${\tilde B_c}M_\psi= M_\theta$, and $ M_\theta\in R^\infty_n\bar\otimes B(\cE)$.
Now,  using Corollary \ref{strict-form} and Theorem \ref{strict} we can complete 
the proof of the theorem.
\end{proof}

 Let $F_s^2(H_n)\subset F^2(H_n)$ be the symmetric Fock space and let us  define
 $W_n^\infty:= P_{F_s^2(H_n)} F_n^\infty|F_s^2(H_n)$.
   We recall that $W_n^\infty$ is the WOT-closed algebra generated by
  $$
  B_i:=P_{F_s^2(H_n)} S_i|F_s^2(H_n), \quad i=1,\ldots, n,
  $$
   and the identity.
  Let $\cH_s\subseteq F_s^2(H_n)\otimes \cY$ be an invariant subspace under each 
  operator $B_i^*\otimes I_\cY$, 
  \ $i=1,\ldots, n$, and let $C:F_s^2(H_n)\otimes \cE\to \cH_s$ be an operator  satisfying
  $\|C\|\leq t$ and 
  \begin{equation}\label{CBi}
  C(B_i\otimes I_\cE)= [P_{\cH_s} (B_i\otimes I_\cY)|\cH_s]C, \quad i=1,\ldots, n.
  \end{equation}

Since $F_s^2(H_n)$ is an invariant   subspace under each $S_i^*$, \ $i=1,\ldots, n$,
 it is easy to see
that $\cH_s$ is also invariant under $S_i^*\otimes I_\cY$, \ $i=1,\ldots, n$.
Setting 
$$
A:=C(P_{F_s^2(H_n)} \otimes I_\cE) : F^2(H_n)\otimes \cE\to F^2(H_n)\otimes \cY,
$$
relation \eqref{CBi} implies
$$
A(S_i\otimes I_\cE)= [P_{\cH_s} (S_i\otimes I_\cY)|\cH_s]A, \quad i=1,\ldots, n.
$$
Let $\tilde B_c$ be the central intertwining lifting of $A$ with respect to 
$\{S_i\otimes I_\cE\}_{i=1}^n $ and $\{S_i\otimes I_\cY\}_{i=1}^n $, i.e.,
$$
\tilde B_c\in R_n^\infty\bar\otimes B(\cE, \cY), \quad P_{\cH_s} \tilde B_c=A,  \quad 
\text{and } \ \|\tilde B_c\|\leq t.
$$
Define 
$$
\tilde C_c:= (P_{F_s^2(H_n)} \otimes I_\cY)\tilde B_c|F_s^2(H_n)\otimes \cE\in
 W_n^\infty\bar\otimes B(\cE, \cY)
 $$
and note that 
$P_{\cH_s} \tilde C_c=C$  and $\|\tilde C_c\|\leq t$.
We call $\tilde C_c$ the central intertwining lifting of $C$
with respect to 
$\{B_i\otimes I_\cE\}_{i=1}^n $ and $\{B_i\otimes I_\cY\}_{i=1}^n $, and tolerance $t$. 
Note that  $\tilde C_c$ is a solution of the Sarason interpolation problem  
for the operator space $W_n^\infty\bar\otimes B(\cE, \cY)$.

 If $\|C\|<t$, then one  can apply 
Theorem \ref{sara} to the  bounded operator $A:=C(P_{F_s^2(H_n)} \otimes I_\cE)$
and  deduce the following   multivariable  commutative 
version for  the Sarason   interpolation problem.

\begin{corollary}\label{sara-com}
Let $\cH_s\subseteq F_s^2(H_n)\otimes \cY$ be an invariant subspace under each 
  operator $B_i^*\otimes I_\cY$, 
  \ $i=1,\ldots, n$, and let $C:F_s^2(H_n)\otimes \cE\to \cH_s$ be  
   such that $\|C\|<t$ and 
  \begin{equation}\label{CBii}
  C(B_i\otimes I_\cE)= [P_{\cH_s}  (B_i\otimes I_\cY)|\cH_s ]C, \quad i=1,\ldots, n.
  \end{equation}
 Let  $\tilde C_c\in W_n^\infty\bar\otimes B(\cE, \cY)$ be defined by
$$
\tilde C_c:=(P_{F_s^2(H_n)} \otimes I_\cY)  M_\theta M_\psi^{-1}|F_s^2(H_n)\otimes \cE,
$$
where $\theta, \psi$ are defined by relation \eqref{psi-theta} and 
$A:=C(P_{F_s^2(H_n)} \otimes I_\cE)$.
Then $ \tilde C_c$ is an intertwining lifting of $C$, i.e.,
$$
P_{\cH_s}\tilde C_c= C \quad \text{ and } \quad \|\tilde C_c\|\leq t.
$$
\end{corollary}

In what follows, we show that, under certain conditions on the 
operator $A$,  there is   an explicit  form of the unique intertwining
 lifting $B$ of $A$ such that $\|A\|=\|B\|$.
We recall that  an operator $T\in B(\cX,\cY)$ attains its norm if there is
 a vector $x\in \cX$ of norm one such that $\|Tx\|=\|T\|$.

\begin{theorem}\label{opti}
Let $\cH'\subset F^2(H_n)\otimes \cY$ be an invariant subspace under each 
operator
$S_i^*\otimes I_\cY$, \ $i=1,\ldots, n$,  and
 let $A: F^2(H_n)\to \cH'$ 
 be a   contraction  with $\|A\|=1$ and  which  attains its norm.  If 
 $$
 AS_i= [P_{\cH'}(S_i\otimes I_\cY)|\cH'] A, \quad i=1,\ldots, n,
 $$
then there is a unique intertwining lifting 
 $M_\varphi \in R_n^\infty\bar\otimes B(\CC, \cY)$ of $A$
  such that $\|A\|=\|M_\varphi\|$.
  Moreover, $M_\varphi$ is an inner operator uniquely defined by the equation
  \begin{equation}\label{mpf}
  M_\varphi g=Ag
  \end{equation}
  and  $g$ is a   vector in $F^2(H_n)$ where $A$ attains its norm.
\end{theorem}
\begin{proof}
Let $g= \chi  f $ be the inner-outer factorization of $g$, where 
$\chi\in F_n^\infty$ is inner and $f\in F^2(H_n)$ is outer.
Since $AS_i =T_i'A$, \ $i=1,\ldots, n$, where 
$$
T_i':= P_{\cH'} (S_i\otimes I_\cY)|\cH',\quad i=1,\ldots, n,
$$
 and
$[T_1'~\cdots ~T_n']$ is a $C_0$-row contraction, 
we can use the $F_n^\infty$-functional
calculus for row contractions (see \cite{Po-funct}) and obtain  
\begin{equation*}\begin{split}
\|A\|&= \|Ag\|=\|A \chi(S_1,\ldots, S_n) f\|\\
&= \|\chi(T_1',\ldots, T_n')Af\|\leq \|Af\|\\
&\leq \|A\| \|f\|= \|A \|.
\end{split}
\end{equation*}
 Here,  we took into account that $\|f\|_2=1$ and $\|\chi\|=1$.
 Therefore, $\|A\|= \|Af\|$.
 By the noncommutative commutant lifting theorem,  
 there is an intertwining lifting $B$ of $A$ such that 
 $$
 BS_i= (S_i\otimes I_\cY)B, \quad i=1,\ldots, n,
 $$
 and $\|A\|=\|B\|$.
 According to \cite{Po-analytic}, there is $M_\varphi\in R_n^\infty\bar\otimes B(\CC, \cY)$
 such that $M_\varphi= B$.
 Since $P_\cH' M_\varphi =A$, we have
 \begin{equation*} \begin{split}
 \|A\|&=\|Af\|=\|P_\cH' M_\varphi f\|\\
 &\leq \|M_\varphi f\|\leq \|M_\varphi\|=\|A\|,
 \end{split}
 \end{equation*}
which implies $Af= P_\cH' M_\varphi f= M_\varphi f$.
Since $M_\varphi$ is  a multi-analytic operator, we deduce that
$$
M_\varphi\left(\sum_{|\alpha|\leq m} a_\alpha S_\alpha f\right)=
\sum_{|\alpha|\leq m} a_\alpha (S_\alpha\otimes I_\cY) Af, \quad a_\alpha\in \CC.
$$
Note that, since $A^*Af= f$, we have
\begin{equation*}\begin{split}
\left\|\sum_{|\alpha|\leq m} a_\alpha (S_\alpha\otimes I_\cY) Af\right\|^2 &= 
\sum_{\alpha>\beta, ~|\alpha|, |\beta|\leq m} 
\left<a_{\alpha\backslash \beta} (S_{\alpha\backslash \beta} \otimes I_\cY) Af, Af\right>\\
&{}\qquad + 
\sum_{\alpha<\beta, ~|\alpha|, |\beta|\leq m} 
\left< Af , a_{\beta\backslash \alpha} (S_{\beta\backslash \alpha} \otimes I_\cY) Af\right>\\
&=
\sum_{\alpha>\beta, ~|\alpha|, |\beta|\leq m} 
\left<a_{\alpha\backslash \beta} T'_{\alpha\backslash \beta}   Af, Af\right>\\
&{}\qquad + 
\sum_{\alpha<\beta, ~|\alpha|, |\beta|\leq m} 
\left< Af , a_{\beta\backslash \alpha} T'_{\beta\backslash \alpha}   Af\right>\\
&=
\sum_{\alpha>\beta, ~|\alpha|, |\beta|\leq m} 
\left<a_{\alpha\backslash \beta} A(S_{\alpha\backslash \beta} \otimes I_\cY) f, Af\right>\\
&{}\qquad + 
\sum_{\alpha<\beta, ~|\alpha|, |\beta|\leq m} 
\left< Af , a_{\beta\backslash \alpha} A(S_{\beta\backslash \alpha} \otimes I_\cY) f\right>\\
 =&
 \sum_{\alpha>\beta, ~|\alpha|, |\beta|\leq m} 
\left<a_{\alpha\backslash \beta} (S_{\alpha\backslash \beta} \otimes I_\cY) f, f\right>\\
&{}\qquad + 
\sum_{\alpha<\beta, ~|\alpha|, |\beta|\leq m} 
\left< f , a_{\beta\backslash \alpha} (S_{\beta\backslash \alpha} \otimes I_\cY) f\right>\\
&= \left\|\sum_{|\alpha|\leq m} a_\alpha S_\alpha f\right\|^2.
\end{split}
\end{equation*} 
Since $f$ is outer, we infer that $M_\varphi$ is an inner operator.

 Assume  now that there is another $M_{\varphi_1}\in R_n^\infty\bar\otimes B(\CC, \cY)$
such that $M_{\varphi_1} g = Ag$.  Then we have $M_\varphi g=M_{\varphi_1} g$.
Using the inner-outer factorization $g=\chi f$, we get $M_\varphi f=M_{\varphi_1} f$.
Since $f$ is outer, we have $M_\varphi =M_{\varphi_1} $. Therefore, $\varphi=\varphi_1$ 
and the proof is complete.
\end{proof}

According to the proof of Theorem \ref{opti}, one can assume that $\|Af\|=1$
  and $M_\varphi f= Af$ for some outer vector $f$  in $F^2(H_n)$. In this
   case, there is a sequence of polynomials $p_k\in F^2(H_n)$ such that 
   $$
   \|p_k(S_1,\ldots, S_n)f-1\|\to 0, \quad \text{ as } k\to\infty.
   $$
   Now, it is clear that 
   $$
   \varphi=\lim_{k\to \infty} p_k(S_1\otimes I_\cY,\ldots, S_n\otimes I_\cY)Af.
   $$
Note that if $\ker (I-A^* A)$ is one dimensional, then any vector 
where $A$ attains its norm 
is outer.

Our commutative version of Theorem \ref{opti} is the following.

\begin{theorem}\label{opti2}
Let $\cH_s\subset F_s^2(H_n)\otimes \cY$ be an invariant subspace under each 
operator
$B_i^*\otimes I_\cY$, \ $i=1,\ldots, n$,  and
 let $C: F_s^2(H_n)\to \cH_s$ 
 be a   contraction  which  attains its norm  $\|C\|=1$ and  
 \begin{equation}\label{cbi}
 CB_i= [P_{\cH_s}(B_i\otimes I_\cY)|\cH_s] C, \quad i=1,\ldots, n.
 \end{equation}
Then there  exists 
 $G \in W_n^\infty\bar\otimes B(\CC, \cY)$  
  such that  $P_{\cH_s } G=C$ and $\|G\|=\|C\|$.
  Moreover,  the operator  $G$ is given  by the equation
  \begin{equation}\label{Gfrac}
   G=\frac {Cg} {g},
  \end{equation}
  where  $g$ is any   vector in $F_s^2(H_n)$ where $C$ attains its norm.
\end{theorem}
\begin{proof}
Let 
$
A:=CP_{F_s^2(H_n)}  : F^2(H_n)\to F^2(H_n)\otimes \cY,
$
and note that 
$$
\|A\|= \|C\|=\|Cg\|=\|Ag\|=1.
$$
 Since  $F_s^2(H_n)$ is  an invariant subspace under each 
operator
$B_i^*$, \ $i=1,\ldots, n$, relation \eqref{cbi} implies
$$
AS_i= [P_{\cH_s} (S_i\otimes I_\cY)|\cH_s]A, \quad i=1,\ldots, n.
$$
According to Theorem \ref{opti}, there exists
$
M_\varphi\in R_n^\infty\bar\otimes B(\CC, \cY)
$ 
such that 
$$
 \ P_{\cH_s} M_\varphi= A, \ \|A\|=\|M_\varphi\|, \ \text{ and  } \ M_\varphi g= Ag.
$$
Hence, we infer that
$$
(P_{F_s^2(H_n) \otimes \cY}M_\varphi |{F_s^2(H_n) }) g =Cg.
$$
 Setting
 $G:= P_{F_s^2(H_n)\otimes \cY }M_\varphi |{F_s^2(H_n) }$ 
 and using the identification 
 of $W_n^\infty$  with the algebra of analytic multipliers of $F_s^2(H_n)$ 
 (see
 \cite{Arv1}), 
 it is clear that 
 $G\in W_n^\infty\bar \otimes B(\CC, \cY)$ has the required properties.
 The proof is complete.
 \end{proof}

We remark that one can obtain versions of Theorem \ref{opti} and Theorem \ref{opti2}
when $\|A\|=t$ and $t>0$.

\bigskip

\subsection{Maximal entropy solution for the  Carath\'eodory-Schur interpolation
 problem for analytic Toeplitz algebras }
\label{Cara}

 Let $q:= \sum\limits_{|\alpha|\leq m-1} R_\alpha \otimes A_\alpha$, 
 \ $A_\alpha\in B(\cE_1, \cE_2)$,   be a polynomial in 
 $R_n^\infty\bar\otimes B(\cE_1, \cE_2) $ and  let
 \begin{equation}\label{Car-opt}
 d_\infty:= \inf\left\{\|G\|:\ G\in R_n^\infty\bar\otimes B(\cE_1, \cE_2) 
 \text{ and } G_\alpha= A_\alpha,  \text{ if } |\alpha| \leq m-1\right\},
 \end{equation}
 where $\{G_\alpha\}_{\alpha\in \FF_n^+}$ are the Fourier coefficients of $G$.
 The Carath\'eodory-Schur (see \cite{Ca}, \cite{Sc}) optimization problem \eqref{Car-opt} is to find 
 $G\in R_n^\infty\bar\otimes B(\cE_1, \cE_2) $ with the smallest norm subject 
 to the constraint
 \begin{equation}\label{constr}
  G_\alpha= A_\alpha  \quad \text{   if  }  \alpha\in \FF_n^+, ~|\alpha| \leq m-1.
  \end{equation}
  The Carath\'eodory-Schur interpolation problem   with tolerance $t>d_\infty$  is to find 
  $G$ with $\|G\|\leq t$ subject 
 to the  the same constraint.
 In \cite{Po-analytic}, we proved that $d_\infty=\|A\|$, where 
 $A:= P_{\cP_{m-1}\otimes \cE_2} q$
 and $\cP_{m-1}$ denotes the set of all polynomials in $F^2(H_n)$ of degree $\leq m-1$.
 
 In what follows, we provide an explicit form of  the  maximal entropy solution for the 
  Carath\'eodory-Schur interpolation  problem   with tolerance $t>d_\infty$.
 For each $j=1,2$, 
 define
 the Fourier transform
 $\cF_j:\ell^2(\FF_n^+)\otimes \cE_j\to F^2(H_n)\otimes \cE_j$ 
 by setting
 \begin{equation}\label{f12}
 \cF_j[(y_\alpha)_{\alpha\in \FF_n^+}]=\sum_{\alpha\in \FF_n^+} e_\alpha\otimes y_\alpha
 \end{equation}
 for any $(y_\alpha)_{\alpha\in \FF_n^+}\in \ell^2(\FF_n^+)\otimes \cE_j$.

 \begin{theorem}\label{C-S}
 Let $q:= \sum\limits_{|\alpha|\leq m-1} R_\alpha \otimes A_\alpha$, 
 \ $A_\alpha\in B(\cE_1, \cE_2)$,   be a polynomial in 
 $R_n^\infty\bar\otimes B(\cE_1, \cE_2) $,  and let 
 $M:= [A_{\alpha, \beta}]_{\alpha|, |\beta|\leq m-1}$ be the operator matrix defined by
 $$
 A_{\alpha, \beta}:= 
 \begin{cases} A_{\alpha\backslash \beta}, & \text{ if  } \alpha \geq \beta\\ 
 0,  &  \text{  otherwise.  }  
\end{cases}
$$

 Let $t>d_\infty=\|M\|$ and define the operators $\psi:\cE_1\to F^2(H_n)\otimes \cE_1$
 and  $\theta:\cE_1\to F^2(H_n)\otimes \cE_2$ by
 \begin{equation}\label{pt}
 \psi:= \cF_1 (t^2I- M^* M)^{-1}\left[  \begin{matrix} 
 I_{\cE_1}\\0\\ \vdots\\0
 \end{matrix} \right] \quad \text{ and }\quad 
 \theta:= \cF_2 M(t^2I- M^* M)^{-1}\left[  \begin{matrix} 
 I_{\cE_1}\\0\\ \vdots\\0
 \end{matrix} \right],
  \end{equation}
 where $\cF_1, \cF_2$ are the Fourier transforms defined  by relation \eqref{f12}.
 Then  the operator
 $$
 M_\varphi:= M_\theta M_\psi^{-1}\in R_n^\infty\bar\otimes B(\cE_1, \cE_2) 
 $$
 is the central solution for the 
  Carath\'eodory-Schur interpolation problem   with tolerance $t$. 
  Moreover, if $\dim \cE_1<\infty$, then $M_\varphi$ is the maximal entropy solution.
 \end{theorem}
 
 \begin{proof}
 Let $\cH':= \cP_{m-1} \otimes\cE_2$ and $T_i'\in B(\cH')$ be defined by 
 $T_i':= P_{\cH'} (S_i\otimes I_{\cE_2}) |\cH'$, \ $i=1,\ldots, n$.
 Setting 
 $$
 A:= P_{\cH'} q: F^2(H_n)\otimes \cE_1\to \cH'\subset F^2(H_n)\otimes \cE_2,
 $$
 we can apply Theorem \ref{sara} to $A$ and find the central intertwining 
 lifting $\tilde B_c$ of  $A$ given by $\tilde B_c=M_\theta M_\psi^{-1}$ 
 in $R_n^\infty\bar\otimes B(\cE_1, \cE_2)$, where $\theta, \psi$ 
 are defined as in \eqref{psi-theta}.
 Since $A\cF_1= \cF_2[M~ 0]$, we have
 \begin{equation*}\begin{split}
 \theta h &= A(t^2I-A^* A)^{-1} (1\otimes h)\\
 &=A\cF_1 \cF_1^{-1}(t^2I-A^* A)^{-1} \cF_1 \cF_1^{-1}(1\otimes h)\\
 &=\cF_2 M(t^2I- M^* M)^{-1}\left[  \begin{matrix} 
 I_{\cE_1}\\0\\ \vdots\\0
 \end{matrix} \right]h 
 \end{split}
 \end{equation*}
 for any $h\in \cE_1$. Similarly, one can get the formula for $\psi$.
 Using again Theorem \ref{sara}, we can complete the proof.
 \end{proof}

  Now, we obtain an  explicit solution for the Carath\'eodory-Schur
   optimization problem \eqref{Car-opt}
  for the analytic Toeplitz algebra  $R_n^\infty$.
  Denote  $$
   N:=\text{\rm card}~\{ \alpha\in \FF_n^+:\ |\alpha|\leq m-1\}.
   $$

 \begin{theorem}\label{c-s}
 Let $q:= \sum_{|\alpha|\leq m-1}  a_\alpha R_\alpha$, 
 \ $a_\alpha\in \CC$,   be a polynomial in 
 $R_n^\infty\  $,  and let 
 $M:= [a_{\alpha, \beta}]_{\alpha|, |\beta|\leq m-1}$ be the  matrix defined by
 $$
 a_{\alpha, \beta}:= 
 \begin{cases} a_{\alpha\backslash \beta}, & \text{ if  } \alpha \geq \beta\\ 
 0,  &  \text{  otherwise.  }  
\end{cases}
$$
Let $y\in \CC^N$ be 
any vector which attains
the norm of $M$. Then there is a unique $\varphi \in R_n^\infty$ such that its first
 $N$ Fourier coefficients are $a_\alpha$, 
$|\alpha|\leq m-1$, and $\|\varphi\|=d_\infty=\|M\|$.
Moreover, $\frac{1} {d_\infty}\varphi$ is inner and $\varphi$ is uniquely determined by 
the equation
\begin{equation}\label{fF}
\varphi \cF y= \cF My,
\end{equation} 
where $\cF:\ell^2(\FF_n^+)\to  F^2(H_n)$ is the   Fourier transform defined by
$$
\cF[(a_\alpha)_{\alpha\in \FF_n^+}]=\sum_{\alpha\in \FF_n^+} a_\alpha e_\alpha 
  $$
 for any $(a_\alpha)_{\alpha\in \FF_n^+}\in \ell^2(\FF_n^+) $. 
 \end{theorem}
 
 \begin{proof}
 Let $\cH':= \cP_{m-1}$, $\cY:= \CC$, and  define the operator
$$
A:= P_{\cH'} q: F^2(H_n) \to \cH'\subset F^2(H_n).
$$
Applying Theorem \ref{opti} in this setting, we find $M_\varphi\in R_n^\infty$   such that 
its first $N$ Fourier coefficients are $a_\alpha$, $|\alpha|\leq m-1,$ and
 $\|M_\varphi\|=d_\infty= \|M\|$. Moreover $\frac{1} {d_\infty} M_\varphi$ is inner and $M_\varphi$ 
 is uniquely determined by the equation
 \begin{equation} \label{Mg}
 M_\varphi g=Ag,
 \end{equation}
 where $g$ is any vector in $F^2(H_n)$ where $A$ attains its norm.
 On the other hand, if $y\in \CC^N$  is 
a vector where $M$  attains
its  norm, then $\cF y$  is 
a vector where $A$  attains
its  norm. Therefore, relation \eqref{Mg} implies \eqref{fF}. The proof is complete.
 \end{proof}
 
 We mention that, using Corollary \ref{sara-com} and Theorem \ref{opti2},  one can 
 obtain 
  multivariable commutative versions 
   of  Theorem \ref{C-S} and Theorem \ref{c-s} for 
   $W_n^\infty\bar\otimes B(\cE_1,\cE_2)$
   and $W_n^\infty$, respectively. On the other hand,  Theorem \ref{c-s} can
    be extended to polynomials 
    $q\in R_n^\infty\bar\otimes B(\CC, \cE_2)$, where $\cE_2$ is an arbitrary 
    Hilbert space.


\bigskip
\subsection{Maximal entropy solution for  the Nevanlinna-Pick interpolation problem 
with operatorial argument
in several variables
}
\label{Nev}

In this section, we obtain   explicit forms for the maximal entropy solution 
of the left tangential Nevanlinna-Pick (see \cite{N}) interpolation problem 
with operatorial argument in several variables (see Theorem \ref{NP1}) as well
 as  for its entropy
(Theorem \ref{ ent-NP}).

As in  \cite{Po-nehari},  the spectral radius associated
 with a  sequence of operators 
$\cZ:=(Z_1,\ldots, Z_n)$,  $Z_i\in B(\cY)$,   is given by
$$
r(\cZ):=\lim_{k\to \infty}\|\sum_{|\alpha|=k} Z_\alpha 
Z_\alpha^*\|^{1/2k}
=\inf_{k\to \infty}\|\sum_{|\alpha|=k} Z_\alpha Z_\alpha^*\|^{1/2k}.
$$
Note that  if $Z_1Z_1^*+\cdots+ Z_n Z_n^*<rI_\cY$~ with $0<r<1$, 
then $r(\cZ)<1$.
Any element $f$ in $F_n^\infty\bar\otimes B(\cH, \cY)$ has a unique Fourier
 representation
$$
f\sim \sum_{\alpha\in \FF_n^+} S_\alpha\otimes A_{(\alpha)}
$$
for some operators $A_{(\alpha)}\in B(\cH, \cY)$ such that
$$
\sum_{\alpha\in\FF_n^+}A_{(\alpha)}^* A_{(\alpha)}\leq \|f\|^2 I.
$$
If 
$r(\cZ)<1$,  it makes sense to define {\it the evaluation} of $f$
at $(Z_1,\ldots, Z_n)$ by setting
\begin{equation}\label{evaluation}
f(Z_1,\ldots, Z_n):=\sum_{k=0}^\infty\sum_{|\alpha|=k} Z_{\tilde\alpha} A_{(\alpha)}.
\end{equation}
Now, using the fact that the spectral radius of $\cZ$ is strictly
 less than 1,  one can prove the norm convergence of the series  \eqref{evaluation}.

Given  $C\in B(\cH, \cY)$, we define 
the controllability operator $W_{\{\cZ,C\}}: F^2(H_n)\otimes \cH\to \cY$ 
 associated with
$\{\cZ,C\}$  by setting
$$
W_{\{\cZ,C\}} \left(\sum_{\alpha \in\FF_n^+} e_\alpha\otimes h_\alpha\right):= 
\sum_{k=0}^\infty \sum_{|\alpha|=k} Z_{\tilde{\alpha}} C h_\alpha,
$$
where $\tilde{\alpha}$ is the reverse of 
$\alpha:=g_{i_1}\cdots g_{i_k}\in \FF_n^+$, i.e., 
$\tilde{\alpha}:=g_{i_k}\cdots g_{i_1}$.
Since  $r(\cZ)<1$, note that $W_{\{\cZ,C\}}$ is  a well-defined  bounded 
operator. 
 We call the positive operator  $G_{\{\cZ,C\}}:=W_{\{\cZ,C\}}W_{\{\cZ,C\}}^*$ 
   the
 controllability grammian for 
$\{\cZ,C\}$.
 It is easy to see that 
\begin{equation}\label{gram}
G_{\{\cZ,C\}}=\sum_{k=0}^\infty\sum_{|\alpha|=k} Z_\alpha CC^*  Z^*_\alpha,
\end{equation}
 where the
series converges in norm. As in the classical case ($n=1$),
 we say that the pair
$\{\cZ,C\}$ is controllable if its grammian $G_{\{\cZ,C\}}$ is strictly positive.
We remark that $G_{\{\cZ,C\}}$ is the unique  positive solution of the Lyapunov equation
\begin{equation}\label{L}
X=\sum_{i=1}^n Z_iXZ_i^*+ CC^*.
\end{equation}

 Given  $\cZ:=[Z_1~\cdots ~Z_n]\in B(\oplus_{i=1}^n \cY, \cY)$ with $r(\cZ)<1$, and 
 the operators $B\in B(\cK, \cY)$ and $C\in B(\cH, \cY)$, we note that
\begin{equation}\begin{split}
W_{\{\cZ,B\}}(R_i\otimes I_\cK)&= Z_i W_{\{\cZ,B\}}\ \text{ and}\\
W_{\{\cZ,C\}}(R_i\otimes I_\cH)&= Z_i W_{\{\cZ,C\}} \ \text{ for any } i=1,\ldots, n .
\end{split}
\end{equation}
  On the other hand, 
 if $\Phi\in F_n^\infty\bar\otimes 
B(\cH,\cK)$ has the Fourier representation 
$$
\Phi\sim 
\sum_{\alpha\in \FF_n^+} S_\alpha\otimes A_{(\alpha)},
$$ then 
\begin{equation}\label{iyc}
[(I \otimes B)\Phi](\cZ)=C
\end{equation}
if and only if
 $$
 \sum_{k=0}^\infty\sum_{|\alpha|=k} Z_{\tilde{\alpha}} B A_{(\alpha)}=C.
 $$
A straightforward computation on  the elements of the form
$e_\beta\otimes h$, ~$h\in \cH$, \ $\beta\in \FF_n^+$,  shows that relation
 \eqref{iyc} holds if and only if 
\begin{equation}\label{wpw}
W_{\{\cZ,B\}}\Phi= W_{\{\cZ,C\}}.
\end{equation}
 
 Consider the optimization problem 
 \begin{equation}\label{opt-NP}
 d_\infty(\cZ, B, C):=\inf\left\{ \|\Theta\|:\ [(I\otimes B)\Theta](\cZ)=C, \
  \Theta\in F_n^\infty\bar\otimes B(\cH, \cK)\right\},
 \end{equation}
 which we call the standard left tangential Nevanlinna-Pick interpolation problem with
  multivariable operatorial argument.
 According to Theorem 7.2 from \cite{Po-nehari}, we have
   $d_\infty:= d_\infty(\cZ, B, C)<\infty$,
 if and only if there exists a bounded operator 
 $$
 A:F^2(H_n)\otimes \cH\to \cH':=\overline{ \text{\rm range }
  W_{\{\cZ,B\}}^*}\subseteq F^2(H_n)\otimes \cK
 $$
 such that 
 \begin{equation}\label{A}
 W_{\{\cZ,B\}}A= W_{\{\cZ,C\}}.
 \end{equation}
 Moreover, in this case $A$ is uniquely determined and there exists  an operator
 $\Theta_{\text{opt}}\in F_n^\infty\bar\otimes B(\cH, \cK)$ such that
 $[(I\otimes B)\Theta_{\text{opt}}](\cZ)=C$ and 
 \begin{equation}\label{opt}
 d_\infty=\|A\|=\|\Theta_{\text{opt}}\|.
 \end{equation}
  Setting $T_i':= P_{\cH'} (R_i\otimes I_\cK)|\cH'$, we have
  $T_i'A= A (R_i\otimes I_\cH)$ for any $i=1,\ldots, n$.
  It is easy to see that $\hat A$ is an intertwining lifting of $A$ with respect to 
  $\{ R_i\otimes I_\cK\}_{i=1}^n$ if and only if  
 $\hat A=\Theta$ for some 
$\Theta \in F_n^\infty\bar\otimes B(\cH, \cK)$ such that
 $W_{\{\cZ,B\}}\Theta= W_{\{\cZ,C\}}$. Using this fact, one can see
  that the optimization problem \eqref{opt-NP} is equivalent to  the following
  \begin{equation}\label{opt2-NP}
  d_\infty(\cZ, B, C):=\inf\left\{ \|\hat A\|:\ \hat A 
  \text{ is an intertwining lifting of } A\right\}.
  \end{equation}
  
 On the other hand, we say that $\Theta _c\in F_n^\infty\bar\otimes B(\cH, \cK)$
 is the central interpolant  for the standard Nevanlinna-Pick problem with tolerance
 $t>0$ if $\Theta_c$ is the central intertwining  lifting of $A$ with tolerance $t$, i.e., 
 $\|\Theta_c\|\leq t$.
  In  \cite{Po-nehari}, we proved that, given $t>0$, there exists 
 $\Theta \in F_n^\infty\bar\otimes B(\cH, \cK)$ satisfying 
 \begin{equation}
 \label{Iten}
 [(I\otimes B)\Theta ](\cZ)=C
 \ \text{ and } \ \|\Theta\|\leq t,
 \end{equation}
 if and only if 
 $$
 t^2 G_{\{\cZ,B\}}-G_{\{\cZ,C\}}\geq 0,
 $$
  where $G_{\{\cZ,B\}}$ and 
 $G_{\{\cZ,C\}}$ are the grammians for $\{\cZ,B\}$ and $\{\cZ,C\}$, respectively.
 Now,  using  the results from previous sections, we deduce that
  $\Theta_c$ is the maximal entropy  of the interpolation problem \eqref{Iten}, 
  if $\dim \cH<\infty$.

 In what follows,  assume that  the grammian
 $G_{\{\cZ,B\}}$ is strictly positive. According to \cite{Po-nehari}, there exists
  $\Theta_{\text{opt}}\in F_n^\infty\bar\otimes B(\cH, \cK)$ solving the problem 
  \eqref{opt-NP}, i.e., 
  \begin{equation}\label{-opt}
  [(I\otimes B)\Theta_{\text{opt}}](\cZ)=C  \ \text{ and } \ 
 d_\infty =\|\Theta_{\text{opt}}\|.
 \end{equation}
 Moreover, the spectral radius of $G_{\{\cZ,C\}}G_{\{\cZ,B\}}^{-1}$ is equal
  to $d_\infty^2$.

 Now, we can obtain an explicit form for the maximal entropy solution 
  of the standard left tangential Nevanlinna-Pick interpolation problem with
 operatorial argument in several variables.

 \begin{theorem}\label{NP1} Let 
  $\cZ:=[Z_1~\cdots ~Z_n]\in B(\oplus_{i=1}^n \cY, \cY)$, 
 $B\in B(\cK, \cY)$,  and $C\in B(\cH, \cY)$  be operators
  such  that  $r(\cZ)<1$ and  the grammian  $G_{\{\cZ,B\}}$ is strictly positive.
  Then $d_\infty<\infty$ and the central interpolant $\Theta_t$ with tolerance
   $t> d_\infty $ for the
  standard left tangential Nevanlinna-Pick problem  with operatorial argument 
  is given by 
  $$
  \Theta_t= \Phi \Psi^{-1},
  $$
  where $\Phi\in F_n^\infty \bar \otimes B(\cH, \cK)$ and 
  $\Psi\in F_n^\infty \bar \otimes B(\cH)$ are analytic Toeplitz operators defined by
  \begin{equation*}
  \begin{split}
  \Phi&:= \sum_{\alpha\in \FF_n^+} S_\alpha\otimes B^* (Z_{\tilde \alpha})^*
  (t^2 G_{\{\cZ,B\}}-G_{\{\cZ,C\}})^{-1} C \quad \text{ and }\\
  \Psi&:= 
  \frac {1} {t^2} \left[
  I+\sum_{\alpha\in \FF_n^+} S_\alpha\otimes C^* (Z_{\tilde \alpha})^*
  (t^2 G_{\{\cZ,B\}}-G_{\{\cZ,C\}})^{-1} C.
  \right]
  \end{split}
  \end{equation*}
 In particular, $\Theta_t\in F_n^\infty \bar \otimes B(\cH, \cK)$ is an analytic 
 Toeplitz operator such that 
 $$
 [(I\otimes B)\Theta_t](Z_1,\ldots, Z_n)=C  \ \text{ and } \ 
 \|\Theta_t\|<t.
 $$
 \end{theorem}
 
 \begin{proof} Since the grammian  $G_{\{\cZ,B\}}$ is strictly
 positive, the operator $W^*_{\{\cZ,B\}}$  has closed range and 
 $$
 P_{\text{\rm range}~ W^*_{\{\cZ,B\}}}= W^*_{\{\cZ,B\}}G_{\{\cZ,B\}}^{-1}W_{\{\cZ,B\}}.
 $$
 Note that the operator 
 \begin{equation}\label{WGW}
 A:= W^*_{\{\cZ,B\}}G_{\{\cZ,B\}}^{-1}W_{\{\cZ,C\}}
 \end{equation}
 satisfies the equation
$W_{\{\cZ,B\}}A= W_{\{\cZ,C\}}$ and 
$$
A(R_i\otimes I_\cH)= A T_i,\quad i=1,\ldots, n,
$$
where 
$$T_i':= P_{\cH'}(R_i\otimes I_\cK)|\cH'\quad \text{and } \ 
\cH':=\text{\rm range}~ W^*_{\{\cZ,B\}}.
$$   
Since $d_\infty=\|A\|$, we have  $\|A\|<t$. 
On the other hand, since 
$$
W_{\{\cZ,B\}} A A^* W_{\{\cZ,B\}}^*= W_{\{\cZ,C\}}W_{\{\cZ,C\}}^*,
$$
 it is easy to see that 
the operator
$t^2 G_{\{\cZ,B\}}-G_{\{\cZ,C\}}$ is strictly positive.
Let $B_c$ be the central intertwining lifting of $A$ with respect 
to $\{R_i\otimes I_\cK\}_{i=1}^n$ and tolerance $t$. Then 
$B_c\in F_n^\infty \bar \otimes B(\cH, \cK)$, $\|B_c\|\leq t$, and 
$W_{\{\cZ,B\}}B_c= W_{\{\cZ,C\}}$.

Now, we use  Theorem \ref{sara}  to find an explicit formula for $B_c$. 
To this end, according to  relation  \eqref{WGW}, straightforward calculations show that
$$D_{A,t}^2:= t^2I-A^* A= t^2I- W^*_{\{\cZ,C\}}G_{\{\cZ,B\}}^{-1}W_{\{\cZ,C\}} 
$$
and 
\begin{equation}\label{DA}
D_{A,t}^{-2}= \frac {1} {t^2} \left[
I+ W^*_{\{\cZ,C\}}\left( t^2 G_{\{\cZ,B\}}-G_{\{\cZ,C\}}\right)^{-1}W_{\{\cZ,C\}}
\right].
\end{equation}
 Hence, and using again  relation \eqref{WGW}, we obtain
 $$
 AD_{A,t}^{-2}= W^*_{\{\cZ,B\}}\left( t^2 G_{\{\cZ,B\}}-W_{\{\cZ,C\}}\right)^{-1}W_{\{\cZ,C\}}.
 $$
 Note that \ $W_{\{\cZ,C\}}(1\otimes h)=Ch$, \ $h\in \cH$, and 
 $$
 W_{\{\cZ,C\}}^*y= \sum_{\alpha\in \FF_n^+} e_\alpha \otimes C^* (Z_{\tilde \alpha})^* y,
  \quad y\in \cY.
 $$
 Therefore, the above equations  imply
 $$
 D_{A,t}^{-2}(1\otimes h)= \frac {1} {t^2} \left[ 1\otimes h+ 
 \sum_{\alpha\in \FF_n^+} e_\alpha \otimes C^* (Z_{\tilde \alpha})^*
 \left( t^2 G_{\{\cZ,B\}}-G_{\{\cZ,C\}}\right)^{-1} Ch\right]
 $$
 and 
 $$
 AD_{A,t}^{-2}(1\otimes h)= \sum_{\alpha\in \FF_n^+} e_\alpha \otimes B^* 
 (Z_{\tilde \alpha})^*
 \left( t^2 G_{\{\cZ,B\}}-G_{\{\cZ,C\}}\right)^{-1} Ch
 $$
 for any $h\in \cH$.
 
 Now   applying  Theorem \ref{sara}   to the operator $A$, which satisfies
 $\|A\|<t$, we find the explicit forms for $\Phi$ and $\Psi$ mentioned in the theorem.
 According to Corollary \ref{strict-form} and Theorem \ref{strict},
 the operator $\Psi$ is invertible in $F^\infty_n \bar\otimes B(\cH)$ and 
 $$
 t^2I-\Theta_t^* \Theta_t= \Psi^* \Psi.
 $$
 This implies $\|\Theta_t\|< t$, and the proof is complete.
 \end{proof}

 \begin{remark}\label{state-sol}
 Under the conditions of Theorem $\ref{NP1}$, the operators  $\Phi$ and $\Psi$  satisfy
    the following  state space formulas
  \begin{equation*}
  \begin{split}
  \Phi &= (I \otimes B^*)\left(I  
  -\sum_{i=1}^n S_i\otimes Z_i^* \right)^{-1}
   \left(I \otimes (t^2 G_{\{\cZ,B\}}-G_{\{\cZ,C\}})^{-1} C\right),\\
   \Psi&= 
   \frac {1} {t^2} \left[
  I+(I \otimes C^*)\left(I  
  -\sum_{i=1}^n S_i\otimes Z_i^* \right)^{-1}
   \left(I \otimes (t^2 G_{\{\cZ,B\}}-G_{\{\cZ,C\}})^{-1} C\right)\right].
 \end{split}
  \end{equation*}
 \end{remark}
 \begin{proof}
 Note that, for any $k=1,2,\ldots,$ the operators $S_\alpha$,  
 ~ $|\alpha|=k$, are isometries with
 orthogonal ranges.
 Now, it is easy to see that
 \begin{equation*}\begin{split}
 \|(S_1\otimes Z_1^*+\cdots + S_n \otimes Z_n^*)^k\|&=
 \left\|\sum_{|\alpha|=k} S_\alpha \otimes Z_{\tilde \alpha}^*\right\| \\
 &=
 \left\|\left(\sum_{|\alpha|=k} S_\alpha^* \otimes Z_{\tilde \alpha}\right)
  \left( \sum_{|\alpha|=k} S_\alpha \otimes Z_{\tilde \alpha}^* \right)\right\| \\
 &=\left\|\sum_{|\alpha|=k} I \otimes Z_{\tilde \alpha}Z_{\tilde \alpha}^*\right\|^{1/2}=
 \left\|\sum_{|\alpha|=k}  Z_{ \alpha}Z_{\alpha}^*\right\|^{1/2}.
 \end{split}
 \end{equation*}
 Since $\text{\rm r}(\cZ)<1$, the root test implies the convergence of the series
 $$
 \sum_{k=0}^\infty \|(S_1\otimes Z_1^*+\cdots + S_n \otimes Z_n^*)^k\|.
 $$
 Hence, it is clear that
 $$
 (I  
  -S_1\otimes Z_1^*-\ldots - S_n\otimes Z_n^*)^{-1}= 
  \sum_{k=0}^\infty\sum_{|\alpha|=k} S_\alpha \otimes Z_{\tilde \alpha}^*.
  $$
  Using Theorem \ref{NP1},  we can complete the proof.
 \end{proof}

 Let $\Theta\in F^\infty_n \bar\otimes B(\cH, \cK)$ be an interpolant for the 
 standard Nevanlinna-Pick problem with tolerance $t$, i.e., 
 \begin{equation} \label{NP<}
  [(I\otimes B)\Theta](\cZ)=C  \ \text{ and } \ 
 \|\Theta\|\leq t.
 \end{equation}
 Define 
 $$
 \left<\Delta(\Theta)x,x\right>=
 \inf\{\left< (t^2I-\Theta^* \Theta)(x-p), x-p\right>:\ 
 p\in F^2(H_n)\otimes \cH ~\text{ and } ~p(0)=0\} 
 $$
   for any  $x\in \cH$.
 If $\cH$ is finite dimensional, then the entropy of $\Theta$ is defined by
 $$
 E(\Theta):= \ln \det \Delta (\Theta).
 $$
 According to Theorem \ref{maxprin2}, we have 
  $\Delta(\Psi)\leq \Delta(\Theta_t)$,
 where $\Theta_t$ is the central interpolant for the problem \eqref{NP<} and $\Psi$
  is any interpolant for the same problem. Moreover, if $\dim \cH< \infty$, then 
  Corollary \ref{en-gen} shows that
  $E(\Theta_t)\geq E(\Psi)$,
  i.e., $\Theta_t$ is a maximal   entropy interpolant.

 Now, we can prove the following result.
 \begin{theorem} \label{ ent-NP}
  Let 
  $\cZ:=[Z_1~\cdots ~Z_n]\in B(\oplus_{i=1}^n \cY, \cY)$, 
 $B\in B(\cK, \cY)$,  and $C\in B(\cH, \cY)$  be operators
  such  that  $r(\cZ)<1$ and  the grammian  $G_{\{\cZ,B\}}$ is strictly positive.
  Then the central interpolant $\Theta_t$ with tolerance $t> d_\infty(\cZ, B, C)$ 
   is the  maximal   interpolant for the standard left Nevanlinna-Pick problem
    with  operatorial argument and tolerance $t$, and 
    \begin{equation}
    \label{dtt}
    \Delta(\Theta_t)=\frac {1} {t^2} \left[I_\cH+ 
    C^*
 \left( t^2 G_{\{\cZ,B\}}-G_{\{\cZ,C\}}\right)^{-1} C
    \right].
    \end{equation}
 If $\dim \cH<\infty$, then the entropy of
  $\Theta_t$ is given by
 \begin{equation}\label{entr-NP}
 E(\Theta_t)=-\ln \det \left\{
 \frac {1} {t^2} \left[I_\cH+ 
    C^*
 \left( t^2 G_{\{\cZ,B\}}-G_{\{\cZ,C\}}\right)^{-1} C\right]
 \right\}
 \end{equation}
 Moreover,  $\Theta_t$ is the 
  maximal entropy interpolant, i.e., 
  $E(\Psi)\leq E(\Theta_t)$ for any interpolant $\Psi$, 
   and 
  $E(\Theta_t)= E(\Psi)$ if and only if $\Theta_t= \Psi$.
   \end{theorem}
 \begin{proof}
 According to Theorem \ref{maxprin2} (when $\|A\|< t$) and Theorem \ref{NP1}, we have
 $$
 \Delta(\Theta_t)= \Delta(A)= P_\cH D_{A,t}^{-1}|\cH.
 $$
 Using relation \eqref{DA} and the definition of $W_{\{\cZ,C\}}$ and $W_{\{\cZ,B\}}$, 
 we get 
 equation \eqref{dtt}. Now, we can use Theorem \ref{strict} to obtain \eqref{entr-NP}.
 The last part of the theorem follows by applying   Corollary \ref{en-gen}
  to our setting.
 The proof is complete.
 \end{proof}

 \smallskip
Let $\cH$,  $\cK$, and  $\cY_i, \ i=1,\ldots, m$,  be Hilbert spaces 
and consider the operators 
\begin{equation} \label{cond}\begin{split}
&B_j:\cK\to \cY_j,\quad  C_j:\cH\to \cY_j,\quad  j=1,\ldots, m,\\
& \Lambda_j:=[\Lambda_{j1}~\cdots ~\Lambda_{jn}]:\oplus_{i=1}^n \cY_j\to \cY_j,\quad 
 j=1,\ldots, m, 
\end{split}
\end{equation}
such that $r(\Lambda_j)<1$ for any $j=1,\ldots, m$.
The left
 tangential Nevanlinna-Pick  interpolation  problem with operatorial 
argument and tolerance $t>0$ for the  tensor product  $F_n^\infty\bar\otimes B(\cH, \cK)$ is 
to find $\Phi$ in $F_n^\infty\bar{\otimes} B(\cH,\cK)$ such that $\|\Phi\|\leq t$
and 
\begin{equation}\label{interp}
[(I\otimes B_j) \Phi](\Lambda_j)=C_j, \quad  j=1,\ldots, m.
\end{equation}

Define the following operators
$$
B:=\left[\begin{matrix} B_1\\\vdots\\B_m \end{matrix}\right]:
\cK\to \oplus_{j=1}^m \cY_j,\quad 
C:=\left[\begin{matrix} C_1\\\vdots\\C_m \end{matrix}\right]:
\cH\to \oplus_{j=1}^m \cY_j,
$$
and $\cZ:=[Z_1~\cdots~ Z_n]$, where $Z_i$ is the diagonal operator defined by
$$
Z_i:= \left[\begin{matrix}
\Lambda_{1i}\\0 \\ \vdots\\ 0\end{matrix}
\begin{matrix}
0\\ \Lambda_{2i}\\\vdots\\ 0\end{matrix}
\begin{matrix}
0\\0\\\vdots\\ \Lambda_{mi}
\end{matrix}\right]: \oplus_{j=1}^m \cY_j \to \oplus_{j=1}^m \cY_j
$$
for any $i=1,\ldots, n$.
 Note that the interpolation  relation \eqref{interp} is equivalent to
  relation  \eqref{iyc}
 and  we have   
 \begin{equation*} \begin{split}
 t^2G_{\{\cZ,B\}} -G_{\{\cZ,C\}} &=
 \sum_{p=0}^\infty\sum_{|\alpha|=p}Z_\alpha [t^2BB^*-CC^*] Z_\alpha^* \\
 &=
 \left[\sum_{p=0}^\infty\sum_{|\alpha|=p} \Lambda_{j\alpha}[t^2B_j B_k^*- C_j
 C_k^*]\Lambda_{k\alpha}^*\right]_{j,k=1}^m
 \end{split}
 \end{equation*}
where $W_{\{\cZ,B\}}$ and $W_{\{\cZ,C\}}$ are the controllability operators
 associated with
$\{\cZ,B\}$ and $\{\cZ,C\}$, respectively.

It    \cite{Po-nehari}, we proved  that 
 the left tangential Nevanlinna-Pick  interpolation  problem with  data 
$\Lambda_j$, $B_j$, and $C_j$, $j=1,\dots, m$, and tolerance $t>0$, 
 has a solution if and only if
the operator matrix
\begin{equation} \label{np}
\left[\sum_{p=0}^\infty\sum_{|\alpha|=p} \Lambda_{j\alpha}[t^2B_j B_k^*- C_j
 C_k^*]\Lambda_{k\alpha}^*\right]_{j,k=1}^m
\end{equation}
is positive semidefinite.

 \begin{remark}\label{rem}
  All the results of  this
 section can be written   for the
 left
 tangential Nevanlinna-Pick  interpolation  problem with operatorial 
argument   with  data 
$\Lambda_j$, $B_j$, and $C_j$, $j=1,\dots, m$, and tolerance $t>0$.
In particular, one can obtain an explicit form of the maximal entropy 
solution of the above mentioned interpolation 
problem. 
 \end{remark}

 \bigskip

 \subsection{Maximal entropy   interpolation on the unit ball
  of $\CC^n$}
  \label{ball}

 In this section,  we present some consequences of the results of this paper 
  to analytic interpolation on the open  unit ball
 $\BB_n$ of $\CC^n$.
 Let $z_j:=(z_{j1},\ldots, z_{jn})$, \ $j=1,\ldots, m$, be distinct
points in $\BB_n$, and let $B_j\in B(\cK, \cY_j)$, $C_j\in B(\cH, \cY_j)$, \ $j=1,\ldots, m$.
The 
left
 tangential Nevanlinna-Pick  interpolation  problem with
 data $z_j\in \BB_n$, $ B_j, C_j$, $j=1,\ldots, m$, and 
  tolerance $t>0$,  is to find 
 $\Theta\in F_n^\infty\bar\otimes B(\cH, \cK)$ such that
 \begin{equation}\label{BTC}
  B_j \Theta(z_j)= C_j, \quad j=1,\ldots, m,
 \end{equation}
 and $\|\Theta\|\leq t$.
  We proved in \cite{Po-interpo}  that this interpolation 
 problem has solution if and only if the operator matrix
 \begin{equation*}  
 \left[\frac {t^2B_j B_k^*- C_j C_k^*} {1-\left< z_j,
z_k\right>}\right]_{j,k=1}^m
 \end{equation*}
 is positive semidefinite.
 Now, as a consequence of the results of Section \ref{Nev}, we can find 
   the maximal entropy solution of the interpolation 
 problem \eqref{BTC}.  Moreover, under certain natural conditions, we obtain
  an explicit form for the unique solution
  of the Nevanlinna-Pick optimization problem (see Theorem \ref{clas-NP-opt}). 
  Finally, we apply our 
  permanence principle to the Nevanlinna-Pick interpolation problem on the unit ball
  (see Theorem \ref{nperm}).

  First,  we obtain an explicit form  for 
  the maximal entropy solution of the left
 tangential Nevanlinna-Pick  interpolation  problem with tolerance $t$ 
  on the unit ball, 
 under certain  natural conditions.

\begin{theorem}\label{clas-NP}
Let $z_j:=(z_{j1},\ldots, z_{jn})$, \ $j=1,\ldots, m$, be distinct
points in $\BB_n$, and let $B_j\in B(\cK, \cY_j)$, $C_j\in B(\cH, \cY_j)$, \ $j=1,\ldots, m$,
be operators such that the grammian
\begin{equation}\label{gra}
\left[\frac{B_j B_k^*}{1-\left< z_j, z_k\right>}\right]_{j,k=1}^m
\end{equation}
is strictly positive.  Let $t> d_\infty$ and  let $\Theta_t$ be the central 
interpolant  for the   left tangential Nevanlinna-Pick 
interpolation problem with tolerance 
$t$.
Then $\Theta_t(\lambda)= \Phi(\lambda) \Psi(\lambda)^{-1}$ where, for any 
 $\lambda=(\lambda_1,\ldots, \lambda_n)\in \BB_n$,
  
$$
\Phi(\lambda_1,\ldots, \lambda_n)=B^*\left(I-
\sum_{i=1}^n \lambda_i Z_i^* \right)^{-1}
\left(\left[\frac{t^2B_j B_k^*-C_j C_k^*}{1-\left< z_j, z_k\right>}\right]_{j,k=1}^m
\right)^{-1} C
$$ and 
  
$$
\Psi(\lambda_1,\ldots, \lambda_n)= \frac{1} {t^2}\left\{I+C^*
\left(I-
\sum_{i=1}^n \lambda_i Z_i^* \right)^{-1}
\left(\left[\frac{t^2B_j B_k^*-C_j C_k^*}{1-\left< z_j, z_k\right>}\right]_{j,k=1}^m
 \right)^{-1} C
\right\},
$$
where  
\begin{equation}\label{Zi}
Z_i:=\left[ \begin{matrix}
z_{1i} I_{\cY_i}& 0&\cdots & 0\\
0& z_{2i} I_{\cY_i}& \cdots 0\\
0& 0&\cdots & z_{mi} I_{\cY_i}\end{matrix}
\right],
\quad 
B=\left[ \begin{matrix}B_1\\
\vdots\\
B_m
\end{matrix}
\right],
\quad 
C=\left[ \begin{matrix}C_1\\
\vdots\\
C_m
\end{matrix}
\right]
\end{equation}
for any $i=1,\ldots, n$.
In particular,  if \ $\dim \cH<\infty$, then $\Theta_t$ 
 is the maximal entropy solution satisfying  $\|\Theta_t\|<t$  and 
$B_j\Theta_t(z_j)= C_j$ for any $j=1,\ldots, m$.
Moreover, the entropy of $\Theta_t$ is given by
$$
E(\Theta_t)= -\ln \det \left\{
\frac{1} {t^2}\left[ I+C^* 
\left(\left[\frac{t^2B_j B_k^*-C_j C_k^*}{1-\left< z_j, z_k\right>}\right]_{j,k=1}^m
\right)^{-1} C\right]
\right\}.
$$
\end{theorem}

\begin{proof}
Let $z_j:=(z_{j1},\ldots, z_{jn}), \ j=1,\ldots, m$,
 be distinct points in $\BB_n$.   
Note that for any $j,k=1,\ldots, m$, we have 
$$
\sum_{\alpha\in \FF_n^+} z_{j\alpha} \bar z_{k\alpha}=
\frac{ 1}{1-\left< z_j, z_k\right>}.
$$
A simple computation shows that
$$
G_{\{\cZ,B\}}=\left[ \sum_{\alpha\in \FF_n^+} z_{j\alpha}
 \bar z_{k\alpha}B_j B_k\right]_{j, k=1}^m=
 \left[\frac{B_j B_k^*}{1-\left< z_j, z_k\right>}\right]_{j,k=1}^m.
$$
Therefore, we have  
$$
 t^2G_{\{\cZ,B\}} -G_{\{\cZ,C\}}
=\left[\frac {t^2B_j B_k^*- C_j C_k^*} {1-\left< z_j,
z_k\right>}\right]_{j,k=1}^m,
$$
where  $\cZ:= [Z_1~ \cdots~ Z_n]$. 
Since $\text{\rm r}(\cZ)<1$ and the grammian  $G_{\{\cZ,B\}}$ is strictly positive,
we have $d_\infty<\infty$ (see Section \ref{Nev}).  Now,  we  can apply Theorem
\ref{NP1} and find the central interpolant $\Theta_t$ ~for the left tangential
Nevanlinna-Pick interpolation problem  with tolerance $t>d_\infty$, on the unit ball.
 Using Remark \ref{state-sol},
 and taking the compression
of $\Theta_t$ to the symmetric Fock space, we  get the corresponding formulas for
 $\Phi(\lambda_1,\ldots, \lambda_n)$ and $\Psi(\lambda_1,\ldots, \lambda_n)$.
 To  complete the proof of the theorem, we need now to use  Theorem \ref{ ent-NP}
  in our particular setting.
\end{proof}

As mentioned in Section \ref{Nev},
  if  the grammian 
 $G_{\{\cZ,B\}}$ is strictly positive, then   there exists
  $\Theta_{\text{opt}}\in F_n^\infty\bar\otimes B(\cH, \cK)$ solving the problem 
  \begin{equation}\label{-opti}
  [(I\otimes B)\Theta_{\text{opt}}](\cZ)=C  \ \text{ and } \ 
 d_\infty =\|\Theta_{\text{opt}}\|.
 \end{equation}
  In what follows,  we find an explicit form  for  the unique solution 
  of the optimization
  problem \eqref{-opti} on the unit ball $\BB_n$.

\begin{theorem}\label{clas-NP-opt}
Let $z_j:=(z_{j1},\ldots, z_{jn})$, \ $j=1,\ldots, m$, be distinct
points in $\BB_n$, and 
$Z_i$ be defined as in \eqref{Zi}. Let  $B_j\in B(\cK, \cY_j)$ and  $C_j\in B(\CC, \cY_j)$,
 \ $j=1,\ldots, m$,
be operators with $\dim\cY_j<\infty$ and  such that
  the grammian
\begin{equation*} 
\left[\frac{B_j B_k^*}{1-\left< z_j, z_k\right>}\right]_{j,k=1}^m
\end{equation*}
is strictly positive.
Then the unique solution $\Theta_{\text {\rm opt}}\in F_n^\infty\bar\otimes B(\CC, \cK)$ 
for the Nevanlinna-Pick optimization problem 
 satisfies the equation 
\begin{equation}\label{fo}
\Theta_{\text {\rm opt}}(\lambda_1,\ldots, \lambda_n)=\frac 
{d_\infty^2 [B_1^*~\cdots B_n^*] (I-\lambda_1 Z_1^*-\cdots -\lambda_n Z_n^*)^{-1}x}
{[C_1^*~\cdots C_n^*] (I-\lambda_1 Z_1^*-\cdots -\lambda_n Z_n^*)^{-1}x},
\end{equation}
where
 $\lambda=(\lambda_1,\ldots, \lambda_n)\in \BB_n$, $Z_i$ is given by  \eqref{Zi},  and 
$x$ is the eigenvector corresponding to the largest eigenvalue $\lambda_{\text{\rm max}}$
of the operator
$$
\left(\left[\frac{B_j B_k^*}{1-\left< z_j, z_k\right>}\right]_{j,k=1 }^m\right)^{-1}
\left[\frac{C_j C_k^*}{1-\left< z_j, z_k\right>}\right]_{j,k=1}^ m.
$$
Moreover, $d_\infty^2=  \lambda_{\text{\rm max}}$ and 
$\frac {1} {d_\infty}\Theta_{\text {\rm opt}}$ is inner in $F_n^\infty\bar\otimes B(\CC, \cK)$.
In particular,
    $$
    \|\Theta_{\text {\rm opt}}\|= d_\infty,  \quad  
B_j\Theta_{\text {\rm opt}}(z_j)= C_j, \quad \text{ for any  }  \ j=1,\ldots, m.
$$  
\end{theorem}
\begin{proof}
Let $\cZ:=[Z_1~\cdots~ Z_n]\in B(\oplus_{i=1}^n\cY, \cY)$,  ~$ B\in B(\cK, \cY)$,
 and $C\in B(\CC, \cY)$
be the data for the standard left Nevanlinna-Pick interpolation problem, 
where $\dim \cY<\infty$.
First we prove that if $\text{\rm r}(\cZ)<1$ 
and $G_{\{Z,B\}}$ is strictly positive, then 
  there is a  unique solution
  $\Theta\in F_n^\infty\bar \otimes B(\CC, \cK)$ 
 for  the  interpolation problem \eqref{-opt},  given by  the equation
 \begin{equation}
 \label{sol-opt}
 \Theta\left(\sum_{\alpha\in \FF_n^+} e_\alpha \otimes C^* Z_{\tilde \alpha}^* x\right)=
  d_\infty^2 \sum_{\alpha\in \FF_n^+} e_\alpha \otimes B^* Z_{\tilde \alpha}^* x,
 \end{equation}
where $x$ is the eigenvector corresponding to the largest eigenvalue 
$\lambda_{\text{\rm max}}$
of the operator
$G_{\{\cZ,B\}}^{-1} G_{\{\cZ,C\}}\in B(\cY)$. Moreover, 
$d_\infty^2=  \lambda_{\text{\rm max}}$ and 
$\frac {1} {d_\infty}\Theta $ is inner.

Indeed, we know that if   $G_{\{Z,B\}}$ is strictly positive,  then the operator
  $$
  A= W_{\{\cZ,B\}}^*G_{\{\cZ,B\}}^{-1} W_{\{\cZ,C\}}
  $$
   satisfies relation \eqref{A}. 
Let $x$ be an  eigenvector corresponding to the largest eigenvalue 
$\lambda_{\text{\rm max}}$   of $G_{\{\cZ,B\}}^{-1}G_{\{\cZ,C\}}$.
Using the above form of $A$, we get 
\begin{equation*}
\begin{split}
A^* AW_{\{\cZ,C\}}^* &= 
A^* A(A^* W_{\{\cZ,B\}}^*x)= A^* W_{\{\cZ,B\}}^* G_{\{\cZ,B\}}^{-1}G_{\{\cZ,C\}}x\\
&= \lambda_{\text{\rm max}}A^* W_{\{\cZ,B\}}^*x=\lambda_{\text{\rm max}}
W_{\{\cZ,C\}}^*x. 
\end{split}
\end{equation*}
On the other hand, note that $W_{\{\cZ,C\}}^*x\neq 0$ because $W_{\{\cZ,B\}}^*x\neq 0$ is the 
eigenvector corresponding to the largest eigenvalue $\lambda_{\text{\rm max}}$ of $AA^*$.
Since $\cY$ is finite dimensional,  the operator $A$ attains its norm at the vector $g:= W_{\{\cZ,C\}}^*x$.

Now, we use Theorem \ref{opti} to solve the optimization problem 
\eqref{opt2-NP}.
 First, note that  the operator
$G_{\{\cZ,B\}}^*$ is one-to-one and onto $\cH'$ 
(see Section \ref{Nev} for the definition of $\cH'$). Now,  one can easily see that
$A A^*$ is similar to 
$G_{\{\cZ,B\}}^{-1}G_{\{\cZ,C\}}$. In particular, if 
 $\lambda_{\text{\rm max}}$ is the largest eigenvalue of 
 $G_{\{\cZ,B\}}^{-1}G_{\{\cZ,C\}}$, then
 $$
 \lambda_{\text{\rm max}}=\|AA^*\|= \|A\|^2 = d_\infty^2.
 $$
 Note that 
 $$
 W_{\{\cZ,C\}}^*x=\sum_{\alpha\in \FF_n^+} e_\alpha \otimes C^* Z_{\tilde \alpha}^* x
 $$
 and 
 \begin{equation*}
 \begin{split}
 AW_{\{\cZ,C\}}^*x &= W_{\{\cZ,B\}}^*G_{\{\cZ,B\}}^{-1}G_{\{\cZ,C\}}x=
  \lambda_{\text{\rm max}} 
 W_{\{\cZ,B\}}^*x\\
 &= \lambda_{\text{\rm max}} \sum_{\alpha\in \FF_n^+} e_\alpha 
 \otimes B^* Z_{\tilde \alpha}^* x.
 \end{split}
 \end{equation*}
Therefore,  according to the $R_n^\infty$-version of  Theorem \ref{opti},
 the equation \eqref{sol-opt} holds.
Since $\text{\rm r}(\cZ)<1$, as in the proof of Remark \ref{state-sol},  the equation 
\eqref{sol-opt}   becomes 
\begin{equation*}\begin{split}
\Theta[(I\otimes C^*)(I-S_1&\otimes Z_1^*-\cdots -S_n\otimes Z_n^*)^{-1} (1\otimes x)\\
&=
d_\infty
(I\otimes B^*)(I-S_1\otimes Z_1^*-\cdots -S_n\otimes Z_n^*)^{-1} (1\otimes x).
\end{split}
\end{equation*}
Since the symmetric Fock space $F_s^2(H_n)$ is invariant under each operator 
$S_1^*,\ldots, S_n^*$, we can multiply  the previous equation   to the left by the 
orthogonal projection $P_{F_s^2(H_n)\otimes \cK}$ and obtain
\begin{equation*}
\begin{split}
\Theta (\lambda_1,\ldots, \lambda_n) 
C^* (I-\lambda_1 Z_1^*&-\cdots -\lambda_n Z_n^*)^{-1}x\\
&=
d_\infty^2 B^* (I-\lambda_1 Z_1^*-\cdots -\lambda_n Z_n^*)^{-1}x
\end{split}
\end{equation*}
for any $(\lambda_1,\ldots, \lambda_n)\in \BB_n$.
Here, we used the identification of $W_n^\infty$ with the algebra of analytic 
multipliers 
of the symmetric Fock space $F_s^2(H_n)$ (see \cite{Arv1}).
Therefore,   setting 
$\Theta_{\text {\rm opt}}:=\Theta$, the relation \eqref{fo} follows.
 The proof is complete.
\end{proof}

In what follows, we present an application of the permanence principle  
to the Nevanlinna-Pick
 interpolation problem 
on the unit ball $\BB_n$.
Let $z_1, \ldots, z_k$ be distinct points in the unit ball $ \BB_n$, and 
let $C_j\in B (\cK, \cK')$, $ j=1,\ldots, k$.  The Nevanlinna-Pick 
 interpolation  problem with tolerance $t$
for $R_n^\infty\bar \otimes B(\cK, \cK')$ is to find 
$\Phi\in R_n^\infty\bar \otimes B(\cK, \cK')$
such that
\begin{equation}\label{NP-clas}
 \|\Phi\|\leq t\ \text{ and } \ \Phi(z_j)= C_j, \ j=1,\ldots, k.
\end{equation}

 According to \cite{Po-central}, this problem has solutions if and only 
$\|P_{\cN'} \Psi\|\leq t$, where $\Psi\in R_n^\infty \bar \otimes B(\cK, \cK')$ is  an 
arbitrary  operator such that $\Psi(z_j)= C_j$, \ $j=1,\ldots, k$,  and 
\begin{equation}\label{n'}
\cN':= \text{\rm span}~\{ f_{z_j}: \ j=1,\ldots, k\}\otimes \cK',
\end{equation}
where
$$
f_{\lambda}:=(I-\bar \lambda_1S_1-\cdots -\bar \lambda_n S_n)^{-1}(1)\in F^2(H_n) 
$$
 for  any $\lambda=(\lambda_1,\ldots, \lambda_n)\in \BB_n$.  

 In what follows we need the following factorization result
  for contractive multi-analytic operators.
 \begin{lemma}\label{pure}
 Let $\Theta\in  R_n^\infty\bar \otimes B(\cE, \cE')$ be a contractive 
 multi-analytic operator. Then $\Theta$ admits  a unique decomposition 
 $\Theta =\Phi\oplus \Lambda$ with the following properties:
 \begin{enumerate}
 \item[(i)] $\Psi\in R_n^\infty\bar \otimes B(\cE_0, \cE_0')$ is purely contractive, i.e.,
 $\|P_{\cE_0'} \Psi h\|< \|h\|$~ for any $h\in \cE_0$, ~$h\neq 0$;
 \item[(ii)] $\Lambda=I\otimes U\in  R_n^\infty\bar \otimes B(\cE_u, \cE_u')$,
  where $U\in B(\cE_u, \cE_u')$ is a unitary operator;
 \item[(iii)]
 $\cE=\cE_0\oplus \cE_u$ and $\cE'=\cE_0'\oplus \cE_u'$.
 \end{enumerate}
 \end{lemma}
 
 \begin{proof}
 Let $\Theta= \sum\limits_{\alpha\in \FF_n^+} R_\alpha\otimes\theta_\alpha$, \ $\theta_\alpha\in
 B(\cE, \cE')$, be the Fourier representation  of $\Theta$.
 It is well-known that any contraction $\theta_{g_0}\in  B(\cE, \cE')$
  admits a unique decomposition $\theta_{g_0}=Z_0\oplus Z_u$, where $Z_0\in  B(\cE_0, \cE_0')$
  is a pure contraction, i.e., $\|Z_0h\|<\|h\|$ for any $h\in \cE_0$, $h\neq 0$,
  the operator $Z_u\in B(\cE_u, \cE_u')$ is unitary, and we have the orthogonal decompositions
  $\cE= \cE_0\oplus \cE_u$ and  $\cE'= \cE_0'\oplus \cE_u'$.
  Since $Z_u$ is unitary and $\Theta$ is contractice, we deduce
  \begin{equation*}
  \begin{split}
  \|h\|^2&= \|Z_uh\|^2 = \|\theta_{g_0}h\|^2\\
  &\leq 
  \sum_{\alpha\in \FF_n^+}  \|\theta_\alpha h\|^2=\|\Theta h\|^2\leq \|h\|^2
  \end{split}
  \end{equation*} 
 for any $h\in \cE_u$. Therefore, we have equality and 
 $$
 \theta_\alpha |\cE_u=0\quad \text{ for any } \alpha\in \FF_n^+,\  |\alpha|\geq 1.
 $$
 Hence, 
 $$
 \Lambda:= \theta|F^2(H_n)\otimes \cE_u= I\otimes Z_u: F^2(H_n)\otimes \cE_u
 \to F^2(H_n)\otimes \cE_u'
 $$
 is a unitary operator.
 Since $\Theta$ is a  multi-analytic operator  from $F^2(H_n)\otimes \cE$ to $F^2(H_n)\otimes \cE'$, we
  infer that
 $$
 \Theta (F^2(H_n)\otimes \cE_0)\subset F^2(H_n)\otimes \cE_0'.
 $$
 Hence, $\Psi:= \Theta |F^2(H_n)\otimes \cE_0$ is purely contractive and 
  $\Theta =\Phi\oplus \Lambda$. The uniqueness part is straightforward, so we omit it.
 \end{proof}

\begin{theorem}\label{nperm}
Let $z_1, \ldots, z_k$ be distinct points in the unit ball $ \BB_n$, and 
let $C_j\in B(\cK, \cK')$, \ $ j=1,\ldots, k$. 
Let 
$$
A:= P_{\cN'}R\in B(F^2(H_n)\otimes \cK, \cN'),
$$
 where 
where $R\in R_n^\infty \bar \otimes B(\cK, \cK')$ is  an 
arbitrary  operator such that $R(z_j)= C_j$, \ $j=1,\ldots, k$
and the subspace $\cN'$ is defined by  relation \eqref{n'}.
If $\dim \cK<\infty$, then the central intertwining lifting 
$\Psi_{\text{\rm max}}\in  R_n^\infty\bar \otimes B(\cK, \cK')$ of $A$
is  the 
 maximal entropy solution  for the Nevanlinna-Pick  interpolation problem.
 
 Let $m>k$ and $z_{k+1},\ldots, z_m\in \BB_n$ be
 such that  $z_1, \ldots, z_k,z_{k+1},\ldots, z_m $  are distinct points, and let
 $$
  C_j:= \Psi_{\text{\rm max}}(z_j), \quad  j=k+1,\ldots, m.
  $$
 Then $\Psi_{\text{\rm max}}$ is also the maximal entropy solution for  the Nevanlinna-Pick
 inerpolation  problem  with the new data $\{z_j\}_{j=1}^m$ and  $\{C_j\}_{j=1}^m$.
\end{theorem}

\begin{proof}
Since $\cN'$ is invariant under each operator $S_i^*\otimes I_{\cK'}$,\ $i=1,\ldots, n$,  we have 
\begin{equation}\label{n}
[F^2(H_n)\otimes \cK']\ominus \cN'= \Phi (F^2(H_n)\otimes \cK_1),
\end{equation}
where $\Phi\in R_n^\infty\bar\otimes B(\cK_1 , \cK')$ is an inner operator.
Moreover,  let us prove that  $\Phi$ is a pure inner operator.
According to Lemma \ref{pure}, we have the decomposition $\Phi=\chi_1\oplus \chi_2$
with the following properties:
 \begin{enumerate}
 \item[(i)] $\chi_1\in R_n^\infty\bar \otimes B(\cE_1, \cE_1')$ is purely contractive;
 \item[(ii)] $\chi_2=I\otimes U$, where $U\in B(\cE_2, \cE_2')$ is a unitary operator;
 \item[(iii)]
 $\cK_1=\cE_1\oplus \cE_2$ and $\cK'=\cE_1'\oplus \cE_2'$.
 \end{enumerate}
This implies that $F^2(H_n)\otimes \cE_2'$ is in the range of $\Phi$.
According to relation \eqref{n}, we have 
$$
F^2(H_n)\otimes \cE_2'\subset F^2(H_n)\otimes \cK']\ominus \cN'
$$ 
and therefore
 $F^2(H_n)\otimes \cE_2'\perp \cN'$.
 Hence, for any $\alpha\in \FF_n^+$, $k'\in \cK'$, and $h\in \cE_2'$, we have 
  $e_\alpha\otimes h\perp f_{z_j} \otimes k'$. Hence, and taking into
   account the definition of 
 $ f_{z_j}$, we get 
 \begin{equation}\label{0}
 z_{i\alpha} \left< h, k'\right>=0
 \end{equation}
 for any  $k'\in \cK'$. 
 Since  $z_1, \ldots, z_k$ are distinct points in the unit ball $ \BB_n$,   we can find 
 $z_{i\alpha}\neq 0$. Therefore, relation \eqref{0} implies $h=0$.
 This proves that $\cE_2'=\{0\}$, which shows that $\Phi$ is a pure inner operator.

Note that if  $\Lambda\in R_n^\infty\bar\otimes B(\cK, \cK')$, 
then $P_{\cN'} \Lambda =0$ if
 and only if $\Lambda (z_j)=0$ for any $j=1,\ldots, m$.
Now, it is clear that $\Psi\in R_n^\infty\bar\otimes B(\cK, \cK')$ is a solution for 
the Nevanlinna-Pick interpolation problem 
with tolerance $t$ if and only if $\Psi$ is a solution  for the problem 
\eqref{sa1}, i.e.,
$$
\Psi = R+\Theta G \ \text{  and } \ \|\Psi\|\leq t,
$$
 where $\Phi$ is given by \eqref{n}, and 
$R\in R_n^\infty\bar\otimes B(\cK, \cK')$ is  such that $R(z_j)=C_j$, $j=1,\ldots, m$
(for the existence of such operator see \cite{ArPo2}).

Using Theorem \ref{sara1}, we see that there is a solution  for 
the Nevanlinna-Pick interpolation problem  if and only if the operator
 $A:= P_{\cN'} R\in B(F^2(H_n)\otimes \cK, \cN')$
 satisfies $\|A\|\leq t$. Moreover, the central intertwining lifting of $A$ is
  the maximal entropy solution  $\Psi_{\text{\rm max}}$ 
 for 
the Nevanlinna-Pick interpolation problem.

Since the subspace 
\begin{equation}\label{n''}
\cN'':= \text{\rm span}~\{ f_{z_j}: \ j=1,\ldots, m\}\otimes \cK'
\end{equation}
is invariant under each operator $S_i^*\otimes I_{\cK'}$,\ $i=1,\ldots, n$,  we have 
\begin{equation*} 
[F^2(H_n)\otimes \cK']\ominus \cN''= \Psi (F^2(H_n)\otimes \cK_2),
\end{equation*}
where $\Psi\in R_n^\infty\bar \otimes B(\cK_2 , \cK')$ is an inner operator.
Since $\cN'\subset \cN''$,  we have
$$
\Psi (F^2(H_n)\otimes \cK_2)\subset \Phi (F^2(H_n)\otimes \cK_1).
$$
Hence, $\Psi = \Phi \Phi_1$ for some inner operator
 $\Psi_1\in R_n^\infty\bar\otimes B(\cK_2 , \cK_1)$
The Nevanlinna-Pick interpolation with data $\{z_j\}_{j=1}^m$ and $\{C_j\}_{j=1}^m$
is equivalent to the interpolation problem 
\eqref{sa2}, i.e.,
$$
\|\Gamma\|\leq t \quad \text{ and }\quad \Gamma= \Psi_{\text{\rm max}}+ \Phi \Phi_1 G.
$$

Now, we can use the permanence principle of Theorem \ref{sara-p} to prove
 the last part of the theorem.
\end{proof}



\clearpage


\begin{thebibliography}{99}


  
\bibitem{AMc2} {\sc J.~Agler and J.~E.~McCarthy}, 
 Complete Nevanlina-Pick kernels, 
%
{\it J. Funct. Anal.}
  {\bf 175} (2000), 111--124. 
 
     
\bibitem{ArPo} {\sc  A.~Arias and G.~Popescu},  
Factorization and reflexivity on Fock spaces, {\it  Integr. Equat. Oper.Th.}
 %
 {\bf  23} (1995),  268--286. 
     
  
  
     
\bibitem{ArPo1} {\sc  A.~Arias and G.~Popescu},  Noncommutative interpolation
 and Poisson transforms II, {\it Houston J. Math.}
 {\bf  25} (1999),  79--98. 

%

\bibitem{ArPo2} {\sc  A.~Arias and G.~Popescu},  Noncommutative interpolation
 and Poisson transforms, {\it Israel J. Math.}
 %
 {\bf  115} (2000),  205--234. 
 
 
 \bibitem{Aro} {\sc A.~Z.~Arov and M.~G.~Krein},
 On computations of entropy functionals and their minima,
 {\it Acta. Sci. Math. (Szeged)}
 {\bf 45} (1983), 51--66.

 
  
     
\bibitem{Arv1} {\sc W.B.~Arveson}, 
 Subalgebras of $C^*$-algebras III: Multivariable operator theory, 
 {\it Acta Math.}
 {\bf 181} (1998),  159--228.   
 
 
    
\bibitem{Arv2} {\sc W.B.~Arveson}, 
  The curvature invariant of a Hilbert module  over $\CC [z_1,\ldots, z_n]$, 
 {\it J. Reine Angew. Math.}
 {\bf 522} (2000),  173--236.

\bibitem{BTV} {\sc J.~A.~Ball, T.~T.~Trent, and V.~Vinikov}, Interpolation and commutant
lifting for multipliers on reproducing kernels Hilbert spaces,
{\it Operator Theory and Analysis: The M.A. Kaashoek Anniversary Volume}, 
pages 89--138,
{\bf OT 122}, Birkhauser-Verlag, Basel-Boston-Berlin, 2001.
     
 
 \bibitem{BB} {\sc J.~A.~Ball and V.~Bolotnikov}, On bitangential interpolation 
 problem for contractive-valued
 functions on the unit ball,
 {\it  Linear Algebra Appl.}  
  {\bf 353} (2002), 107--147. 
 
\bibitem{B} {\sc J.~W.~Bunce},
 Models for n-tuples of noncommuting operators, {\it J. Funct. Anal.}
  {\bf 57}(1984), 21--30. 
 
 
 \bibitem{Bu} {\sc J. Burg},
 {\em Maximum entropy spectral analysis}, 
 Ph.D. dissertation, Stanford University, Stanford, 1975.
     
\bibitem{Ca} {\sc C.~Carath\' eodory},
 \" Uber den Variabilit\" atsbereich der Koeffzienten von
  Potenzreinen die gegebene Werte nicht annehmen,
{\it Math. Ann.}
 {\bf 64}  (1907), 95--115.
 

      
\bibitem{DKP}  {\sc K.~R.~Davidson, E.~Katsoulis, and D.~Pitts},      
  The structure of free semigroup algebras,
 {\it J. Reine Angew. Math.} 
  {\bf 533} (2001), 99--125.


\bibitem{DP1}  {\sc K.~R.~Davidson and D.~Pitts}, 
 Invariant subspaces and hyper-reflexivity for free semigroup algebras,
{\it  Proc. London Math. Soc.}
 {\bf 78} (1999),  401--430. 
 
 
 
      
 
 
\bibitem{DP2} {\sc K.~R.~Davidson and D.~Pitts},
%
 Automorphisms and representations of the noncommutative
  analytic Toeplitz algebras,
{\it  Math. Ann.}
   {\bf 311} (1998),  275--303. 
 

 
\bibitem{DP} {\sc K.~R.~Davidson and D.~Pitts},     
 Nevanlinna-Pick interpolation for noncommutative analytic Toeplitz algebras,
 {\it  Integr. Equat. Oper.Th.}
{\bf 31} (1998), 321--337.   





 \bibitem{DG1} {\sc H.~Dym and I.~ Gohberg},
  A maximum entropy principle for
  contractive interpolants,
  {\it J. Funct. Anal.} {\bf 65} (1986), 83--125.

\bibitem{DG2} {\sc H.~Dym and I.~ Gohberg},
  A new class of contractive interpolants
  and maximum entropy principle,
   %
  {\it Topics in operator theory and interpolation}, 117--150,
  Oper. Theory Adv. Appl., 29, Birkhäuser, Basel, 1988.

\bibitem{EGL} {\sc R.~L.~Elis, I.~Gohberg, and D.~Lay},
 Band extensions maximum entropy and the permanence principle,
 {\it Maximum Entropy and Bayesian Methods in Applied Statistics},
 Cambridge University Press, Cambridge, 1986.




\bibitem{EP} {\sc J.~Eschmeier and M.~Putinar}, Spherical contractions and interpolation problems
on the unit ball,
{\it J. Reine Angew. Math.}
{\bf 542} (2002), 219--236.
  
\bibitem{FF-book}  {\sc C.~Foias, A.~E.~Frazho},
 {\it The commutant lifting approach to interpolation problems},
 Operator Theory: Advances and Applications,
Birh\"auser Verlag,  Bassel,
 1990.
 
 



\bibitem{FFG}  {\sc C.~Foias, A.~E.~Frazho, and I.~Gohberg}, 
Central intertwining lifting, maximum entropy and their permanence, 
{\it  Integr. Equat. Oper. Th.}
  {\bf 18} (1994),  166--201.

 %

 
\bibitem{FFGK-book}  {\sc C.~Foias, A.~E.~Frazho, I.~Gohberg, and M.~Kaashoek},
{\em  Metric constrained interpolation, commutant lifting and systems}, 
  Operator Theory: Advances and Applications vol. 100,  Birh\"auser Verlag,
 Bassel, 1998.

 \bibitem{F} {\sc  A.~E.~Frazho},
 Models for noncommuting operators, {\it J. Funct. Anal.}
  {\bf 48}(1982), 1--11. 


\bibitem{Fra} {\sc B.~A.~Francis},
{\em A Course in $H^\infty$ Control Theory},
Lecture Notes in Control and Information Sciences, Springer-Verlag, New York, 1987.

\bibitem{GM1} {\sc K.~Glover and D.~Mustafa},
Derivation of the maximum entropy
 $H_\infty$-controller and a state-space formula for its entropy,
 {\it Int. J. Cont.}
 {\bf 50} (1989), 899--916.
 
 \bibitem{GM2} {\sc K.~Glover and D.~Mustafa},
 {\em Minimum Entropy  $H_\infty$-control},
 Lecture Notes in Control and Information Sciences, 
 Springer-Verlag, New York, 1990.
 
 
 
 

\bibitem{GSz} {\sc U.~Grenander and G.~Szeg\" o},   
 {\em  Toeplitz forms and their applications}, 
 Univ. of California Press,
 Berkeley, California, 1958. 
 

\bibitem{HL1}  {\sc H.~Helson and D.~Lowdenslager},  
Prediction theory and Fourier series in several variables,
   {\it  Acta Math.} 
{\bf 99} (1958),  165--202.
%
 

\bibitem{HL2} { \sc H.~Helson and D.~Lowdenslager}, 
 Prediction theory and Fourier series in several variables II,
{\it Acta Math.} 
 {\bf 106} (1961), 175--213.   
 
 
\bibitem{H} {\sc K.~Hoffman},  {\em Banach Spaces of Analytic Functions},
Englewood Cliffs: Prentice-Hall, 1962.

\bibitem{K1}  {\sc A.~N.~Kolmogorov},   
 Sur l'interpolation et l'extrapolation des suites stationaire,
{\it C. R. Acad. Sci. (Paris)}
 {\bf 208} (1939), 2043--2045.


\bibitem{K2}  {\sc A.~N.~Kolmogorov },   
  Stationary sequences in Hilbert spaces,
{\it Bull. Math. Univ. Moscow}
  {\bf 2} (1941).
 


%


\bibitem{Ne} {\sc Z.~Nehari},
On bounded bilinear forms,
{\it Ann. of Math.}
{\bf 65} (1957), 153--162.


 
\bibitem{N}  {\sc R.~Nevanlinna}, 
 \"Uber  beschr\"ankte Functionen, die in gegebenen Punkten 
vorgeschribene Werte annehmen,
%
{\it Ann. Acad. Sci. Fenn. Ser A}
  {\bf 13} (1919),  7--23.
 

   
\bibitem{PPoS} {\sc V.~I.~Paulsen, G.~Popescu,  and D.~Singh},   
  On Bohr's inequality,
 {\it Proc. London Math. Soc.} 
 {\bf 85} (2002), 493--512.
 
\bibitem{Po-models} {\sc G.~Popescu}, Models for infinite sequences of
noncommuting operators, {\it Acta. Sci. Math. (Szeged)} {\bf
53}(1989), 355--368.


\bibitem{Po-isometric} {\sc G.~Popescu}, Isometric dilations for infinite
sequences of noncommuting operators, {\it Trans. Amer. Math. Soc.}
{\bf 316}(1989), 523--536.

\bibitem{Po-charact} {\sc G.~Popescu}, Characteristic functions for infinite
sequences of noncommuting operators, {\it J. Operator Theory}
{\bf 22}(1989), 51--71.

 
\bibitem{Po-multi} {\sc G.~Popescu},
 Multi-analytic operators and some  factorization theorems,
 {\it Indiana Univ. Math.~J.} 
 {\bf 38} (1989),   693--710.
 
  
\bibitem{Po-intert} {\sc G.~Popescu},  On intertwining dilations for
 sequences of noncommuting operators, {\it J. Math. Anal. Appl.}
{\bf 167} (1992), 382--402.

 
\bibitem{Po-von} {\sc G.~Popescu},
 Von Neumann inequality for $(B(H)^n)_1$, {\it  Math. Scand.} 
{\bf 68} (1991), 292--304.

\bibitem{Po-funct} {\sc G.~Popescu},  
 Functional calculus for noncommuting operators, {\it  Michigan Math. J.}
 {\bf 42}  (1995),    345--356.  
  
\bibitem{Po-analytic} {\sc G.~Popescu}, 
 Multi-analytic operators on Fock spaces, {\it  Math. Ann.}
{\bf 303} (1995),  31--46.

\bibitem{Po-disc} {\sc G.~Popescu}, Noncommutative disc algebras and their
  representations, {\it Proc. Amer. Math.} {\bf 124} (1996), 2137--2148. 
  
\bibitem{Po-interpo} {\sc G.~Popescu},  
  Interpolation problems in several variables, {\it  J. Math Anal. Appl.}
 {\bf 227}  (1998),    227--250.  

      
\bibitem{Po-poisson} {\sc  G.~Popescu}, 
 Poisson transforms on some $C^*$-algebras generated by isometries,
 {\it J. Funct. Anal.}
  {\bf 161} (1999),  27--61.  
 


\bibitem{Po-structure} {\sc G.~Popescu},  
Structure and entropy for positive definite Toeplitz kernels on free semigroups,
   {\it  J. Math. Anal. Appl.} {\bf 254} (2001), 191-218. 
 
\bibitem{Po-spectral}  {\sc G.~Popescu},  
  Spectral liftings in Banach algebras 
and interpolation in several variables,
{\it  Trans. Amer. Math. Soc.}
%
  {\bf 353} (2001),   2843--2857. 
 

  
 \bibitem{Po-tensor} {\sc  G.~Popescu},
 Commutant lifting, tensor algebras, and functional calculus,
 {\it Proc.  
      Edinb. Math. Soc.} {\bf 44} (2001), 389-406.
      
\bibitem{Po-curvature} {\sc  G.~Popescu}, 
  Curvature invariant for Hilbert modules over free semigroup algebras,
   {\it Adv. Math.}
 {\bf 158} (2001), 264--309. 
 
%
\bibitem{Po-lifting} {\sc G.~Popescu},  
A lifting theorem for symmetric commutants,
{\it Proc. Amer. Math.}
{\bf 129} (2001), 1705--1711.
  
\bibitem{Po-central} {\sc G.~Popescu},  
  Central intertwining lifting, suboptimization, and interpolation
  in several variables, {\it  J. Funct. Anal.}  {\bf 189}  (2002), 132--154.
   

        
\bibitem{Po-meromorphic} {\sc G.~Popescu},  
   Meromorphic interpolation in several variables,
   {\it  Linear Algebra Appl.}  
  {\bf 357} (2002), 173--196. 
  
  
       
\bibitem{Po-nehari} {\sc G.~Popescu},  
  Multivariable Nehari problem and interpolation,
   {\it  J. Funct. Anal.} 
   {\bf 200} (2003), 536--581.
   
   \bibitem{Po-similarity} {\sc G.~Popescu},  
    Similarity and ergodic theory of positive linear maps,
   {\it J. Reine Angew. Math.} 
  {\bf 561} (2003), 87--129.
  
  \bibitem{R}  {\sc E.~A.~Robinson},
  {\em Random Wavelets and Cybernetic Systems}, Griffin, London, 1962. 
    
   
  \bibitem{RSW} {\sc L.~Rodman , I.~M.~Spitkovsky, and H.~J.~Woerdeman},
  Abstract band
  method via factorization, positive and band extensions of multivariable almost periodic
  matrix functions, and spectral estimation,
  {\it  Mem. Amer. Math. Soc} {\bf 160} (2002), no. 762.
   
   
\bibitem{RR}  {\sc M.~Rosenblum and  J.~Rovnyak},    
{\em Hardy classes and operator theory},
 Oxford University Press-New York,  1985.  
 
 
\bibitem{S} {\sc D.~Sarason}, 
 Generalized interpolation in $H^\infty$,
{\it Trans. Amer. Math. Soc.} 
  {\bf 127} (1967),   179--203.  
   
  
 
\bibitem{Sc}  {\sc I.~Schur},  
%
 \"Uber Potenzreihen die im innern des Einheitshreises beschr\"ankt sind,
{\it  J. Reine Angew. Math.}
 {\bf 148} (1918), 122--145.
 

 
\bibitem{Sz1} {\sc G.~Szeg\" o},   
%
  Ein Grenzwertsvatz uber die Toeplitzschen Determinanten, 
  einer reellen positiven Funktion,  
{\it  Math. Ann.}
 {\bf  76} (1915), 490--503.
 
 
\bibitem{Sz2} {\sc G.~Szeg\" o},   
%
  Orthogonal polynomials,  
{\it Colloquium Publications }
{\bf 23} (1939), 
 Amer. Math. Soc. Providence, Rhode Island. 
  

 
\bibitem{Sz3}  {\sc G.~Szeg\" o},  
%
  {\em On certain Hermitian forms associated with the Fourier series
 of a positive function}, pp. 222--238, 
 Festskrift Marcel Riesz,  Lund, 1952.
  
\bibitem{SzF} {\sc B.~Sz.-Nagy and C.~Foia\c{s}},
Dilation des commutants d'op\' erateurs,
{\it C. R. Acad. Sci. Paris, s\' erie A}
{\bf 266} (1968), 493--495.

\bibitem{SzF-book} {\sc B.~Sz.-Nagy and C.~Foia\c{s}}, {\em Harmonic
Analysis of Operators on Hilbert Space}, North Holland, New York 1970.



\bibitem{W}  {\sc N.~Wiener},  
{\em  Extrapolation, Interpolation, and Smoothing of stationary time series},
 Wiley, New York, 1949.
 

 \bibitem{WM1} {\sc N.~Wiener and P.~R.~Masani},  111--150
The prediction theory of multivariate stochastic processes I,
{\it  Acta Math.}
  {\bf 98} (1957),  111--150.
 

 \bibitem{WM2} {\sc N.~Wiener and P.~R.~Masani},   
 The prediction theory of multivariate stochastic processes II,
{\it  Acta Math.}
  {\bf  99} (1958),  93--139.
 






\end{thebibliography}
\end{document}